\documentclass[9pt]{article}
\usepackage{mathrsfs}
\usepackage{amsthm}
\usepackage{amssymb}
\usepackage{amsmath}
\usepackage{graphicx}
\usepackage{color}
\usepackage{amsfonts}
\usepackage{float}
\usepackage{cite}
\usepackage [latin1]{inputenc}
\usepackage[text={140mm,210mm},left=35mm,vmarginratio=1:1]{geometry}
\newtheorem{theorem}{Theorem}[section]

\newtheorem{lemma}[theorem]{Lemma}
\newtheorem{proposition}[theorem]{Proposition}

\numberwithin{equation}{section}
\normalsize

\begin{document}
\title{\textbf{Stationary fluctuations for occupation times of the long-range voter models on lattices}}

\author{Xiaofeng Xue \thanks{\textbf{E-mail}: xfxue@bjtu.edu.cn \textbf{Address}: School of Mathematics and Statistics, Beijing Jiaotong University, Beijing 100044, China.}\\ Beijing Jiaotong University}

\date{}
\maketitle

\noindent {\bf Abstract:} In this paper, we are concerned with the long-range voter model on lattices. We prove a stationary fluctuation theorem for the occupation time of the model under a proper time-space scaling. In several cases, the fluctuation limits are driven by fractional Brownian motions with Hurst parameters in $(1/2, 1)$. The proof of our main result utilizes the martingale decomposition strategy introduced in \cite{Kipnis1987}. A local central limit theorem of the long-range random walk, the duality relationship between the model and the long-range coalescing random walk and a fluctuation theorem of the empirical density field of the model play the key roles in the proof.

\quad



\section{Introduction}\label{section one}
\subsection{The model}\label{subsection 1.1}
In this paper, we are concerned with the voter model on a lattice. We aim to extend the central limit theorem of the occupation time of the nearest-neighbor voter model given in \cite{Cox1983} to the long-range case. We first recall the definition of the long-range voter model. For given integer $d\geq 1$ and constant $\alpha>0$, the long-range voter model $\{\eta_t\}_{t\geq 0}$ on the lattice $\mathbb{Z}^d$ with parameter $\alpha$ is a continuous-time Markov process with state space $\{0, 1\}^{\mathbb{Z}^d}$ and generator $\mathcal{G}$ given by
\begin{equation}\label{equ 1.1 generator of long-range voter model}
\mathcal{G}f(\eta)=\sum_{x\in \mathbb{Z}^d}\sum_{y:y\neq x}\|y-x\|_2^{-(d+\alpha)}\left(f(\eta^{x,y})-f(\eta)\right)
\end{equation}
for any $\eta\in \{0, 1\}^{\mathbb{Z}^d}$ and local $f:\{0, 1\}^{\mathbb{Z}^d}\rightarrow \mathbb{R}$, where
$\eta^{x, y}$ is defined as
\[
\eta^{x, y}(z)=
\begin{cases}
\eta(z) & \text{~if\quad} z\neq x,\\
\eta(y) & \text{~if\quad} z=x
\end{cases}
\]
and $\|\cdot\|_2$ is the $l_2$-norm on $\mathbb{Z}^d$, i.e.,
\[
\|x\|_2=\sqrt{\sum_{j=1}^dx_j^2}
\]
for any $x=(x_1, \ldots, x_d)\in \mathbb{Z}^d$. Note that $\eta^{x, y}\neq \eta^{y, x}$, i.e., we care about the order of $x$ and $y$.

The voter model describes the variations of the respective supporters of two opposite opinions of a topic. In detail,  $\{x: \eta(x)=0\}$ are individuals with opinion `$0$' and $\{x: \eta(x)=1\}$ are individuals with opinion `$1$'. For any $x\in \mathbb{Z}^d$ and $y\neq x$, $x$ adopts the opinion of $y$ at rate $\|y-x\|_2^{-(d+\alpha)}$. Note that, denoting by $\mathbf{0}$ the origin of $\mathbb{Z}^d$, we have
\[
\sum_{y\neq \mathbf{0}}\|y\|_2^{-(d+\alpha)}<+\infty,
\]
which ensures that the process with generator given by \eqref{equ 1.1 generator of long-range voter model} is well-defined. Readers could see Chapter 5 of \cite{Lig1985} and Part {\rm \uppercase\expandafter{\romannumeral2}} of \cite{Lig1999} for a comprehensive survey of basic properties of the voter model.

In this paper, we are concerned with the occupation time of the long-range voter model. For any $t\geq 0$, the occupation time on the origin $\mathbf{0}$ during $[0, t]$ is defined as
\[
\xi_t=\int_0^t\eta_u(O)du.
\]
Since 1980s, the limit theorem of the occupation time of the voter model is a popular research topic. For the nearest-neighbor voter model on $\mathbb{Z}^d$, it is proved in \cite{Cox1983} that, when $d\geq 2$ and the process starts from the product measure with density $p$, $\frac{1}{A_d(t)}\left(\xi_t-tp\right)$ converges weakly to a Gaussian distribution as $t\rightarrow+\infty$, where
\[
A_d(t)=
\begin{cases}
\frac{t}{\sqrt{\log t}} & \text{~if\quad} d=2,\\
t^{\frac{3}{4}} & \text{~if\quad} d=3,\\
\sqrt{t\log t} & \text{~if\quad} d=4,\\
\sqrt{t} & \text{~if\quad} d\geq 5.
\end{cases}
\]
Reference \cite{Xue2024} extends the above result to a functional central limit theorem when $d\geq 3$. In detail, for the process in the nearest-neighbor case starting from the product measure with density $p$, it is shown in \cite{Xue2024} that, for any given $T>0$,
\[
\left\{\frac{1}{A_d(N)}(\xi_{tN}-tNp):~0\leq t\leq T\right\}
\]
converges weakly, under the uniform topology, to $\{\beta_t^d\}_{0\leq t\leq T}$ as $N\rightarrow+\infty$, where $\{\beta_t^d\}_{t\geq 0}$ is a Brownian motion when $d\geq 4$ or a Gaussian process with continuous sample path, mean zero and covariance functions given by
\[
{\rm Cov}\left(\beta_s^3, \beta_t^3\right)=s^{3/2}+t^{3/2}-\frac{1}{2}|t-s|^{3/2}-\frac{1}{2}(t+s)^{3/2}
\]
when $d=3$. Large deviation principles are also discussed for the occupation time of the nearest-neighbor voter model on $\mathbb{Z}^d$ when $d\geq 2$, see
\cite{Bramson1988} and \cite{Maillard2009}. In this paper, we extend the main results in \cite{Cox1983} and \cite{Xue2024} to the long-range voter model starting from a stationary distribution. For mathematical details, see Section \ref{section two}.

\subsection{The long-range random walk on $\mathbb{Z}^d$}\label{subsection 1.2}
The duality relationship between the voter model and the coalescing random walk plays the key role in the research of the limit behaviors of the voter models. We will recall this duality relationship in the next subsection. In this subsection, we recall basic properties of the long-range random walk as a preliminary.

For given integer $d\geq 1$ and constant $\alpha>0$, we denote by $\{X_t^{d, \alpha}\}_{t\geq 0}$ the long-range random walk on $\mathbb{Z}^d$ with generator $\mathcal{L}^{d, \alpha}$ given by
\begin{equation}\label{equ 1.2 generator longRange RW}
\mathcal{L}^{d, \alpha}h(x)=\sum_{y\in \mathbb{Z}^d\setminus\{\mathbf{0}\}}\|y\|_2^{-(d+\alpha)}\left(h(x+y)-h(x)\right)
\end{equation}
for any $x\in \mathbb{Z}^d$ and $h: \mathbb{R}^d\rightarrow \mathbb{R}$. That is to say, conditioned on $X_\cdot^{d, \alpha}=x$, the walk jumps to $y$ at rate $\|y-x\|_2^{-(d+\alpha)}$. We denote by $\{p_t^{d, \alpha}\left(\cdot, \cdot\right)\}_{t\geq 0}$ the transition probabilities of $\{X_t^{d, \alpha}\}_{t\geq 0}$, i.e.,
\[
p_t^{d, \alpha}(x, y)=\mathbb{P}\left(X_t^{d, \alpha}=y\Big|X_0^{d, \alpha}=x\right)
\]
for any $x, y\in \mathbb{Z}^d$.

We first recall the scaling limit of $\{X_t^{d, \alpha}\}_{t\geq 0}$. For $t>0$ and $\alpha>0$, we define
\[
h_\alpha(t)=
\begin{cases}
\sqrt{t} & \text{~if\quad} \alpha>2,\\
\sqrt{t \log t} & \text{~if\quad} \alpha=2,\\
t^{\frac{1}{\alpha}} & \text{~if\quad} 0<\alpha<2.
\end{cases}
\]
Then for given $T>0$, condition on $X_0^{d, \alpha}=\mathbf{0}$,
\[
\left\{\frac{1}{h_\alpha(N)}X_{tN}^{d, \alpha}:~0\leq t\leq T\right\}
\]
converges weakly, under the Skorohod topology, to a $\mathbb{R}^d$-valued stable process $\{\hat{X}_t^{d, \alpha}\}_{0\leq t\leq T}$ starting from $\mathbf{0}$ as $N\rightarrow+\infty$. The generator $\hat{\mathcal{L}}^{d, \alpha}$ of $\{\hat{X}_t^{d, \alpha}\}_{t\geq 0}$ is given as follows. When $0<\alpha<2$,
\[
\hat{\mathcal{L}}^{d, \alpha}f(u)=\int_{\mathbb{R}^d}\frac{f(v+u)-f(u)}{\|v\|^{d+\alpha}}dv
\]
for any $u\in \mathbb{R}^d$ and $f\in C_c^2(\mathbb{R}^d)$. Note that $0<\alpha<2$ ensures that $\int_{\mathbb{R}^d}\frac{|f(v+u)-f(u)|}{\|v\|^{d+\alpha}}dv<+\infty$. When $\alpha>2$,
\[
\hat{\mathcal{L}}^{d, \alpha}f(u)=\frac{1}{2}\sum_{i=1}^d\sum_{j=1}^d\left(\sum_{y\in \mathbb{Z}^d\setminus \mathbf{\{0\}}}\frac{y_iy_j}{\|y\|_2^{d+\alpha}}\right)\partial_{u_iu_j}^2f(u)
\]
for any $u\in \mathbb{R}$ and $f\in C_c^2(\mathbb{R}^d)$, where $u_i$ is the $i$th coordinate of $u$. Note that $\alpha>2$ ensures that $\sum_{y\in \mathbb{Z}^d\setminus \mathbf{0}}\frac{|y_iy_j|}{\|y\|_2^{d+\alpha}}<+\infty$. When $\alpha=2$,
\[
\hat{\mathcal{L}}^{d, 2}f(u)=\frac{1}{2}\sum_{i=1}^d\sum_{j=1}^dK_{d, 2}^{i, j}\partial_{u_iu_j}^2f(u),
\]
where
\[
K_{d, 2}^{i, j}=\frac{1}{2}\int_{\{v\in \mathbb{R}^d:~\|v\|_2=1\}}v_iv_jdS
\]
and $dS$ represents the surface integral. Consequently, when $\alpha\geq 2$, $\{\hat{X}_t^{d, \alpha}\}_{t\geq 0}$ equals a linear transformation on the $d$-dimension standard Brownian motion.

The conclusion that $\frac{1}{h_\alpha(N)}X_{tN}^{d, \alpha}$ converges weakly to $\hat{X}_t^{d, \alpha}$ is well-known, so here we only give an outline of how to check this weak convergence for readers not familiar with this conclusion. Let $\mathcal{L}_N^{d, \alpha}$ be the generator of $\{\frac{1}{h_\alpha(N)}X_{tN}^{d, \alpha}\}_{t\geq 0}$, then we only need to check that $\mathcal{L}_N^{d, \alpha}f\rightarrow \hat{\mathcal{L}}^{d, \alpha}f$ as $N\rightarrow+\infty$. By \eqref{equ 1.2 generator longRange RW},
\[
\mathcal{L}_N^{d, \alpha}f(u)=N\sum_{y:y\in \mathbb{Z}^d\setminus\{\mathbf{0}\}}\|y\|_2^{-(d+\alpha)}\left(f\left(u+\frac{1}{h_\alpha(N)}y\right)-f(u)\right).
\]

When $\alpha<2$, $\mathcal{L}_N^{d, \alpha}f\rightarrow \hat{\mathcal{L}}^{d, \alpha}f$ follows from the definition of Riemann integral. When $\alpha>2$, for sufficiently small $\epsilon$ and $y$ such that $\|y\|_2\leq \sqrt{N}\epsilon$,
\[
f\left(u+\frac{1}{\sqrt{N}}y\right)-f(u)\approx \sum_{i=1}^d\partial_{u_i}f(u)\frac{y_i}{\sqrt{N}}+\frac{1}{2N}\sum_{i=1}^d\sum_{j=1}^dy_iy_j\partial^2_{u_iu_j} f(u)
\]
by Taylor's expansion up to the second order. Since $(-y)_i=-y_i$ and $\|y\|_2=\|-y\|_2$,
\begin{align*}
&N\sum_{y\in \mathbb{Z}^d\setminus \{\mathbf{0}\}, \|y\|_2\leq \sqrt{N}\epsilon}\|y\|_2^{-(d+\alpha)}\left(f\left(u+\frac{1}{\sqrt{N}}y\right)-f(u)\right)\\
&\approx \frac{1}{2}\sum_{i=1}^d\sum_{j=1}^d\left(\sum_{y\in \mathbb{Z}^d\setminus \{\mathbf{0}\}, \|y\|_2\leq \sqrt{N}\epsilon}\frac{y_iy_j}{\|y\|_2^{d+\alpha}}\right)\partial^2_{u_iu_j}f(u)\\
&=o(1)+\frac{1}{2}\sum_{i=1}^d\sum_{j=1}^d\left(\sum_{y\in \mathbb{Z}^d\setminus \{\mathbf{0}\}}\frac{y_iy_j}{\|y\|_2^{d+\alpha}}\right)\partial^2_{u_iu_j}f(u)\\
&=o(1)+\hat{\mathcal{L}}^{d, \alpha}f(u).
\end{align*}
Note that, when $\alpha>2$,
\begin{align*}
&\sum_{y\in \mathbb{Z}^d, \|y\|_2>\sqrt{N}\epsilon}\|y\|_2^{-(d+\alpha)}\left(f\left(u+\frac{1}{\sqrt{N}}y\right)-f(u)\right)\\
&=O(1)\int_{v:\|v\|_2\geq \epsilon \sqrt{N}}\frac{1}{\|v\|_2^{d+\alpha}}dv=O(1)\int_{\epsilon\sqrt{N}}^{+\infty}\frac{r^{d-1}}{r^{d+\alpha}}dr
=O(1)N^{-\frac{\alpha}{2}}=o(1)N^{-1}.
\end{align*}
In conclusion, when $\alpha>2$,
\[
\mathcal{L}_N^{d, \alpha}f(u)=o(1)+\hat{\mathcal{L}}^{d, \alpha}f(u)+NN^{-1}o(1)=\hat{\mathcal{L}}^{d, \alpha}f(u)+o(1)
\]
and hence $\mathcal{L}_N^{d, \alpha}f\rightarrow \hat{\mathcal{L}}^{d, \alpha}f$.

When $\alpha=2$, according to an argument similar to that in the case of $\alpha>2$,
\begin{align*}
&N\sum_{y\in \mathbb{Z}^d\setminus \{\mathbf{0}\}, \|y\|_2\leq \sqrt{N\log N}\epsilon}\|y\|_2^{-(d+2)}\left(f\left(u+\frac{1}{\sqrt{N\log N}}y\right)-f(u)\right)\\
&\approx N\sum_{i=1}^d\sum_{j=1}^d\partial^2_{u_iu_j}f(u)\left(\sum_{y\in \mathbb{Z}^d\setminus \{\mathbf{0}\}, \|y\|_2\leq \sqrt{N\log N}\epsilon}\frac{y_iy_j}{\|y\|_2^{d+2}}\frac{1}{2N\log N}\right)\\
&=N\frac{1}{2N\log N}(1+o(1))\sum_{i=1}^d\sum_{j=1}^d\partial^2_{u_iu_j}f(u)\int_1^{\epsilon\sqrt{N\log N}}\frac{1}{r^{d+2}}\left(\int_{\{v\in \mathbb{R}^d:~\|v\|_2=r\}}v_iv_jdS\right)dr\\
&=\frac{1+o(1)}{2\log N}\sum_{i=1}^d\sum_{j=1}^d\partial^2_{u_iu_j}f(u)\int_1^{\epsilon\sqrt{N\log N}}\frac{r^{d-1}r^2}{r^{d+2}}\left(\int_{\{v\in \mathbb{R}^d:~\|v\|_2=1\}}v_iv_jdS\right)dr\\
&=\frac{1+o(1)}{2}\sum_{i=1}^d\sum_{j=1}^d\partial_{u_iu_j}^2f(u)2K_{d, 2}^{i, j}\left(\frac{1}{\log N}\int_1^{\epsilon\sqrt{N\log N}}r^{-1}dr\right)\\
&=(1+o(1))\sum_{i=1}^d\sum_{j=1}^dK_{d, 2}^{i, j}\partial_{u_iu_j}^2f(u)\frac{\frac{1}{2}(1+o(1))\log N}{\log N}=o(1)+\hat{\mathcal{L}}^{d, 2}f(u)
\end{align*}
and
\begin{align*}
&\sum_{y\in \mathbb{Z}^d, \|y\|_2>\sqrt{N\log N}\epsilon}\|y\|_2^{-(d+2)}\left(f\left(u+\frac{1}{\sqrt{N\log N}}y\right)-f(u)\right)\\
&=O(1)\int_{v:\|v\|_2\geq \epsilon \sqrt{N\log N}}\frac{1}{\|v\|_2^{d+2}}dv=O(1)\int_{\epsilon\sqrt{N\log N}}^{+\infty}\frac{r^{d-1}}{r^{d+2}}dr\\
&=O(1)\left(N\log N\right)^{-1}=o(1)N^{-1}.
\end{align*}
Therefore, $\mathcal{L}_N^{d, 2}f(u)=\hat{\mathcal{L}}^{d, 2}f(u)+o(1)$.

Now we recall the characteristic functions of $X_t^{d, \alpha}$ and $\hat{X}_t^{d, \alpha}$, which are important tools for us to analyse the limit behaviors of $X_t^{d, \alpha}$ and $\hat{X}_t^{d, \alpha}$. We denote by $\psi_t^{d, \alpha}$ and $\Psi_t^{d, \alpha}$ the characteristic functions of $X_t^{d, \alpha}$ and $\hat{X}_t^{d, \alpha}$ starting from $\mathbf{0}$ respectively, i.e.,
\[
\psi_t^{d, \alpha}(\theta)=\mathbb{E}\left(\exp\left\{\sqrt{-1}\theta\cdot X_t^{d, \alpha} \right\}\Big|X_0^{d, \alpha}=\mathbf{0}\right)
\]
and
\[
\Psi_t^{d, \alpha}(\theta)=\mathbb{E}\left(\exp\left\{\sqrt{-1}\theta\cdot \hat{X}_t^{d, \alpha} \right\}\Big|\hat{X}_0^{d, \alpha}=\mathbf{0}\right)
\]
for any $\theta\in \mathbb{R}^d$. By \eqref{equ 1.2 generator longRange RW} and Kolmogorov-Chapman equation,
\begin{align*}
\frac{d}{dt}\psi_t^{d, \alpha}(\theta)&=\psi_t^{d, \alpha}(\theta)\left\{\sum_{y\in\mathbb{Z}^d\setminus\{\mathbf{0}\}}\frac{e^{\sqrt{-1}\theta\cdot y}-1}{\|y\|_2^{d+\alpha}}\right\}\\
&=\psi_t^{d, \alpha}(\theta)\left\{\sum_{y\in\mathbb{Z}^d\setminus\{\mathbf{0}\}}\frac{\cos(\theta\cdot y)-1}{\|y\|_2^{d+\alpha}}\right\}
\end{align*}
and hence
\begin{equation}\label{equ 1.3 characteristic of X_t}
\psi_t^{d, \alpha}(\theta)=\exp\left\{t\sum_{y\in\mathbb{Z}^d\setminus\{\mathbf{0}\}}\frac{\cos(\theta\cdot y)-1}{\|y\|_2^{d+\alpha}}\right\}.
\end{equation}
According to similar arguments, when $0<\alpha<2$,
\begin{equation}\label{equ 1.3 two char hat X case 1}
\Psi_t^{d, \alpha}(\theta)=\exp\left\{t\int_{\mathbb{R}^d}\frac{\cos(\theta\cdot v)-1}{\|v\|_2^{d+\alpha}}dv\right\}.
\end{equation}
When $\alpha>2$,
\begin{equation*}
\Psi_t^{d, \alpha}(\theta)=\exp\left\{-\frac{t}{2}\sum_{i=1}^d\sum_{j=1}^d\left(\sum_{y\in \mathbb{Z}^d\setminus \mathbf{\{0\}}}\frac{y_iy_j}{\|y\|_2^{d+\alpha}}\right)\theta_i\theta_j\right\}.
\end{equation*}
When $\alpha=2$,
\begin{equation}\label{equ 1.3 three char hat X case 3}
\Psi_t^{d, 2}(\theta)=\exp\left\{-\frac{t}{2}\sum_{i=1}^d\sum_{j=1}^dK_{d, 2}^{i, j}\theta_i\theta_j\right\}.
\end{equation}
Let $f_t^{d, \alpha}$ be the probability density of $\hat{X}_t^{d, \alpha}$ starting from $\mathbf{0}$. By Kolmogorov-Chapman equation, $f_t^{d, \alpha}$ satisfies
\[
\frac{d}{dt}f_t^{d, \alpha}=\hat{\mathcal{L}}^{d, \alpha}f_t^{d, \alpha}.
\]
According to the inversion formula,
\begin{equation}\label{equ 1.4 inverse formula random walk}
p_t^{d, \alpha}(\mathbf{0}, x)=\frac{1}{(2\pi)^d}\int_{[-\pi, \pi]^d}\psi_t^{d, \alpha}(\theta)e^{-\sqrt{-1}x\cdot \theta}d\theta
\end{equation}
for any $x\in \mathbb{Z}^d$ and
\begin{equation}\label{equ 1.5 inverse formula stable law}
f_t^{d, \alpha}(u)=\frac{1}{(2\pi)^d}\int_{\mathbb{R}^d}e^{-\sqrt{-1}\theta\cdot u}\Psi_t^{d, \alpha}(\theta)d\theta
\end{equation}
for any $u\in \mathbb{R}^d$.
By \eqref{equ 1.4 inverse formula random walk} and \eqref{equ 1.5 inverse formula stable law}, we have the following lemma.
\begin{lemma}\label{lemma LCLT}
For any $\alpha>0$ and $t>0$,
\[
\lim_{N\rightarrow+\infty}\sup_{x\in \mathbb{Z}^d}\left|\left(h_\alpha(N)\right)^dp_{tN}^{d, \alpha}\left(\mathbf{0}, x\right)-f_t^{d, \alpha}\left(\frac{x}{h_\alpha(N)}\right)\right|=0.
\]
\end{lemma}
Lemma \ref{lemma LCLT} in the case $\alpha>2$ is a corollary of Theorem 2.1.3 in Chapter 2 of \cite{Lawler2010}. We prove Lemma \ref{lemma LCLT} in the case $0<\alpha\leq 2$ in Section \ref{section three}.

By Lemma \ref{lemma LCLT},
\begin{equation}\label{equ 1.6}
p_{t}^{d, \alpha}(\mathbf{0}, \mathbf{0})=\Theta(1)\left(h_\alpha(t)\right)^{-d}.
\end{equation}
Since $\{X_t^{d, \alpha}\}_{t\geq 0}$ is transient when and only when
\[
\int_0^{+\infty}p_t^{d, \alpha}(\mathbf{0}, \mathbf{0})dt<+\infty,
\]
we have the following conclusion by \eqref{equ 1.6}. When $d=1, 2$, $\{X_t^{d, \alpha}\}_{t\geq 0}$ is transient when and only when $0<\alpha<d$. When $d\geq 3$, $\{X_t^{d, \alpha}\}_{t\geq 0}$ is transient for all $\alpha>0$.

\subsection{Duality}\label{subsection 1.3}
In this subsection, we recall the duality relationship between the voter model and the coalescing random walk. For given $x_1, x_2, \ldots, x_m\in \mathbb{Z}^d$ and $t_1, t_2, \ldots, t_m\geq 0$, the coalescing long-range random walk
\[
\left\{\left(\mathcal{X}_{t, x_1}^{d, \alpha, t_1}, \mathcal{X}_{t, x_2}^{d, \alpha, t_2}, \ldots, \mathcal{X}_{t, x_m}^{d, \alpha, t_m}\right)\right\}_{t\geq 0}
\]
on $\mathbb{Z}^d$ with parameter $\alpha>0$ is defined as follows. For each $1\leq k\leq m$, $\mathcal{X}_{t, x_k}^{d, \alpha, t_k}=x_k$ when $0\leq t\leq t_k$ and $\{\mathcal{X}_{t_k+s, x_k}^{d, \alpha, t_k}\}_{s\geq 0}$ is a copy of the long-range random walk $\{X_s^{d, \alpha}\}_{s\geq 0}$ starting from $x_k$. We say $\mathcal{X}_{t, x_k}^{d, \alpha, t_k}$ active when and only when $t>t_k$. The joint distribution of active walks is defined as follows. All active walks are independent until collisions occurs. When two active walks collide with each other, they are coalesced into a single active walk.

For later use, we denote by $\mathbb{P}$ the probability measure of the long-range random walks and the coalescing long-range walks. We denote by $\mathbb{E}$ the expectation with respect to $\mathbb{P}$. For any probability measure $\mu$ on $\{0, 1\}^{\mathbb{Z}^d}$, we denote by $\hat{\mathbb{P}}_\mu$ the probability measure of the long-range voter model $\{\eta_t\}_{t\geq 0}$ starting from $\mu$. We denote by $\hat{\mathbb{E}}_\mu$ the expectation with respect to $\hat{\mathbb{P}}_\mu$. When $\mu$ is the Dirac measure concentrated on $\eta\in \{0, 1\}^{\mathbb{Z}^d}$, we write $\hat{\mathbb{P}}_\mu$ and $\hat{\mathbb{E}}_\mu$ as $\hat{\mathbb{P}}_\eta$ and $\hat{\mathbb{E}}_\eta$ respectively.

We have the following duality relationship between the voter model and the coalescing long-range random walk.

\begin{proposition}\label{proposition voter dual}
For any $\eta\in \{0, 1\}^{\mathbb{Z}^d}$, $x_1, \ldots, x_m\in \mathbb{Z}^d$, $s_1, s_2, \ldots, s_m\geq 0$,
\[
u\geq \max\left\{s_1, \ldots, s_m\right\}
\]
and $h: \{0, 1\}^m\rightarrow \mathbb{R}$, we have
\begin{align}\label{equ voter dual}
&\hat{\mathbb{E}}_\eta h\left(\eta_{s_1}(x_1), \eta_{s_2}(x_2), \ldots, \eta_{s_m}(x_m)\right)\notag\\
&=\mathbb{E}h\left(\eta\left(\mathcal{X}_{u, x_1}^{d, \alpha, u-s_1}\right),
\eta\left(\mathcal{X}_{u, x_2}^{d, \alpha, u-s_2}\right), \ldots, \eta\left(\mathcal{X}_{u, x_m}^{d, \alpha, u-s_m}\right)\right).
\end{align}
\end{proposition}

Here we recall the outline of a proof of Proposition \ref{proposition voter dual}. For details, see Section 3.4 of \cite{Lig1985}. For $0<s\leq t$ and $x\in \mathbb{Z}^d$, let $\mathcal{Z}_{s, t}^x\in \mathbb{Z}^d$ be the provider of $\eta_t(x)$ at moment $t-s$.
Then,
\[
\eta_{t-s}\left(\mathcal{Z}_{s, t}^x\right)=\eta_t(x)
\]
ad hence
\[
h\left(\eta_{t}(x_1), \eta_{t}(x_2), \ldots, \eta_{t}(x_m)\right)=h\left(\eta_0(\mathcal{Z}_{t, t}^{x_1}), \ldots, \eta_0(\mathcal{Z}_{t, t}^{x_m})\right).
\]
For $r<t$, if $\mathcal{Z}_{r-, t}^x$ adopts $y$'s opinion at moment $t-r$, then $\mathcal{Z}_{r, t}^x=y$, i.e, the provider of $\eta_t(x)$ jumps from $\mathcal{Z}_{r-, t}^x$ to $y$. Hence, $\{Z_{s, t}^x\}_{0\leq s\leq t}$ is a copy of $\{X_s^{d, \alpha}\}_{0\leq s\leq t}$ starting from $x$. If $y$ is the common provider of $\eta_t(x)$ and $\eta_t(w)$ at moment $t-s$, then for any $u<t-s$,
\[
\mathcal{Z}_{t-u, t}^x=\mathcal{Z}_{t-u, t}^w=\mathcal{Z}_{t-s-u, t-s}^y.
\]
Therefore, $\{\left(\mathcal{Z}_{s, t}^{x_1}, \ldots, \mathcal{Z}_{s, t}^{x_m}\right)\}_{0\leq s\leq t}$ evolves as the coalescing long-range random walk. Consequently, Proposition \ref{proposition voter dual} holds when $s_1=s_2=\ldots=s_m=t$ and $u=t$. For $s_1, s_2,\ldots, s_m$ which are not the same, we can choose $u\geq \max\{s_1, \ldots, s_m\}$ and define
\[
\hat{\mathcal{Z}}_{s, u}^{x_k, s_k}=
\begin{cases}
x_k & \text{~if~}s\leq u-s_k,\\
\mathcal{Z}_{s-(u-s_k), s_k}^{x_k} & \text{~if~}u-s_k\leq s\leq u
\end{cases}
\]
for $1\leq k\leq m$. Then, similar to the case $s_1=\ldots=s_m$,
\[
h\left(\eta_{s_1}(x_1), \eta_{s_2}(x_2), \ldots, \eta_{s_m}(x_m)\right)=h\left(\eta_0(\hat{\mathcal{Z}}_{u, u}^{x_1, s_1}), \ldots, \eta_0(\hat{\mathcal{Z}}_{u, u}^{x_m, s_m})\right)
\]
and $\left(\hat{\mathcal{Z}}_{u, u}^{x_1, s_1}, \ldots, \hat{\mathcal{Z}}_{u, u}^{x_m, s_m}\right)$ is a copy of $\left(\mathcal{X}_{u, x_1}^{d, \alpha, u-s_1}, \ldots, \mathcal{X}_{u, x_m}^{d, \alpha, u-s_m}\right)$, which completes the proof of Proposition \ref{proposition voter dual} in all cases.

For $0<p<1$, let $\mu_p$ be the product measure on $\mathbb{Z}^d$ with density $p$, i.e., $\{\eta(x)\}_{x\in \mathbb{Z}^d}$ are independent under $\mu_p$ and
\[
\mu_p\{\eta: \eta(x)=1\}=p
\]
for all $x\in \mathbb{Z}^d$. By Proposition \ref{proposition voter dual}, for any $x_1,\ldots, x_m\in \mathbb{Z}^d$ and $t\geq 0$,
\begin{equation}\label{equ 1.10}
\hat{\mathbb{P}}_{\mu_p}\left(\eta_t(x_i)=1\text{~for all~}1\leq i\leq m\right)
=\mathbb{E}p^{\mathcal{A}_t^{d, \alpha, x_1, \ldots, x_m}},
\end{equation}
where $\mathcal{A}_t^{d, \alpha, x_1, \ldots, x_m}$ is the cardinality of $\{\mathcal{X}_{t, x_1}^{d, \alpha, 0}, \mathcal{X}_{t, x_2}^{d, \alpha, 0}, \ldots, \mathcal{X}_{t, x_m}^{d, \alpha, 0}\}$. Since $\mathcal{A}_t^{d, \alpha, x_1, \ldots, x_m}$ is decreasing in $t$, there exists probability measure $\nu_p^{d, \alpha}$ on $\{0, 1\}^{\mathbb{Z}^d}$ such that $\{\eta_t\}_{t\geq 0}$ starting from $\mu_p$ converges weakly to $\nu_p^{d, \alpha}$ as $t\rightarrow+\infty$ and
\[
\nu_p^{d, \alpha}\left(\eta: \eta(x_i)=1\text{~for all~}1\leq i\leq m\right)=\lim_{t\rightarrow+\infty}\mathbb{E}p^{\mathcal{A}_t^{d, \alpha, x_1, \ldots, x_m}} \]
for any $x_1, \ldots, x_m\in \mathbb{Z}^d$. If $d, \alpha$ makes $\{X_t^{d, \alpha}\}_{t\geq 0}$ recurrent, then $\mathcal{X}_{t, x_1}^{d, \alpha, 0}, \mathcal{X}_{t, x_2}^{d, \alpha, 0}, \ldots, \mathcal{X}_{t, x_m}^{d, \alpha, 0}$ are coalesced into a single random walk finally with probability one and hence
\[
\nu_p^{d, \alpha}\left(\eta: \eta(x_i)=1\text{~for all~}1\leq i\leq m\right)=p.
\]
So in this case, $\nu_p^{d, \alpha}$ is a convex combination of $\delta_0$ and $\delta_1$, where $\delta_a$ is the Dirac measure concentrated on the configuration where all vertices take value $a$. In this paper, we are concerned with the long-range voter model starting from $\nu_p^{d, \alpha}$. To avoid above trivial case, we require that $(d, \alpha)$ makes $\{X_t^{d, \alpha}\}_{t\geq 0}$ transient, i.e,
\begin{equation}\label{equ 1.11}
(d, \alpha)\in \{1\}\times (0, 1)\bigcup \{2\}\times (0, 2)\bigcup \{3, 4, \ldots\}\times (0, +\infty)
\end{equation}
as we have introduced in Subsection \ref{subsection 1.2}.

For later use, we recall some properties of $\nu_p^{d, \alpha}$. According to the definition of $\nu_p^{d, \alpha}$, for $x, y\in \mathbb{Z}^d$,
\[
\hat{\mathbb{E}}_{\nu_p^{d, \alpha}}\left((\eta_0(x)-\eta_0(y))^2\right)=\lim_{t\rightarrow+\infty}\hat{\mathbb{E}}_{\mu_p}\left((\eta_t(x)-\eta_t(y))^2\right).
\]
By Proposition \ref{proposition voter dual},
\[
\hat{\mathbb{E}}_{\mu_p}\left((\eta_t(x)-\eta_t(y))^2\right)=2p(1-p)(1-0)^2\mathbb{P}(\tau_{x, y}>t),
\]
where
\[
\tau_{x, y}=\inf\left\{t>0:~\mathcal{X}_{t, x}^{d, \alpha, 0}=\mathcal{X}_{t, y}^{d, \alpha, 0}\right\}.
\]
Before $\tau_{x, y}$, $\mathcal{X}_{t, x}^{d, \alpha, 0}-\mathcal{X}_{t, y}^{d, \alpha, 0}$ evolves as a copy of $X_{2t}^{d, \alpha}$ starting from $x-y$. Hence,
\[
\lim_{t\rightarrow+\infty}\mathbb{P}(\tau_{x, y}>t)=\Phi^{d, \alpha}(x-y),
\]
where
\begin{equation}\label{equ 1.12}
\Phi^{d, \alpha}(x)=\mathbb{P}\left(X_t^{d, \alpha}\neq \mathbf{0} \text{~for all~}t\geq 0\Big|X_0^{d, \alpha}=x\right).
\end{equation}
Therefore,
\begin{equation}\label{equ 1.13}
\hat{\mathbb{E}}_{\nu_p^{d, \alpha}}\left((\eta_0(x)-\eta_0(y))^2\right)=2p(1-p)\Phi^{d, \alpha}(x-y).
\end{equation}
By \eqref{equ 1.10},
\begin{equation}\label{equ 1.14}
\hat{\mathbb{E}}_{\nu_p}\eta_0(x)=\nu_p^{d, \alpha}\left(\eta:~\eta(x)=1\right)=\lim_{t\rightarrow+\infty}p^1=p.
\end{equation}
Then, by \eqref{equ 1.13},
\begin{equation*}
\hat{\mathbb{E}}_{\nu_p^{d, \alpha}}\left(\eta_0(x)\eta_0(y)\right)=\frac{2p-2p(1-p)\Phi^{d, \alpha}(x-y)}{2}=p\left(1-\Phi^{d, \alpha}(x-y)\right)
+p^2\Phi^{d, \alpha}(x-y)
\end{equation*}
and hence
\begin{equation}\label{equ 1.15}
{\rm Cov}_{\nu_p^{d, \alpha}}(\eta_0(x), \eta_0(y))=p(1-p)\left(1-\Phi^{d, \alpha}(x-y)\right).
\end{equation}

\section{Main results}\label{section two}
In this section, we give our main result. As we have introduced in Subsection \ref{subsection 1.3}, from now on we assume that $(d, \alpha)$ satisfies \eqref{equ 1.11}. We first introduce some notations and definitions for later use. We denote by $\{\mathcal{W}_t\}_{t\geq 0}$ the standard Brownian motion starting from $0$. For any $\beta\in (0, 1)$, we denote by $\{B_t^\beta\}_{t\geq 0}$ the fractional Brownian motion with Hurst parameter $\beta$, i.e., $\{B_t^\beta\}_{t\geq 0}$ is the Gaussian process with continuous sample path, mean zero and covariance functions given by
\[
{\rm Cov}\left(B_t^\beta, B_s^\beta\right)=\frac{1}{2}\left(t^{2\beta}+s^{2\beta}-|t-s|^{2\beta}\right)
\]
for any $t, s>0$. Note that $\{\mathcal{W}_t\}_{t\geq 0}=\{B_t^{1/2}\}_{t\geq 0}$.

For $p\in (0, 1)$, let $\nu_p^{d, \alpha}$ be defined as in Subsection \ref{subsection 1.3}. By \eqref{equ 1.14}, for $\{\eta_t\}_{t\geq 0}$ starting from $\nu_p^{d, \alpha}$, the centered occupation time on $\mathbf{0}$ during $[0, t]$ is defined as
\[
\int_0^t\left(\eta_s(\mathbf{0})-p\right)ds.
\]
For $(d, \alpha)$ satisfies \eqref{equ 1.11},
we define
\[
\mathcal{C}_{d, \alpha}=\frac{1}{\int_0^{+\infty}p_s^{d, \alpha}(\mathbf{0}, \mathbf{0})ds},
\]
where $\{p_t^{d, \alpha}(\cdot, \cdot)\}_{t\geq 0}$ are transition probabilities of $\{X_t^{d, \alpha}\}_{t\geq 0}$ defined as in Subsection \ref{subsection 1.2}. Note that $\mathcal{C}_{d, \alpha}>0$ since $\{X_t^{d, \alpha}\}_{t\geq 0}$ is transient when $(d, \alpha)$ satisfies \eqref{equ 1.11}.
Now we give our main result in the case $d\in\{1, 2, 3\}$ and $0<\alpha<\min\{2, d\}$.

\begin{theorem}\label{theorem 2.1 main result d=1, 2}
Let $d\in \{1, 2, 3\}$, $p\in (0, 1)$, $\alpha\in \left(0, \min\{d, 2\}\right)$ and $T>0$. If $\eta_0$ is distributed with $\nu_p^{d, \alpha}$, then
\[
\left\{\frac{1}{\Lambda_{d, \alpha}(N)}\int_0^{tN}\left(\eta_s(\mathbf{0})-p\right)ds:~0\leq t\leq T\right\}
\]
converges weakly, with respect to the uniform topology of $C[0, T]$, to
\[
\begin{cases}
\left\{\sqrt{2f_1^{d, \alpha}(\mathbf{0})\frac{\alpha^3\mathcal{C}_{d, \alpha}p(1-p)}{(d-\alpha)(2\alpha-d)(3\alpha-d)}}B_t^{\frac{3}{2}-\frac{d}{2\alpha}}:~0\leq t\leq T\right\}
\text{\quad if\quad} \frac{d}{2}<\alpha<\min\{2, d\},\\
\left\{\sqrt{2f_1^{d, \frac{d}{2}}(\mathbf{0})\mathcal{C}_{d, \frac{d}{2}}p(1-p)}\mathcal{W}_t:~0\leq t\leq T\right\}  \text{\quad if \quad} \alpha=\frac{d}{2},\\
\left\{\sqrt{2\int_0^{+\infty}\int_0^{+\infty}p_{u+r}^{d, \alpha}(\mathbf{0}, \mathbf{0})dudr\mathcal{C}_{d, \alpha}p(1-p)}\mathcal{W}_t:~0\leq t\leq T\right\}
 \text{\quad if\quad} 0<\alpha<\frac{d}{2}
\end{cases}
\]
as $N\rightarrow+\infty$, where $f_1^{d, \alpha}$ is the probability density of $\hat{X}_1^{d, \alpha}$ defined as in Subsection \ref{subsection 1.2} and
\[
\Lambda_{d, \alpha}(N)=
\begin{cases}
N^{\frac{3}{2}-\frac{d}{2\alpha}} & \text{~if\quad} \frac{d}{2}<\alpha<\min\{d, 2\},\\
\sqrt{N\log N} & \text{~if\quad} \alpha=\frac{d}{2},\\
\sqrt{N} & \text{~if\quad} 0<\alpha<\frac{d}{2}.
\end{cases}
\]
\end{theorem}

\quad

By Proposition \ref{proposition voter dual}, it is not difficult to show that
\[
{\rm Var}_{\nu_p^{d, \alpha}}\left(\int_0^t\left(\eta_s(\mathbf{0})-p\right)ds\right)
=\Theta(1)\int_0^t\left(\int_0^s\left(\int_0^{+\infty}p_{u+r}^{d, \alpha}(\mathbf{0}, \mathbf{0})dr\right)du\right)ds.
\]
For mathematical details of how to check the above equation, see Section \ref{section four}. Then,
\[
{\rm Var}_{\nu_p^{d, \alpha}}\left(\int_0^{tN}\left(\eta_s(\mathbf{0}-p)\right)ds\right)=\Theta(N)
\]
when and only when
\[
\int_0^{+\infty}\int_0^{+\infty}p_{u+r}^{d, \alpha}(\mathbf{0}, \mathbf{0})dudr<+\infty.
\]
By Lemma \ref{lemma LCLT}, when $d\in \{1, 2, 3\}$,
\[
\int_0^{+\infty}\int_0^{+\infty}p_{u+r}^{d, \alpha}(\mathbf{0}, \mathbf{0})dudr<+\infty
\]
when and only when $0<\alpha<\frac{d}{2}$. Hence, the phase transition in Theorem \ref{theorem 2.1 main result d=1, 2} occurs at $\alpha=\frac{d}{2}$.

Now we give our main result in the case $d=3$ and $\alpha\geq 2$. Note that the case $d=3$ and $\alpha\in (0, 2)$ is included in Theorem \ref{theorem 2.1 main result d=1, 2}.
\begin{theorem}\label{theorem 2.2 main result d=3}
Let $d=3$, $p\in (0, 1)$, $\alpha\in [2, +\infty)$ and $T>0$. If $\eta_0$ is distributed with $\nu_p^{3, \alpha}$, then
\[
\left\{\frac{1}{\Lambda_{3, \alpha}(N)}\int_0^{tN}\left(\eta_s(\mathbf{0})-p\right)ds:~0\leq t\leq T\right\}
\]
converges weakly, with respect to the uniform topology of $C[0, T]$, to
\[
\left\{\sqrt{\frac{16}{3}f_1^{3, \alpha}(\mathbf{0})\mathcal{C}_{3, \alpha}p(1-p)}B_t^{\frac{3}{4}}:~0\leq t\leq T:~0\leq t\leq T\right\}
\]
as $N\rightarrow+\infty$, where $f_1^{3, \alpha}$ is the probability density of $\hat{X}_1^{3, \alpha}$ defined as in Subsection \ref{subsection 1.2} and
\[
\Lambda_{3, \alpha}(N)=
\begin{cases}
N^{\frac{3}{4}} & \text{~if\quad} \alpha>2,\\
N^{\frac{3}{4}}\left(\log N\right)^{-\frac{3}{4}} & \text{~if\quad} \alpha=2.
\end{cases}
\]
\end{theorem}

\quad

The phase transition in Theorem \ref{theorem 2.2 main result d=3} occurs at $\alpha=2$ since the phase transition in Lemma \ref{lemma LCLT} occurs at $\alpha=2$.

Now we give our main result in the case $d=4$.
\begin{theorem}\label{theorem 2.3 main result d=4}
Let $d=4$, $p\in (0, 1)$, $\alpha\in (0, +\infty)$ and $T>0$. If $\eta_0$ is distributed with $\nu_p^{4, \alpha}$, then
\[
\left\{\frac{1}{\Lambda_{4, \alpha}(N)}\int_0^{tN}\left(\eta_s(\mathbf{0})-p\right)ds:~0\leq t\leq T\right\}
\]
converges weakly, with respect to the uniform topology of $C[0, T]$, to
\[
\begin{cases}
\left\{\sqrt{2f_1^{4, \alpha}(\mathbf{0})\mathcal{C}_{4, \alpha}p(1-p)}\mathcal{W}_t:~0\leq t\leq T:~0\leq t\leq T\right\}  \text{\quad if\quad} \alpha>2,\\
\left\{\sqrt{2\int_0^{+\infty}\int_0^{+\infty}p_{u+r}^{4, \alpha}(\mathbf{0}, \mathbf{0})dudr\mathcal{C}_{4, \alpha}p(1-p)}\mathcal{W}_t:~0\leq t\leq T\right\}
 \text{\quad if\quad} 0<\alpha\leq 2
\end{cases}
\]
as $N\rightarrow+\infty$, where $f_1^{4, \alpha}$ is the probability density of $\hat{X}_1^{4, \alpha}$ defined as in Subsection \ref{subsection 1.2} and
\[
\Lambda_{4, \alpha}(N)=
\begin{cases}
\sqrt{N\log N} & \text{~if\quad} \alpha>2,\\
\sqrt{N} & \text{~if\quad} 0<\alpha\leq 2.
\end{cases}
\]
\end{theorem}

\quad

By Lemma \ref{lemma LCLT}, when $d=4$,
\[
\int_0^{+\infty}\int_0^{+\infty}p_{u+r}^{4, \alpha}(\mathbf{0}, \mathbf{0})dudr<+\infty
\]
when and only when $0<\alpha\leq 2$. Hence, the phase transition in Theorem \ref{theorem 2.3 main result d=4} occurs at $\alpha=2$.

At last, we give our main result in the case $d\geq 5$.
\begin{theorem}\label{theorem 2.4 main result d geq 5}
Let $d\geq 5$, $p\in (0, 1)$, $\alpha\in (0, +\infty)$ and $T>0$. If $\eta_0$ is distributed with $\nu_p^{d, \alpha}$, then
\[
\left\{\frac{1}{\sqrt{N}}\int_0^{tN}\left(\eta_s(\mathbf{0})-p\right)ds:~0\leq t\leq T\right\}
\]
converges weakly, with respect to the uniform topology of $C[0, T]$, to
\[
\left\{\sqrt{2\int_0^{+\infty}\int_0^{+\infty}p_{u+r}^{d, \alpha}(\mathbf{0}, \mathbf{0})dudr\mathcal{C}_{d, \alpha}p(1-p)}\mathcal{W}_t:~0\leq t\leq T\right\}
\]
as $N\rightarrow+\infty$, where $f_1^{d, \alpha}$ is the probability density of $\hat{X}_1^{d, \alpha}$ defined as in Subsection \ref{subsection 1.2}.
\end{theorem}

\quad

By Lemma \ref{lemma LCLT}, when $d\geq 5$,
\[
\int_0^{+\infty}\int_0^{+\infty}p_{u+r}^{d, \alpha}(\mathbf{0}, \mathbf{0})dudr<+\infty
\]
for all $\alpha>0$. Hence, there is no phase transition in Theorem \ref{theorem 2.4 main result d geq 5}.

The proofs of Theorems \ref{theorem 2.1 main result d=1, 2}-\ref{theorem 2.4 main result d geq 5} utilize a strategy similar to that given in \cite{Xue2024} to deal with the nearest-neighbor case, where the resolvent strategy introduced in \cite{Kipnis1987} and the Poisson flow strategy introduced in \cite{Birkner2007} are applied and the heat estimation Lemma \ref{lemma LCLT} plays the key role. By the resolvent strategy, we can decompose the centered occupation time as a martingale plus a remainder term. By the Poisson flow strategy, we can write the martingale term in the above decomposition as a Lebesgue-Stieltjes integral with respect to a series of centered Poisson processes. Then we can show that this martingale term converges weakly to a Gaussian process by computing the quadratic variation process. The limit behavior of the remainder term in the aforesaid decomposition is a main difference between this paper and \cite{Xue2024}. In \cite{Xue2024}, the nearest-neighbor voter model starts from the product measure, which ensures that the remainder term converges weakly to $0$. In this paper, however, the long-range voter model starts from the stationary distribution $\nu_p^{d, \alpha}$, due to which we need to show that the remainder term converges weakly to a Gaussian distribution in several cases. To deal with the remainder term, we give a Gaussian limit of the fluctuation density field of the long-range voter model, which is an analogue of nearest-neighbor case given in \cite{Presutti1983}. Then, we can approximately write the remainder term as the integral of a Schwartz function with respect to the  fluctuation density field and consequently the conclusion that the remainder term converges weakly to a Gaussian distribution follows from our fluctuation limit theorem. We can further show that the joint distribution of the aforesaid martingale term and remainder term converges weakly to the product measure of the respective Gaussian limits of these two terms. As a result, the centered occupation time converges weakly to a Gaussian process and then a computing of the covariance functions completes the proof. For mathematical details, see Sections \ref{section three} to \ref{section six}.

\section{Proof of Lemma \ref{lemma LCLT}}\label{section three}
In this section, we prove Lemma \ref{lemma LCLT}. As we have explained in Section \ref{section one}, we only need to deal with $\alpha\in (0, 2]$. We first give the proof of Lemma \ref{lemma LCLT} in the case $0<\alpha<2$.

\proof
[Proof of Lemma \ref{lemma LCLT} in the case $0<\alpha<2$]
When $0<\alpha<2$, by \eqref{equ 1.3 characteristic of X_t} and \eqref{equ 1.4 inverse formula random walk},
\begin{align}\label{equ 3.1}
N^{\frac{d}{\alpha}}p_{tN}^{d, \alpha}\left(\mathbf{0}, x\right)
&=\frac{N^{\frac{d}{\alpha}}}{(2\pi)^d}\int_{[-\pi, \pi]^d}\psi_{tN}^{d, \alpha}(\theta)e^{-\sqrt{-1}x\cdot \theta}d\theta\notag\\
&=\frac{1}{(2\pi)^d}\int_{[-N^{\frac{1}{\alpha}}\pi, N^{\frac{1}{\alpha}}\pi]^d}\psi_{tN}^{d, \alpha}(\theta/N^{\frac{1}{\alpha}})e^{-\sqrt{-1}\frac{x}{N^{\frac{1}{\alpha}}}\cdot \theta}d\theta,
\end{align}
where
\[
\psi_{tN}^{d, \alpha}(\theta/N^{\frac{1}{\alpha}})=e^{-t\|\theta\|_2^\alpha G^{d, \alpha}\left(\frac{\theta}{N^{\frac{1}{\alpha}}}\right)}
\]
with
\[
G^{d, \alpha}(u)=\sum_{y\in\mathbb{Z}^d\setminus\{\mathbf{0}\}}\frac{1-\cos\left(u\cdot y\right)}{\|y\|_2^{d+\alpha}\|u\|_2^{\alpha}}
\]
for any $u\in \mathbb{R}^d\setminus\{\mathbf{0}\}$. By \eqref{equ 1.3 two char hat X case 1}, when $0<\alpha<2$,
\[
\Psi_t^{d, \alpha}(\theta)=\exp\left\{-t\|\theta\|_2^\alpha\tilde{G}^{d, \alpha}\right\},
\]
where
\[
\tilde{G}^{d, \alpha}=\int_{\mathbb{R}^d}\frac{\left(1-\cos(\varsigma\cdot v)\right)}{\|v\|_2^{d+\alpha}}dv
\]
with $\varsigma$ being any element in $\mathbb{R}^d$ such that $\|\varsigma\|_2=1$. According to the definition of Riemann integral, for any given $\theta\in \mathbb{R}^d$,
\[
\lim_{N\rightarrow+\infty}G^{d, \alpha}\left(\theta/N^{\frac{1}{\alpha}}\right)=\tilde{G}^{d, \alpha}.
\]
Therefore, by the dominated convergence theorem, for any $M\geq 0$,
\begin{align*}
\lim_{N\rightarrow+\infty}\int_{[-M, M]^d}\left|\psi_{tN}^{d, \alpha}(\theta/N^{\frac{1}{\alpha}})-\Psi_t^{d, \alpha}(\theta)\right|e^{-\sqrt{-1}\frac{x}{N^{\frac{1}{\alpha}}}\cdot \theta}d\theta=0.
\end{align*}
Then, by \eqref{equ 1.5 inverse formula stable law} and \eqref{equ 3.1}, to complete the proof we only need to show that
\begin{equation}\label{equ 3.2}
\lim_{M\rightarrow+\infty}\limsup_{N\rightarrow+\infty}\int_{\theta:M<\|\theta\|_\infty\leq N^{\frac{1}{\alpha}}\pi}e^{-t\|\theta\|_2^\alpha G^{d, \alpha}\left(\frac{\theta}{N^{\frac{1}{\alpha}}}\right)}d\theta=0,
\end{equation}
where $\|\cdot\|_{\infty}$ is the $l_\infty$-norm on $\mathbb{R}^d$. For any $\varpi\in \mathbb{R}^d$ such that $\|\varpi\|_\infty=1$ and $s>0$, we define
\[
H_1(\varpi, s)=s^d\sum_{y\in\mathbb{Z}^d\setminus\{\mathbf{0}\}}\frac{1-\cos\left(\varpi\cdot(sy)\right)}{\|sy\|_2^{d+\alpha}}.
\]
Then,
\[
G^{d, \alpha}\left(\frac{\theta}{N^{\frac{1}{\alpha}}}\right)=H_1\left(\frac{\theta}{\|\theta\|_\infty}, \frac{\|\theta\|_\infty}{N^{\frac{1}{\alpha}}}\right)\left(\frac{\|\theta\|_\infty}{\|\theta\|_2}\right)^\alpha.
\]
Since
\[
\lim_{s\rightarrow 0}H_1(\varpi, s)=\int_{\mathbb{R}^d}\frac{1-\cos(\varpi\cdot v)}{\|v\|_2^{d+\alpha}}dv=\tilde{G}^{d, \alpha}\|\varpi\|_2^\alpha,
\]
we supplementarily define $H_1(\varpi, 0)=\tilde{G}^{d, \alpha}\|\varpi\|_2^\alpha$ and then $H_1$ is a continuous function on
\[
\left\{\varpi\in \mathbb{R}^d:~\|\varpi\|_\infty=1\right\}\times [0, +\infty).
\]
Note that $H_1(\varpi, 0)=\tilde{G}^{d, \alpha}\|\varpi\|_2^\alpha>0$ for any $\varpi$ such that $\|\varpi\|_\infty=1$ and $H_1(\varpi, s)\geq 0$ for any $\varpi$ such that $\|\varpi\|_\infty=1$ and $s>0$. Furthermore, when $s>0$, $H_1(\varpi, s)=0$ when and only when $\cos(\varpi\cdot (sy))=1$ for all $y\in \mathbb{Z}^d\setminus\{\mathbf{0}\}$. Hence, when $0<s\leq \pi$, $H_1(\varpi, s)>0$ for all $\varpi$ such that $\|\varpi\|_\infty=1$. Therefore,
\[
\min\left\{H_1(\varpi, s):~\|\varpi\|_\infty=1, 0\leq s\leq \pi\right\}>0.
\]
As a result, for $\theta$ such that $M\leq \|\theta\|_\infty\leq N^{\frac{1}{\alpha}}\pi$,
\begin{equation}\label{equ 3.3}
\exp\left\{-t\|\theta\|_2^\alpha G^{d, \alpha}\left(\frac{\theta}{N^{\frac{1}{\alpha}}}\right)\right\}
\leq \exp\left\{-t\|\theta\|_2^\alpha Q_1^{d, \alpha}\left(\frac{1}{\sqrt{d}}\right)^\alpha\right\},
\end{equation}
where
\[
Q_1^{d, \alpha}=\min\left\{H_1(\varpi, s):~\|\varpi\|_\infty=1, 0\leq s\leq \pi\right\}>0.
\]
Equation \eqref{equ 3.2} follows from \eqref{equ 3.3} and the proof is complete.
\qed


Now we give the proof of Lemma \ref{lemma LCLT} in the case $\alpha=2$.
\proof
[Proof of Lemma \ref{lemma LCLT} in the case $\alpha=2$] According to an analysis similar to that leading to $\mathcal{L}_N^{d, 2}\rightarrow \hat{\mathcal{L}}^{d, 2}$ as given in Section \ref{section one}, we have
\[
\lim_{N\rightarrow+\infty}\psi_{tN}^{d, 2}\left(\theta\big/\sqrt{N\log N}\right)=\Psi_t^{d, 2}(\theta)
\]
for all $\theta\in \mathbb{R}^d$. Then by inversion formula, to complete the proof we only need to check the following analogue of \eqref{equ 3.2},
\begin{equation}\label{equ 3.4}
\lim_{M\rightarrow+\infty}\limsup_{N\rightarrow+\infty}\int_{\theta:M<\|\theta\|_\infty\leq \sqrt{N\log N}\pi}\psi_{tN}^{d, 2}\left(\theta\big/\sqrt{N\log N}\right)d\theta=0.
\end{equation}
Note that
\[
\psi_{tN}^{d, 2}\left(\theta\big/\sqrt{N\log N}\right)=\exp\left\{-tN\mathcal{R}\left(\theta\big/\sqrt{N\log N}\right)\right\},
\]
where
\[
\mathcal{R}(u)=\sum_{y\in \mathbb{Z}^d\setminus\{\mathbf{0}\}}\frac{1-\cos(u\cdot y)}{\|y\|_2^{d+2}}
\]
for all $u\in \mathbb{R}^d$. Note that $\mathcal{R}(u)\geq 0$ for all $u\in \mathbb{R}^d$. Furthermore, $\mathcal{R}(u)=0$ when and only when $\cos(u\cdot y)=1$ for all $y\in \mathbb{Z}^d\setminus\{\mathbf{0}\}$. Therefore, for any sufficiently small $\epsilon>0$, we have $\mathcal{R}(u)>0$ when $\epsilon\leq \|u\|_\infty\leq \pi$ and hence
\[
Q_2^\epsilon:=\min\left\{\mathcal{R}(u):~\epsilon \leq \|u\|_\infty\leq \pi\right\}>0.
\]
As a result,
\begin{align*}
&\int_{\theta:\sqrt{N\log N}\epsilon<\|\theta\|_\infty\leq \sqrt{N\log N}\pi}\psi_{tN}^{d, 2}\left(\theta\big/\sqrt{N\log N}\right)d\theta\\
&\leq \int_{\theta:\sqrt{N\log N}\epsilon<\|\theta\|_\infty\leq \sqrt{N\log N}\pi}e^{-tNQ_2^\epsilon}d\theta\rightarrow 0
\end{align*}
as $N\rightarrow+\infty$. Hence, to check \eqref{equ 3.4}, we only need to show that
\begin{equation}\label{equ 3.5}
\lim_{M\rightarrow+\infty}\limsup_{N\rightarrow+\infty}\int_{\theta:M<\|\theta\|_\infty\leq \epsilon\sqrt{N\log N}}\psi_{tN}^{d, 2}\left(\theta\big/\sqrt{N\log N}\right)d\theta=0
\end{equation}
for any sufficiently small $\epsilon>0$. For any $\varpi\in \mathbb{R}^d$ such that $\|\varpi\|_\infty=1$ and $s>0$, we define
\[
H_2(\varpi, s)=s^d\sum_{y\in \mathbb{Z}^d\setminus\{\mathbf{0}\}}\frac{1-\cos\left(\varpi\cdot(sy)\right)}{\|sy\|_2^{d+2}}.
\]
Then,
\begin{equation}\label{equ 3.5 two}
\psi_{tN}^{d, 2}\left(\theta\big/\sqrt{N\log N}\right)=\exp\left\{-\frac{t}{\log N}H_2\left(\frac{\theta}{\|\theta\|_\infty}, \frac{\|\theta\|_\infty}{\sqrt{N\log N}}\right)\|\theta\|_\infty^2\right\}.
\end{equation}
For any $s>0$,
\[
H_2(\varpi, s)\geq s^d\sum_{y\in \mathbb{Z}^d\setminus\{\mathbf{0}\}: \|y\|_2\leq \frac{1}{2\sqrt{d}s}}\frac{1-\cos\left(\varpi\cdot(sy)\right)}{\|sy\|_2^{d+2}}.
\]
When $0\leq a\leq \frac{1}{2}$, there exists constant $C_1>0$ independent of $a$ such that $1-\cos(a)\geq C_1a^2$.
When $\|y\|_2\leq \frac{1}{2s\sqrt{d}}$, $|\varpi\cdot(sy)|\leq \frac{1}{2}$ and hence
\[
1-\cos\left(\varpi\cdot(sy)\right)\geq C_1|\varpi\cdot(sy)|^2.
\]
Therefore, for sufficiently small $s>0$,
\begin{align*}
H_2(\varpi, s)&\geq s^d\sum_{y\in \mathbb{Z}^d\setminus\{\mathbf{0}\}: \|y\|_2\leq \frac{1}{2\sqrt{d}s}}\frac{1-\cos\left(\varpi\cdot(sy)\right)}{\|sy\|_2^{d+2}}\\
&\geq C_1s^d\sum_{y\in \mathbb{Z}^d\setminus\{\mathbf{0}\}: \|y\|_2\leq \frac{1}{2\sqrt{d}s}}\frac{|\varpi\cdot(sy)|^2}{\|sy\|_2^{d+2}}\\
&=C_1\sum_{y\in \mathbb{Z}^d\setminus\{\mathbf{0}\}: \|y\|_2\leq \frac{1}{2\sqrt{d}s}}\frac{|\varpi\cdot y|^2}{\|y\|_2^{d+2}}\\
&=C_1\Theta(1)\int_1^{\frac{1}{2\sqrt{d}s}}\frac{r^{d-1}r^2}{r^{d+2}}dr\int_{v:\|v\|_2=1}\frac{|\varsigma\cdot v|^2}{\|v\|_2^{d+2}}dS\\
&=\Theta(1)\log (s^{-1}),
\end{align*}
where $\varsigma$ can be any element in $\mathbb{R}^d$ such that $\|\varsigma\|_2=1$. Hence, when $\frac{M}{\sqrt{N\log N}}\leq s\leq (N\log N)^{-\frac{1}{4}}$,
\[
H_2(\varpi, s)\geq \Theta(1)\log N
\]
and then, by \eqref{equ 3.5 two},
\begin{align}\label{equ 3.6}
&\int_{\theta: M<\|\theta\|_\infty\leq (N\log N)^{1/4}}\psi_{tN}^{d, 2}\left(\theta\big/\sqrt{N\log N}\right)d\theta\notag\\
&\leq \int_{\theta: M<\|\theta\|_\infty\leq (N\log N)^{1/4}}e^{-t\Theta(1)\|\theta\|_\infty^2}d\theta.
\end{align}
When $(N\log N)^{-\frac{1}{4}}<s\leq \epsilon$, $H_2(\varpi, s)\geq \Theta(1)\log (\epsilon^{-1})\geq \Theta(1)$. Hence, by \eqref{equ 3.5 two},
\[
\psi_{tN}^{d, 2}\left(\theta\big/\sqrt{N\log N}\right)\leq \exp\left\{-t\Theta(1)\frac{\sqrt{N}}{\sqrt{\log N}}\right\}
\]
when $(N\log N)^{1/4}<\|\theta\|_\infty\leq \epsilon\sqrt{N\log N}$. Therefore,
\begin{align}\label{equ 3.8}
&\int_{\theta:(N\log N)^{1/4}<\|\theta\|_\infty\leq \epsilon\sqrt{N\log N}}\psi_{tN}^{d, 2}\left(\theta\big/\sqrt{N\log N}\right)d\theta\notag\\
&\leq \int_{\theta: (N\log N)^{1/4}<\|\theta\|_\infty\leq \epsilon\sqrt{N\log N}}\exp\left\{-t\Theta(1)\frac{\sqrt{N}}{\sqrt{\log N}}\right\}d\theta.
\end{align}
Equation \eqref{equ 3.5} follows from \eqref{equ 3.6} and \eqref{equ 3.8}. Since \eqref{equ 3.5} holds, the proof is complete.
\qed

\section{Proof of Theorem \ref{theorem 2.1 main result d=1, 2}}\label{section four}
In this section, we prove Theorem \ref{theorem 2.1 main result d=1, 2}. As we have introduced in Section \ref{section two}, the strategy of the proof is similar to that of the nearest-neighbor case given in \cite{Xue2024}. Hence, to avoid repeating similar details, we will emphasize the main difference between the arguments in this paper and \cite{Xue2024} and only give an outline of which has been explained precisely in \cite{Xue2024} for the nearest-neighbor case and can be directly extended to the long-range case.
\subsection{Proof of Theorem \ref{theorem 2.1 main result d=1, 2}: the case of $\alpha\in \left(d/2, \min\{2, d\}\right)$}\label{subsection 4.1}
We first deal with the case of $\alpha\in \left(d/2, \min\{2, d\}\right)$.
For later use, we introduce some notations and definitions. For any $x, y\in \mathbb{Z}^d$ such that $x \neq y$, let $\{N_t^{x, y}\}_{t\geq 0}$ be the Poisson process with rate $\|x-y\|_2^{-(d+\alpha)}$ such that $x$ adopts the opinion of $y$ at every event moments of $\{N_t^{x, y}\}_{t\geq 0}$. Then, all these Poisson processes are independent and
\begin{equation}\label{equ 4.1}
\eta_t(x)=\eta_0(x)+\sum_{y\in \mathbb{Z}^d, y\neq x}\int_0^t\left(\eta_{s-}(y)-\eta_{s-}(x)\right)dN_s^{x, y}
\end{equation}
for any $t\geq 0$, since
\[
\eta_s(x)=\eta_{s-}(y)=\eta_{s-}(x)+\left(\eta_{s-}(y)-\eta_{s-}(x)\right)
\]
at each event moment $s$ of $\{N^{x, y}_t\}_{t\geq 0}$. For any $t\geq 0$, we define
\begin{equation}\label{equ centered Poisson}
\hat{N}_t^{x, y}=N_t^{x, y}-\|x-y\|_2^{-(d+\alpha)}t,
\end{equation}
then $\{\hat{N}^{x, y}_t\}_{t\geq 0}$ is a martingale with quadratic variation process
\[
\langle \hat{N}^{x, y}\rangle_t=\|x-y\|_2^{-(d+\alpha)}t.
\]
Let $\{p_t^{d, \alpha}\}_{t\geq 0}$ be the transition probabilities of $\{X_t^{d, \alpha}\}_{t\geq 0}$ defined as in Subsection \ref{subsection 1.2}. We define
\[
v(t, x)=\int_0^tp_s^{d, \alpha}(\mathbf{0}, x)ds
\]
for any $x\in \mathbb{Z}^d$. For $0\leq s\leq t$, we define
\[
V^t(s, \eta)=\sum_{x\in \mathbb{Z}^d}\left(\eta(x)-p\right)v(t-s, x).
\]
and
\[
M_s^t=V^t(s, \eta_s)-V^t(0, \eta_0)-\int_0^s\left(\partial_r+\mathcal{G}\right)V^t(r, \eta_r)dr,
\]
where $\mathcal{G}$ is the generator of the long-range voter model given by \eqref{equ 1.1 generator of long-range voter model}. According to the Dynkin's martingale formula, $\{M_s^t\}_{0\leq s\leq t}$ is a martingale. According to the expression of $\mathcal{G}$ and Kolmogorov-Chapman equation, it is not difficult to check that
\begin{equation}\label{equ 4.2}
\left(\partial_r+\mathcal{G}\right)V^t(r, \eta_r)=-\left(\eta_r(\mathbf{0})-p\right).
\end{equation}
References \cite{Kipnis1987}, \cite{Birkner2007} and \cite{Xue2024} give precise check of analogues of \eqref{equ 4.2} for exclusion processes, critical branching random walks and nearest-neighbor voter models, which can be extended to the long-range voter model directly. So here we omit the detailed proof of \eqref{equ 4.2}.

Since $V^t(t, \eta)=0$, by \eqref{equ 4.2},
\begin{equation}\label{equ 4.3 martingale decomposition}
\int_0^{tN}\left(\eta_s(\mathbf{0})-p\right)ds=M_{tN}^{tN}+V^{tN}(0, \eta_0).
\end{equation}
We first deal with the term $V^{tN}(0, \eta_0)$ in decomposition \eqref{equ 4.3 martingale decomposition}. We have the following lemma.

\begin{lemma}\label{lemma 4.1}
Let
\begin{equation}\label{equ 4.4}
(d, \alpha)\in \{1\}\times \left(\frac{1}{2}, 1\right)\bigcup \{2\}\times \left(1, 2\right)\bigcup \{3\}\times \left(\frac{3}{2}, 2\right).
\end{equation}
For any integer $m\geq 1$ and $0<t_1<t_2<\ldots<t_m$, if $\eta_0$ is distributed with $\nu_p^{d, \alpha}$, then
\[
\frac{1}{N^{\frac{3}{2}-\frac{d}{2\alpha}}}\left(V^{t_1N}(0, \eta_0), V^{t_2N}(0, \eta_0), \ldots, V^{t_mN}(0, \eta_0)\right)
\]
converges weakly to a $\mathbb{R}^m$-valued Gaussian random variable  $\left(\Gamma_{t_1}, \ldots, \Gamma_{t_m}\right)$ as $N\rightarrow+\infty$.
\end{lemma}

To prove Lemma \ref{lemma 4.1}, we first give a central limit theorem of the fluctuation density field of the long-range voter model. In detail, for $(d, \alpha)$ satisfying \eqref{equ 4.4}, we define
\[
\mathcal{Y}_t^\epsilon(du)=\epsilon^{\frac{d+\alpha}{2}}\sum_{x\in \mathbb{Z}^d}\left(\eta_{t\epsilon^{-\alpha}}(x)-p\right)\delta_{\epsilon x}(du)
\]
for any $\epsilon>0$, where $\delta_a(du)$ is the Dirac measure concentrated on $a$. That is to say, for any Schwartz function $f$ on $\mathbb{R}^d$,
\[
\mathcal{Y}_t^\epsilon(f)=\epsilon^{\frac{d+\alpha}{2}}\sum_{x\in \mathbb{Z}^d}\left(\eta_{t\epsilon^{-\alpha}}(x)-p\right)f(\epsilon x).
\]
We call $\mathcal{Y}_t^\epsilon$ the fluctuation density field of the voter model. To give the limit theorem of $\mathcal{Y}_t^{\epsilon}$ as $\epsilon\rightarrow 0$, we introduce the following generalized Ornstein-Uhlenbeck process
\begin{equation}\label{equ 4.5 generalized O-U process}
d\mathcal{Y}_t=\hat{\mathcal{L}}^{d, \alpha}\mathcal{Y}_tdt+Cd\mathcal{B}_t,
\end{equation}
where $\hat{\mathcal{L}}^{d, \alpha}$ is the generator of the stable process $\{\hat{X}_t^{d, \alpha}\}_{t\geq 0}$ defined as in Subsection \ref{subsection 1.2}, $\{\mathcal{B}_t\}_{t\geq 0}$ is the time-space white noise on $\mathbb{R}^d$ and $C\neq 0$ is a real constant. Note that the time-space white noise $\{\mathcal{B}_t\}_{t\geq 0}$ satisfies that
$\{\mathcal{B}_t(f)\}_{t\geq 0}$ is a copy of
\[
\left\{\sqrt{\int_{\mathbb{R}^d}f^2(u)du}\mathcal{W}_t\right\}_{t\geq 0}
\]
for any Schwartz function $f$ on $\mathbb{R}^d$, where $\{\mathcal{W}_t\}_{t\geq 0}$ is the real-valued standard Brownian motion starting from $0$ defined as in Section \ref{section two}. The rigorous definition of the solution $\{\mathcal{Y}_t\}_{t\geq 0}$ to \eqref{equ 4.5 generalized O-U process} is as follows. For any Schwartz function $f$ on $\mathbb{R}^d$ and $G\in C_c^\infty(\mathbb{R})$,
\begin{align*}
&\Bigg\{G(\mathcal{Y}_t(f))-G(\mathcal{Y}_0(f))-\int_0^tG^\prime(\mathcal{Y}_s(f))\mathcal{Y}_s\left(\hat{\mathcal{L}}^{d, \alpha}f\right)ds\\
&\text{\quad\quad\quad}-\frac{C^2}{2}\int_{\mathbb{R}^d}f^2(u)du\int_0^tG^{\prime\prime}(\mathcal{Y}_s(f))ds\Bigg\}_{t\geq 0}
\end{align*}
is a martingale. By solving Equation \eqref{equ 4.5 generalized O-U process} via the integrating factor method, an equivalent definition of $\{\mathcal{Y}_t\}_{t\geq 0}$ is given by
\begin{equation}\label{equ 4.6}
\mathcal{Y}_t=\mathcal{P}_t\mathcal{Y}_0+C\int_0^t\mathcal{P}_{t-s}d\mathcal{B}_s,
\end{equation}
where $\{\mathcal{P}_t\}_{t\geq 0}$ is the Markov semigroup with respect to $\hat{\mathcal{L}}^{d, \alpha}$, i.e.,
\[
\mathcal{P}_tf(u)=\mathbb{E}\left(f\left(\hat{X}_t^{d,\alpha}\right)\big|\hat{X}_0^{d, \alpha}=u\right)
\]
and
\[
(\mathcal{P}_t\mu)(f)=\mu\left(\mathcal{P}_tf\right)
\]
for any $u\in \mathbb{R}^d$, Schwartz function $f$ and tempered distribution $\mu$. By \eqref{equ 4.6}, since a generator is a negative operator, $\mathcal{Y}_t$ converges weakly, as $t\rightarrow +\infty$, to a tempered-distribution-valued Gaussian random variable $\mathcal{Y}_\infty$ such that $\mathcal{Y}_\infty(f)$ follows the normal distribution with mean $0$ and variance
\[
C^2\int_0^{+\infty}\int_{\mathbb{R}^d}\left(\mathcal{P}_sf(u)\right)^2dsdu
\]
for any Schwartz function $f$, no matter what is the distribution of $\mathcal{Y}_0$. Hence, Equation \eqref{equ 4.5 generalized O-U process} has a unique stationary distribution $\pi^C$ such that $\mu(f)$ follows the normal distribution with mean $0$ and variance $C^2\int_0^{+\infty}\int_{\mathbb{R}^d}\left(\mathcal{P}_sf(u)\right)^2duds$ for any Schwartz function $f$ when the random tempered distribution $\mu$ is distributed with $\pi^C$.

Let $\Phi^{d,\alpha}(x)$ be defined as in \eqref{equ 1.12}, now we give a lemma of the central limit theorem of $\{\mathcal{Y}_t^\epsilon\}_{t\geq 0}$.

\begin{lemma}\label{lemma 4.2 fluctuation of long-range voter model}
For $(d, \alpha)$ satisfying \eqref{equ 4.4} and any $T>0$, as $\epsilon\rightarrow 0$,
\[
\left\{\mathcal{Y}_t^\epsilon:~0\leq t\leq T\right\}
\]
converges weakly, with respect to the Skorohod topology, to $\left\{\mathcal{Y}_t:~0\leq t\leq T\right\}$, where $\left\{\mathcal{Y}_t\right\}_{t\geq 0}$ is the solution to Equation \eqref{equ 4.5 generalized O-U process} with
\[
C=\sqrt{2p(1-p)\sum_{y\neq \mathbf{0}}\frac{\Phi^{d, \alpha}(y)}{\|y\|^{d+\alpha}}}
\]
and $\mathcal{Y}_0$ being distributed with $\pi^C$.
\end{lemma}

Lemma \ref{lemma 4.2 fluctuation of long-range voter model} is an analogue of Theorem 2 of \cite{Presutti1983}, which gives the stationary fluctuation limit of the empirical density field of the nearest-neighbor voter model. We prove Lemma \ref{lemma 4.2 fluctuation of long-range voter model} later. Now we utilize Lemma \ref{lemma 4.2 fluctuation of long-range voter model} to prove Lemma \ref{lemma 4.1}.

\proof
[Proof of Lemma \ref{lemma 4.1}]
For any $N\geq 1$, $u\in \mathbb{R}^d$ and $\alpha>0$, we denote by $u_{\alpha, N}$ the element in $\mathbb{Z}^d/(N^{\frac{1}{\alpha}})$ such that
\[
u-u_{\alpha, N}\in \left(-\frac{1}{2N^{\frac{1}{\alpha}}}, \frac{1}{2N^{\frac{1}{\alpha}}}\right]^d.
\]
For all $N\geq 1$ and $u\in \mathbb{R}^d$, we define
\[
\tilde{J}^{N}_t(u)=\int_0^tN^{\frac{d}{\alpha}}p_{sN}^{d, \alpha}(\mathbf{0}, N^{\frac{1}{\alpha}}u_{\alpha, N})ds
\]
and
\[
J_t(u)=\int_0^tf^{d, \alpha}_s(u)ds,
\]
where $f^{\alpha, d}_s$ is defined as in Subsection \ref{subsection 1.2}. By Lemma \ref{lemma LCLT},
\[
\lim_{N\rightarrow+\infty}\tilde{J}^N_t(u)=J_t(u)
\]
for any $u\in \mathbb{R}^d$. Furthermore, since
\[
\int_{\mathbb{R}^d}\tilde{J}^{N}_t(u)du=t=\int_{\mathbb{R}^d}J_t(u)du,
\]
we have
\begin{equation}\label{equ 4.7}
\lim_{N\rightarrow+\infty}\int_{\mathbb{R}^d}\left|\tilde{J}^N_t(u)-J_t(u)\right|du=0.
\end{equation}
We claim that, for $(d, \alpha)$ satisfies \eqref{equ 4.4},
\begin{equation}\label{equ 4.8}
\lim_{N\rightarrow+\infty}\int_{\mathbb{R}^d}\left(\tilde{J}^N_t(u)\right)^2du=\int_{\mathbb{R}^d}J_t^2(u)du<+\infty.
\end{equation}
We check \eqref{equ 4.8} in Appendix \ref{Appendix A.1}. By \eqref{equ 4.8}, we further have
\begin{equation}\label{equ 4.9}
\lim_{N\rightarrow+\infty}\int_{\mathbb{R}^d}\left|\tilde{J}^N_t(u)-J_t(u)\right|^2du=0.
\end{equation}
According to the definition of $V^t(s, \eta)$ and $\mathcal{Y}_t^\epsilon$,
\begin{align}\label{equ 4.9 two}
\frac{1}{N^{\frac{3}{2}-\frac{d}{2\alpha}}}V^{tN}(0, \eta_0)&=\frac{1}{N^{\frac{1}{2}+\frac{d}{2\alpha}}}\sum_{x\in \mathbb{Z}^d}\left(\eta_0(x)-p\right)\int_0^tN^{\frac{d}{\alpha}}p_{sN}^{d, \alpha}(0,x)ds\notag\\
&=\mathcal{Y}_0^{N^{-\frac{1}{\alpha}}}(\tilde{J}_t^N).
\end{align}
By Lemma \ref{lemma 4.2 fluctuation of long-range voter model}, as $N\rightarrow+\infty$,
\[
\left(\mathcal{Y}_0^{N^{-\frac{1}{\alpha}}}(J_{t_1}), \mathcal{Y}_0^{N^{-\frac{1}{\alpha}}}(J_{t_2}), \ldots, \mathcal{Y}_0^{N^{-\frac{1}{\alpha}}}(J_{t_m})\right)
\]
converges weakly to
\[
\left(\mathcal{Y}_0(J_{t_1}), \mathcal{Y}_0(J_{t_2}), \ldots, \mathcal{Y}_0(J_{t_m})\right),
\]
where $\mathcal{Y}_0$ is distributed with $\pi^C$ and hence $\left(\mathcal{Y}_0(J_{t_1}), \mathcal{Y}_0(J_{t_2}), \ldots, \mathcal{Y}_0(J_{t_m})\right)$
is a $\mathbb{R}^m$-valued Gaussian random variable such that
\[
{\rm Cov}\left(\mathcal{Y}_0(J_{t_i}), \mathcal{Y}_0(J_{t_j})\right)
=\left(2p(1-p)\sum_{y\neq \mathbf{0}}\frac{\Phi^{d, \alpha}(y)}{\|y\|^{d+\alpha}}\right)\int_0^{+\infty}\int_{\mathbb{R}^d}\mathcal{P}_sJ_{t_i}(u)\mathcal{P}_sJ_{t_j}(u)dsdu
\]
for any $1\leq i, j\leq m$.
Then by \eqref{equ 4.9 two}, to complete this proof, we only need to show that $\mathcal{Y}_0^{N^{-\frac{1}{\alpha}}}\left(J_t-\tilde{J}_t^N\right)$ converges weakly to $0$ as $N\rightarrow+\infty$.

For any Schwartz function $f$, by \eqref{equ 1.14} and \eqref{equ 1.15},
\begin{align}\label{equ 4.9 three}
&\hat{\mathbb{E}}_{\nu_p^{d, \alpha}}\left(\left(\mathcal{Y}_0^{N^{-\frac{1}{\alpha}}}(f)\right)^2\right)\notag\\
&=N^{-1-\frac{d}{\alpha}}\sum_{x\in \mathbb{Z}^d}\sum_{y\in \mathbb{Z}^d}{\rm Cov}_{\nu_p^{d,\alpha}}\left(\eta_0(x), \eta_0(y)\right)f\left(N^{-\frac{1}{\alpha}}x\right)f\left(N^{-\frac{1}{\alpha}}y\right)\notag\\
&=p(1-p)N^{-1-\frac{d}{\alpha}}\sum_{x\in \mathbb{Z}^d}\sum_{y\in \mathbb{Z}^d}\left(1-\Phi(x-y)\right)f\left(N^{-\frac{1}{\alpha}}x\right)f\left(N^{-\frac{1}{\alpha}}y\right).
\end{align}
According to the Strong Markov property of the random walk,
\begin{align}\label{equ 4.9 four}
1-\Phi(x-y)&=\mathcal{C}_{d, \alpha}\int_0^{+\infty}p_s^{d, \alpha}(\mathbf{0}, x-y)ds,
\end{align}
where $\mathcal{C}_{d, \alpha}$ is defined as in Section \ref{section two}. For $(d, \alpha)$ satisfying \eqref{equ 4.4}, we claim that there exists a constant $C_2\in (0, +\infty)$ independent of $x\in \mathbb{Z}^d$ such that
\begin{equation}\label{equ 4.10}
\lim_{\|x\|_2\rightarrow+\infty}\frac{\int_0^{+\infty}p_s^{d, \alpha}(\mathbf{0}, x)ds}{\|x\|_2^{\alpha-d}}=C_2.
\end{equation}
We check \eqref{equ 4.10} in Appendix \ref{Appendix A.2}. By \eqref{equ 4.10},
\begin{align*}
&N^{-1-\frac{d}{\alpha}}\sum_{x\in \mathbb{Z}^d}\sum_{y\in \mathbb{Z}^d}\left(1-\Phi(x-y)\right)f\left(N^{-\frac{1}{\alpha}}x\right)f\left(N^{-\frac{1}{\alpha}}y\right)\\
&=O(1)N^{-1-\frac{d}{\alpha}}N^{1-\frac{d}{\alpha}}\sum_{x\in \mathbb{Z}^d}\sum_{y\in \mathbb{Z}^d, y\neq x}\left\|\frac{x}{N^{\frac{1}{\alpha}}}-\frac{y}{N^{\frac{1}{\alpha}}}\right\|_2^{\alpha-d}
f\left(N^{-\frac{1}{\alpha}}x\right)f\left(N^{-\frac{1}{\alpha}}y\right)\\
&\text{\quad\quad}+O(1)N^{-1-\frac{d}{\alpha}}\sum_{x\in \mathbb{Z}^d}f^2\left(N^{-\frac{1}{\alpha}}x\right)\\
&=O(1)\int_{u\in \mathbb{R}^d}\int_{v\in \mathbb{R}^d}\|u-v\|_2^{\alpha-d}f(u)f(v)dudv+N^{-1}O(1)\int_{\mathbb{R}^d}f^2(u)du.
\end{align*}
Since $\alpha-d<0$,
\[
\int\int_{u,v\in \mathbb{R}^d:\|u-v\|_2\geq 1}\|u-v\|_2^{\alpha-d}f(u)f(v)dudv\leq \left(\int_{\mathbb{R}^d}|f(u)|du\right)^2.
\]
By Cauchy-Schwarz inequality,
\begin{align*}
&\int\int_{u,v\in \mathbb{R}^d:\|u-v\|_2\leq 1}\|u-v\|_2^{\alpha-d}f(u)f(v)dudv\\
&=\int_{\theta\in \mathbb{R}^d:\|\theta\|_2\leq 1}\|\theta\|_2^{\alpha-d}\left(\int_{\mathbb{R}^d}f(v)f(v+\theta)dv\right)d\theta\\
&\leq \int_{\theta\in \mathbb{R}^d:\|\theta\|_2\leq 1}\|\theta\|_2^{\alpha-d}d\theta\int_{\mathbb{R}^d}f^2(u)du.
\end{align*}
In conclusion, for sufficiently large $N$,
\begin{align}\label{equ 4.12}
&N^{-1-\frac{d}{\alpha}}\sum_{x\in \mathbb{Z}^d}\sum_{y\in \mathbb{Z}^d}\left(1-\Phi(x-y)\right)f\left(N^{-\frac{1}{\alpha}}x\right)f\left(N^{-\frac{1}{\alpha}}y\right)\notag\\
&\leq O(1)\left(\left(\int_{\mathbb{R}^d}|f(u)|du\right)^2+\int_{\mathbb{R}^d}f^2(u)du\right).
\end{align}
By \eqref{equ 4.9 three} and \eqref{equ 4.12},
\begin{equation}\label{equ 4.12 second}
\hat{\mathbb{E}}_{\nu_p^{d, \alpha}}\left(\left(\mathcal{Y}_0^{N^{-\frac{1}{\alpha}}}(f)\right)^2\right)
=O(1)\left(\left(\int_{\mathbb{R}^d}|f(u)|du\right)^2+\int_{\mathbb{R}^d}f^2(u)du\right).
\end{equation}
By \eqref{equ 4.7}, \eqref{equ 4.9}, \eqref{equ 4.12 second} and Chebyshev inequality, $\mathcal{Y}_0^{N^{-\frac{1}{\alpha}}}\left(J_t-\tilde{J}_t^N\right)$
 converges weakly to $0$ as $N\rightarrow+\infty$ and the proof is complete. \qed

Now we prove Lemma \ref{lemma 4.2 fluctuation of long-range voter model} to complete the proof of Lemma \ref{lemma 4.1}. The proof of Lemma \ref{lemma 4.2 fluctuation of long-range voter model} follows from the classic Dynkin's martingale strategy and hence we only highlight several crucial steps.

\proof[Proof of Lemma \ref{lemma 4.2 fluctuation of long-range voter model}]
For any Schwartz function $f$ on $\mathbb{R}^d$ and $G\in C_c^\infty(\mathbb{R})$, by Dynkin's martingale formula,
\[
\left\{G(\mathcal{Y}_t^\epsilon(f))-G(\mathcal{Y}_0^\epsilon(f))-\int_0^t\epsilon^{-\alpha}\mathcal{G}G(\mathcal{Y}_s^\epsilon(f))ds\right\}_{t\geq 0}
\]
is a martingale. According to the definition of $\mathcal{G}$ given in \eqref{equ 1.1 generator of long-range voter model},
\begin{align*}
&\mathcal{G}G(\mathcal{Y}_s^\epsilon(f))=\\
&\sum_{x\in \mathbb{Z}^d}\sum_{y\in \mathbb{Z}^d\setminus\{\mathbf{0}\}}\|y\|_2^{-(d+\alpha)}
\left(G\left(\mathcal{Y}_s^\epsilon(f)+\epsilon^{\frac{d+\alpha}{2}}\left(\eta_{s\epsilon^{-\alpha}}(y+x)-\eta_{s\epsilon^{-\alpha}}(x)\right)f(\epsilon x)\right)-G(\mathcal{Y}_s^\epsilon(f))\right).
\end{align*}
By Taylor's expansion formula up to the second order,
\begin{align*}
&G\left(\mathcal{Y}_s^\epsilon(f)+\epsilon^{\frac{d+\alpha}{2}}\left(\eta_{s\epsilon^{-\alpha}}(y+x)-\eta_{s\epsilon^{-\alpha}}(x)\right)f(\epsilon x)\right)-G(\mathcal{Y}_s^\epsilon(f))\\
&=G^\prime(\mathcal{Y}_s^\epsilon(f))\epsilon^{\frac{d+\alpha}{2}}\left(\eta_{s\epsilon^{-\alpha}}(y+x)-\eta_{s\epsilon^{-\alpha}}(x)\right)f(\epsilon x)\\
&\text{\quad}+\frac{G^{\prime\prime}(\mathcal{Y}_s^\epsilon(f))}{2}\epsilon^{d+\alpha}
\left(\eta_{s\epsilon^{-\alpha}}(y+x)-\eta_{s\epsilon^{-\alpha}}(x)\right)^2f^2(\epsilon x)+o(\epsilon^{d+\alpha})f^3(\epsilon x).
\end{align*}
Note that
\begin{align*}
&\epsilon^{-\alpha}\sum_{x\in \mathbb{Z}^d}\sum_{y\in\mathbb{Z}^d\setminus\{\mathbf{0}\}}\|y\|_2^{-(d+\alpha)}\epsilon^{\frac{d+\alpha}{2}}\left(\eta_{s\epsilon^{-\alpha}}(y+x)-\eta_{s\epsilon^{-\alpha}}(x)\right)f(\epsilon x)\\
&=\epsilon^{-\alpha}\sum_{x\in \mathbb{Z}^d}\sum_{y\in \mathbb{Z}^d\setminus\{\mathbf{0}\}}\|y\|_2^{-(d+\alpha)}\epsilon^{\frac{d+\alpha}{2}}\left((\eta_{s\epsilon^{-\alpha}}(y+x)-p)-(\eta_{s\epsilon^{-\alpha}}(x)-p)\right)f(\epsilon x)\\
&=\epsilon^{\frac{d+\alpha}{2}}\sum_{x\in \mathbb{Z}^d}(\eta_{s\epsilon^{-\alpha}}(x)-p)\left(\tilde{\mathcal{L}}^\epsilon f(\epsilon x)\right)=\mathcal{Y}_s^\epsilon(\tilde{\mathcal{L}}^\epsilon f),
\end{align*}
where
\[
\tilde{\mathcal{L}}^\epsilon f(u)=\epsilon^{-\alpha}\sum_{y\in\mathbb{Z}^d\setminus\{\mathbf{0}\}}\|y\|_2^{-(d+\alpha)}\left(f(u+\epsilon y)-f(u)\right).
\]
According to the definition of Riemann integral, $\tilde{\mathcal{L}}^\epsilon f(u)\rightarrow \hat{\mathcal{L}}^{d, \alpha}f(u)$ as $\epsilon\rightarrow 0$. Furthermore, since
\[
\lim_{\|u\|_2\rightarrow+\infty}\|u\|_2^m|f(u)|=0
\]
and
\[
\lim_{\|u\|_2\rightarrow+\infty}\|u\|_2^m\left|\partial_{u_iu_j}^2f(u)\right|=0
\]
for any $m\geq 1$ and $1\leq i, j\leq d$, it is easy to check that
\[
\lim_{\epsilon\rightarrow 0}\int_{\mathbb{R}^d}\left|\tilde{\mathcal{L}}^\epsilon f(u)-\hat{\mathcal{L}}^{d, \alpha}f(u)\right|^2du=0
\]
and
\[
\lim_{\epsilon\rightarrow 0}\int_{\mathbb{R}^d}\left|\tilde{\mathcal{L}}^\epsilon f(u)-\hat{\mathcal{L}}^{d, \alpha}f(u)\right|du=0
\]
according to a routine truncation argument.
Then, by \eqref{equ 4.12 second} and the invariance of $\nu_p^{d, \alpha}$, $\mathcal{Y}_s^\epsilon(\tilde{\mathcal{L}}^\epsilon f)-\mathcal{Y}_s^\epsilon(\hat{\mathcal{L}}^{d, \alpha}f)$ converges to $0$ in $L^2$ as $\epsilon\rightarrow 0$. For any $T>0$, we claim that
\[
\left\{\mathcal{Y}_t^\epsilon:~0\leq t\leq T\right\}_{\epsilon>0}
\]
are tight under the Skorohod topology. By \eqref{equ 4.12 second}, the above tightness follows from an argument utilizing the Aldous' criterion similar to that given in \cite{Presutti1983} for the nearest-neighbor case, the detail of which we omit here.

By \eqref{equ 1.13} and the invariance of $\nu_p^{d, \alpha}$,
\[
\hat{\mathbb{E}}_{\nu_p^{d, \alpha}}\left(\left(\eta_{s\epsilon^{-\alpha}}(y+x)-\eta_{s\epsilon^{-\alpha}}(x)\right)^2\right)=2p(1-p)\Phi^{d, \alpha}(y)
\]
and hence
\begin{align}\label{equ 4.15}
&\lim_{\epsilon\rightarrow 0}\epsilon^{-\alpha}\sum_{x\in \mathbb{Z}^d}\sum_{y\in \mathbb{Z}^d\setminus\{\mathbf{0}\}}\|y\|_2^{-(d+\alpha)}\epsilon^{d+\alpha}
\mathbb{E}_{\nu_p^{d, \alpha}}\left(\left(\eta_{s\epsilon^{-\alpha}}(y+x)-\eta_{s\epsilon^{-\alpha}}(x)\right)^2\right)f^2(\epsilon x)\notag\\
&=C^2\int_{\mathbb{R}^d}f^2(u)du,
\end{align}
where
\[
C=\sqrt{2p(1-p)\sum_{y\neq \mathbf{0}}\frac{\Phi^{d, \alpha}(y)}{\|y\|^{d+\alpha}}}.
\]
We claim that
\begin{equation}\label{equ 4.16}
\lim_{\epsilon\rightarrow 0}{\rm Var}_{\nu_p^{d, \alpha}}\left(\epsilon^d\sum_{x\in \mathbb{Z}^d}\sum_{y\in \mathbb{Z}^d\setminus\{\mathbf{0}\}}\|y\|_2^{-(d+\alpha)}\left(\eta_{s\epsilon^{-\alpha}}(y+x)-\eta_{s\epsilon^{-\alpha}}(x)\right)^2f^2(\epsilon x)\right)=0.
\end{equation}
We check \eqref{equ 4.16} in Appendix \ref{Appendix A.3}. By \eqref{equ 4.15} and \eqref{equ 4.16},
\begin{align*}
&\lim_{\epsilon\rightarrow 0}\epsilon^d\sum_{x\in \mathbb{Z}^d}\sum_{y\in \mathbb{Z}^d\setminus\{\mathbf{0}\}}\|y\|_2^{-(d+\alpha)}\left(\eta_{s\epsilon^{-\alpha}}(y+x)-\eta_{s\epsilon^{-\alpha}}(x)\right)^2f^2(\epsilon x)\\
&=C^2\int_{\mathbb{R}^d}f^2(u)du \text{\quad in~}L^2.
\end{align*}
Consequently, for any weak limit $\left\{\tilde{\mathcal{Y}}_t:~0\leq t\leq T\right\}$ of a subsequence of $\left\{\mathcal{Y}_t^\epsilon:~0\leq t\leq T\right\}_{\epsilon>0}$,
\begin{align*}
&\Bigg\{G(\tilde{\mathcal{Y}}_t(f))-G(\tilde{\mathcal{Y}}_0(f))-\int_0^tG^\prime(\tilde{\mathcal{Y}}_s(f))\tilde{\mathcal{Y}}_s\left(\hat{\mathcal{L}}^{d, \alpha}f\right)ds\\
&\text{\quad\quad\quad}-\frac{C^2}{2}\int_{\mathbb{R}^d}f^2(u)du\int_0^tG^{\prime\prime}(\tilde{\mathcal{Y}}_s(f))ds\Bigg\}_{0\leq t\leq T}
\end{align*}
is a weak limit of a subsequence of
\[
\left\{G(\mathcal{Y}_t^\epsilon(f))-G(\mathcal{Y}_0^\epsilon(f))-\int_0^t\epsilon^{-\alpha}\mathcal{G}G(\mathcal{Y}_s^\epsilon(f))ds:~0\leq t\leq T\right\}_{\epsilon>0}
\]
and hence is a martingale. Note that the existence of $\left\{\tilde{\mathcal{Y}}_t:~0\leq t\leq T\right\}$ follows from the aforesaid tightness.  As a result,
$\left\{\tilde{\mathcal{Y}}_t:~0\leq t\leq T\right\}$ is a solution to Equation \eqref{equ 4.5 generalized O-U process}. According to the invariance of $\nu_p^{d, \alpha}$, the distribution of $\tilde{\mathcal{Y}}_t$ is independent of $t$. Hence, $\tilde{\mathcal{Y}}_0$ is distributed with a stationary distribution of \eqref{equ 4.5 generalized O-U process}. As we have explained, the stationary distribution of \eqref{equ 4.5 generalized O-U process} is unique, i.e., is $\pi^C$, which completes the proof.
\qed

\quad

Now we deal with the martingale term $M_{tN}^{tN}$ in decomposition \eqref{equ 4.3 martingale decomposition}. We have the following lemma.

\begin{lemma}\label{lemma 4.3}
Let $(d, \alpha)$ satisfy \eqref{equ 4.4}. For any integer $m\geq 1$ and $0<t_1<t_2<\ldots<t_m$, if $\eta_0$ is distributed with $\nu_p^{d, \alpha}$, then
\[
\frac{1}{N^{\frac{3}{2}-\frac{d}{2\alpha}}}\left(M_{t_1N}^{t_1N}, \ldots, M_{t_mN}^{t_mN}\right)
\]
converges weakly, as $N\rightarrow+\infty$, to a $\mathbb{R}^m$-valued Gaussian random variable $(\Xi_{t_1},\ldots, \Xi_{t_m})$.
\end{lemma}
As a preliminary of the proof of Lemma \ref{lemma 4.3}, we introduce some notations. As we have defined in the proof of Lemma \ref{lemma 4.1}, for any $u\in \mathbb{R}^d$, we denote by $u_{\alpha, N}$ the element in $\mathbb{Z}^d/(N^{\frac{1}{\alpha}})$ such that
\[
u-u_{\alpha, N}\in \left(-\frac{1}{2N^{\frac{1}{\alpha}}}, \frac{1}{2N^{\frac{1}{\alpha}}}\right]^d.
\]
For $T>0$ and $N\geq 1$, we denote by $\mathcal{U}^N$ the random measure on $[0, T]\times \mathbb{R}^d$ such that
\begin{align*}
&\mathcal{U}^N(H)= \\
&N^{\frac{d}{2\alpha}-\frac{1}{2}}\int_{\mathbb{R}^d}\sum_{y\in \mathbb{Z}^d\setminus\{N^{\frac{1}{\alpha}}u_{\alpha, N}\}}\left(\int_0^T H(s, u)\left(\eta_{sN-}(y)-\eta_{sN-}(N^{\frac{1}{\alpha}}u_{\alpha, N})\right)d\hat{N}_{Ns}^{(N^{\frac{1}{\alpha}}u_{\alpha, N}), y}\right)du\\
&=N^{\frac{d}{2\alpha}-\frac{1}{2}}\sum_{x\in \mathbb{Z}^d}\sum_{y\in \mathbb{Z}^d\setminus\{x\}}\int_0^T\left(\int_{\frac{x}{N^{\frac{1}{\alpha}}}+\left(-\frac{1}{2N^{\frac{1}{\alpha}}}, \frac{1}{2N^{\frac{1}{\alpha}}}\right]^d}H(s, u)du\right)(\eta_{sN-}(y)-\eta_{sN-}(x))d\hat{N}_{Ns}^{x, y}
\end{align*}
for any $H\in C_c([0, T]\times \mathbb{R}^d)$, where $\hat{N}_s^{x, y}$ is defined as in \eqref{equ centered Poisson}. For each $N\geq 1$ and $t\leq T$,
we define
\[
b_t^N(s, u)=N^{\frac{d}{\alpha}-1}\sum_{x\in \mathbb{Z}^d}v(N(t-s), x)1_{\left\{u\in \frac{x}{N^{\frac{1}{\alpha}}}+\left(-\frac{1}{2N^{\frac{1}{\alpha}}}, \frac{1}{2N^{\frac{1}{\alpha}}}\right]^d, s\leq t\right\}}=\tilde{J}_{t-s}^N(u)1_{\{s\leq t\}}
\]
for any $(s, u)\in [0, T]\times \mathbb{R}^d$, where $1_A$ is the indicator function of the set $A$ and $\tilde{J}_t^N$ is defined as in the proof of Lemma \ref{lemma 4.1}. By It\^{o}'s formula, it is easy to check that
\begin{equation}\label{equ 4.18}
M_s^t=\sum_{x\in \mathbb{Z}^d}\sum_{y\in \mathbb{Z}^d\setminus\{x\}}\int_0^sv(t-r, x)\left(\eta_{r-}(y)-\eta_{r-}(x)\right)d\hat{N}_r^{x, y}.
\end{equation}
Lemma 4.1 of \cite{Xue2024} gives an analogue of \eqref{equ 4.18} for the nearest-neighbor case, the proof of which given in \cite{Xue2024} still applies in the long-range case, so we omit a detailed check of \eqref{equ 4.18} here. According to \eqref{equ 4.18} and the definition of $\mathcal{U}^N$,
\begin{equation}\label{equ 4.19}
\frac{1}{N^{\frac{3}{2}-\frac{d}{2\alpha}}}M_{tN}^{tN}=\mathcal{U}^N(b_t^N).
\end{equation}

We denote by $\mathcal{U}$ the Gaussian time-space white noise on $[0, T]\times \mathbb{R}^d$ such that $\mathcal{U}(H)$ follows the normal distribution with mean zero and variance
\[
2p(1-p)\left(\sum_{y\in \mathbb{Z}^d\setminus\{\mathbf{0}\}}\frac{\Phi^{d, \alpha}(y)}{\|y\|_2^{d+\alpha}}\right)\int_0^T\int_{\mathbb{R}^d}H^2(s, u)dsdu
\]
for any $H\in C_c([0, T]\times \mathbb{R}^d)$. For later use, here we point out that
\begin{equation}\label{equ 4.20}
\sum_{y\in \mathbb{Z}^d\setminus\{\mathbf{0}\}}\frac{\Phi^{d, \alpha}(y)}{\|y\|_2^{d+\alpha}}=\mathcal{C}_{d, \alpha},
\end{equation}
where $\mathcal{C}_{d, \alpha}$ is defined as in Section \ref{section two}. We check \eqref{equ 4.20} in Appendix \ref{Appendix A.4}. The following lemma is crucial for us to prove Lemma \ref{lemma 4.3}.

\begin{lemma}\label{lemma 4.4}
Let $(d, \alpha)$ satisfy \eqref{equ 4.4} and $\eta_0$ be distributed with $\nu_p^{d, \alpha}$. For any
\[
H\in C_c^\infty([0, T]\times\mathbb{R}^d),
\]
$\mathcal{U}^N(H)$ converges weakly to $\mathcal{U}(H)$ as $N\rightarrow+\infty$.
\end{lemma}

\proof[Proof of Lemma \ref{lemma 4.4}]
For given $H\in C_c^\infty([0, T]\times\mathbb{R}^d)$ and $0\leq t\leq T$, we define
\[
H^t(s, u)=H(s, u)1_{\{s\leq t\}}
\]
for any $(s, u)\in [0, T]\times \mathbb{R}^d$.
Let $\varpi_t^{N,H}=\mathcal{U}^N(H^t)$, then to complete the proof, we only need to show that
\[
\left\{\varpi_t^{N, H}:~0\leq t\leq T\right\}
\]
converges weakly, with respect to the Skorohod topology, to
\[
\left\{\sqrt{2p(1-p)\left(\sum_{y\in \mathbb{Z}^d\setminus\{\mathbf{0}\}}\frac{\Phi^{d, \alpha}(y)}{\|y\|_2^{d+\alpha}}\right)}
\int_0^t\sqrt{\int_{\mathbb{R}^d}H^2(s, u)du}\text{~}dW_s:~0\leq t\leq T\right\}
\]
as $N\rightarrow+\infty$, where $\{\mathcal{W}\}_{t\geq 0}$ is the standard Brownian motion defined as in Section \ref{section two}. According to the definition of $\mathcal{U}^N$, $\{\varpi_t^{N, H}\}_{t\geq 0}$ is a martingale with quadratic variation process $\{\langle\varpi^{N, H}\rangle_t\}_{t\geq 0}$ given by
\begin{align*}
\langle\varpi^{N, H}\rangle_t&=N^{\frac{d}{\alpha}-1}\sum_{x\in \mathbb{Z}^d}\sum_{y\in\mathbb{Z}^d\setminus\{x\}}\\
&\text{\quad}\int_0^t\left(\int_{\frac{x}{N^{\frac{1}{\alpha}}}+\left(-\frac{1}{2N^{\frac{1}{\alpha}}}, \frac{1}{2N^{\frac{1}{\alpha}}}\right]^d}H(s, u)du\right)^2\left(\eta_{sN-}(y)-\eta_{sN-}(x)\right)^2\frac{N}{\|y-x\|_2^{d+\alpha}}ds.
\end{align*}
Hence, by \eqref{equ 1.13},
\begin{equation}\label{equ 4.21}
\lim_{N\rightarrow+\infty}\hat{\mathbb{E}}_{\nu_p^{d, \alpha}}\langle\varpi^{N, H}\rangle_t
=2p(1-p)\left(\sum_{y\in \mathbb{Z}^d\setminus\{\mathbf{0}\}}\frac{\Phi^{d, \alpha}(y)}{\|y\|_2^{d+\alpha}}\right)\int_0^t\int_{\mathbb{R}^d}H^2(s, u)dsdu.
\end{equation}
By Theorem 1.4 in Chapter 7 of \cite{Ethier1986}, to complete the proof, we only need to show that
\[
\lim_{N\rightarrow+\infty}\langle\varpi^{N, H}\rangle_t
=2p(1-p)\left(\sum_{y\in \mathbb{Z}^d\setminus\{\mathbf{0}\}}\frac{\Phi^{d, \alpha}(y)}{\|y\|_2^{d+\alpha}}\right)\int_0^t\int_{\mathbb{R}^d}H^2(s, u)dsdu
\]
in probability. Then, by \eqref{equ 4.21}, we only need to show that
\begin{equation}\label{equ 4.22}
\lim_{N\rightarrow+\infty}{\rm Var}_{\nu_p^{d, \alpha}}\left(\langle\varpi^{N, H}\rangle_t\right)=0.
\end{equation}
According to the bilinear property of the covariance and the dominated convergence theorem, to check \eqref{equ 4.22}, we only need to show that
\begin{equation}\label{equ 4.23}
\lim_{t\rightarrow+\infty}\left|{\rm Cov}_{\nu_p^{d, \alpha}}\left(\left(\eta_t(y_1)-\eta_t(x_1)\right)^2, \left(\eta_0(y_2)-\eta_0(x_2)\right)^2\right)\right|=0
\end{equation}
for any given $x_1, y_1, x_2, y_2\in \mathbb{Z}^d$. According to the definition of $\nu_p^{d, \alpha}$,
\begin{align}\label{equ 4.23 two}
&{\rm Cov}_{\nu_p^{d, \alpha}}\left(\left(\eta_t(y_1)-\eta_t(x_1)\right)^2, \left(\eta_0(y_2)-\eta_0(x_2)\right)^2\right)\notag\\
&=\lim_{s\rightarrow+\infty}{\rm Cov}_{\mu_p}\left(\left(\eta_{t+s}(y_1)-\eta_{t+s}(x_1)\right)^2, \left(\eta_s(y_2)-\eta_s(x_2)\right)^2\right).
\end{align}
By Proposition \ref{proposition voter dual},
\[
{\rm Cov}_{\mu_p}\left(\left(\eta_{t+s}(y_1)-\eta_{t+s}(x_1)\right)^2, \left(\eta_s(y_2)-\eta_s(x_2)\right)^2\right)
={\rm \uppercase\expandafter{\romannumeral3}}-{\rm \uppercase\expandafter{\romannumeral4}},
\]
where
\[
{\rm \uppercase\expandafter{\romannumeral3}}=\mathbb{E}\left(\hat{\mathbb{E}}_{\mu_p}
\left(\left(\eta_0\left(\mathcal{X}_{t+s, y_1}^{d, \alpha, 0}\right)-\eta_0\left(\mathcal{X}_{t+s, x_1}^{d, \alpha, 0}\right)\right)^2\left(\eta_0\left(\mathcal{X}_{t+s, y_2}^{d, \alpha, t}\right)-\eta_0\left(\mathcal{X}_{t+s, x_2}^{d, \alpha, t}\right)\right)^2\right)\right)
\]
and
\begin{align*}
&{\rm \uppercase\expandafter{\romannumeral4}}=\\
&\mathbb{E}\left(\hat{\mathbb{E}}_{\mu_p}\left(\left(\eta_0\left(\mathcal{X}_{t+s, y_1}^{d, \alpha, 0}\right)-\eta_0\left(\mathcal{X}_{t+s, x_1}^{d, \alpha, 0}\right)\right)^2\right)\right)\mathbb{E}\left(\hat{\mathbb{E}}_{\mu_p}\left(\left(\eta_0\left(\mathcal{X}_{t+s, y_2}^{d, \alpha, t}\right)-\eta_0\left(\mathcal{X}_{t+s, x_2}^{d, \alpha, t}\right)\right)^2\right)\right).
\end{align*}
According to an argument similar to that given in the proof of \eqref{equ 4.16},
\begin{equation}\label{equ 4.24}
\left|{\rm \uppercase\expandafter{\romannumeral3}}-{\rm \uppercase\expandafter{\romannumeral4}}\right|
\leq 2\mathbb{P}\left(\hat{\tau}_{x_1, y_1, x_2, y_2}^t\leq t+s\right),
\end{equation}
where
\[
\hat{\tau}_{x_1, y_1, x_2, y_2}^t=\inf\left\{r:~r\geq t, \left\{\mathcal{X}_{r, y_1}^{d, \alpha, 0},
\mathcal{X}_{r, x_1}^{d, \alpha, 0}\right\}\bigcap\left\{\mathcal{X}_{r, y_2}^{d, \alpha, t}, \mathcal{X}_{r, x_2}^{d, \alpha, t}\right\}\neq \emptyset\right\}.
\]
For any $z, w\in \mathbb{Z}^d$,
\begin{align*}
&\mathbb{P}\left(\hat{\tau}_{x_1, y_1, x_2, y_2}^t\leq t+s\Big|\left(\mathcal{X}_{t, y_1}^{d, \alpha, 0},
\mathcal{X}_{t, x_1}^{d, \alpha, 0}\right)=(z, w)\right)\\
&\leq (1-\Phi^{d, \alpha}(z-y_2))+(1-\Phi^{d, \alpha}(z-x_2))+(1-\Phi^{d, \alpha}(w-x_2))+(1-\Phi^{d, \alpha}(w-y_2)).
\end{align*}
For $(d, \alpha)$ satisfying \eqref{equ 4.4}, $\{X_t^{d, \alpha}\}_{t\geq 0}$ is transient. Hence, for any $\epsilon>0$, there exists $M>0$ such that
\[
1-\Phi^{d, \alpha}(x)\leq \epsilon
\]
when $\|x\|_2>M$. Consequently, by the total probability formula,
\begin{align*}
\mathbb{P}\left(\hat{\tau}_{x_1, y_1, x_2, y_2}^t\leq t+s\right)
\leq 4\epsilon \mathbb{P}\left(E_1\right)+1-\mathbb{P}(E_1)\leq 4\epsilon+1-\mathbb{P}(E_1),
\end{align*}
where $E_1$ is the event that
\[
\min\left\{\|\mathcal{X}_{t, y_1}^{d, \alpha, 0}-y_2\|_2, \|\mathcal{X}_{t, y_1}^{d, \alpha, 0}-x_2\|_2,
\|\mathcal{X}_{t, x_1}^{d, \alpha, 0}-y_2\|_2, \|\mathcal{X}_{t, x_1}^{d, \alpha, 0}-x_2\|_2\right\}\geq M.
\]
Then, by  \eqref{equ 4.23 two} and \eqref{equ 4.24},
\begin{equation}\label{equ 4.26}
\left|{\rm Cov}_{\nu_p^{d, \alpha}}\left(\left(\eta_t(y_1)-\eta_t(x_1)\right)^2, \left(\eta_0(y_2)-\eta_0(x_2)\right)^2\right)\right|
\leq 8\epsilon+2\left(1-\mathbb{P}(E_1)\right).
\end{equation}
By \eqref{equ 1.4 inverse formula random walk}, $p_t^{d, \alpha}(x, y)\leq p_t^{d, \alpha}(\mathbf{0}, \mathbf{0})$ for any $x, y\in \mathbb{Z}^d$. Hence
\[
1-\mathbb{P}(E_1)\leq 4(2M+1)^d p_t^{d, \alpha}(\mathbf{0}, \mathbf{0}).
\]
By Lemma \ref{lemma LCLT}, $\lim_{t\rightarrow+\infty}p_t^{d, \alpha}(\mathbf{0}, \mathbf{0})=0$. Then, by \eqref{equ 4.26},
\[
\limsup_{t\rightarrow+\infty}\left|{\rm Cov}_{\nu_p^{d, \alpha}}\left(\left(\eta_t(y_1)-\eta_t(x_1)\right)^2, \left(\eta_0(y_2)-\eta_0(x_2)\right)^2\right)\right|\leq 8\epsilon.
\]
Since $\epsilon$ is arbitrary, let $\epsilon\rightarrow 0$ and then \eqref{equ 4.23} holds, which completes the proof.
\qed

Now we prove Lemma \ref{lemma 4.3}.

\proof[Proof of Lemma \ref{lemma 4.3}]
For $0\leq s, t\leq T$ and $u\in \mathbb{R}^d$, we define
\[
b_t(s, u)=J_{t-s}(u)1_{\{s\leq t\}},
\]
where $J_t$ is defined as in the proof of Lemma \ref{lemma 4.1}. By Lemma \ref{lemma 4.3}, for any $c_1, c_2, \ldots, c_m\in \mathbb{R}$,
\[
\mathcal{U}^N(\sum_{i=1}^mc_ib_{t_i})=\sum_{i=1}^mc_i\mathcal{U}^N(b_{t_i})
\]
converges weakly to
\[
\mathcal{U}(\sum_{i=1}^mc_ib_{t_i})=\sum_{i=1}^mc_i\mathcal{U}(b_{t_i})
\]
as $N\rightarrow+\infty$. Hence,
\[
\left(\mathcal{U}^N(b_{t_1}), \ldots, \mathcal{U}^N(b_{t_m})\right)
\]
converges weakly to the $\mathbb{R}^m$-valued Gaussian random variable
\[
\left(\mathcal{U}(b_{t_1}), \ldots, \mathcal{U}(b_{t_m})\right)
\]
as $N\rightarrow+\infty$. By \eqref{equ 4.9}, $\mathcal{U}^N(b_{t_k})-\mathcal{U}^N(b^N_{t_k})$ converges weakly to $0$ as $N\rightarrow+\infty$. Conseqently, by \eqref{equ 4.19}, Lemma \ref{lemma 4.3} holds with
\[
(\Xi_{t_1},\ldots, \Xi_{t_m})=\left(\mathcal{U}(b_{t_1}), \ldots, \mathcal{U}(b_{t_m})\right).
\]
\qed

Now we deal with the joint distribution of $M_{tN}^{tN}$ and $V^{tN}(0, \eta_0)$.
\begin{lemma}\label{lemma 4.5}
Let $(d, \alpha)$ satisfy \eqref{equ 4.4} and $\eta_0$ be distributed with $\nu_p^{d, \alpha}$. For any $0\leq t_1<t_2<\ldots<t_m$,
the joint distribution of
\[
\frac{1}{N^{\frac{3}{2}-\frac{d}{2\alpha}}}\left(V^{t_1N}(0, \eta_0), V^{t_2N}(0, \eta_0), \ldots, V^{t_mN}(0, \eta_0)\right)
\]
and
\[
\frac{1}{N^{\frac{3}{2}-\frac{d}{2\alpha}}}\left(M_{t_1N}^{t_1N}, \ldots, M_{t_mN}^{t_mN}\right)
\]
converges weakly, as $N\rightarrow+\infty$, to the independent coupling of
\[
\left(\Gamma_{t_1}, \ldots, \Gamma_{t_m}\right)
\text{\quad and\quad}
(\Xi_{t_1},\ldots, \Xi_{t_m}),
\]
which are given in Lemmas \ref{lemma 4.1} and \ref{lemma 4.3} respectively.
\end{lemma}

\proof[Proof of Lemma \ref{lemma 4.5}]
For any $\eta\in \{0, 1\}^{\mathbb{Z}^d}$, $t>0$, $H\in C_c^\infty([0, T]\times \mathbb{R}^d)$, and integers $N, M\geq 1$, we denote by
$\mathcal{E}(\eta, t, H, M, N)$ the term
\[
\hat{\mathbb{P}}_\eta\left(\left|\varpi_t^{N, H}-k_t^H\right|\geq \frac{1}{M}\right),
\]
where
\[
k_t^H=2p(1-p)\left(\sum_{y\in \mathbb{Z}^d\setminus\{\mathbf{0}\}}\frac{\Phi^{d, \alpha}(y)}{\|y\|_2^{d+\alpha}}\right)\int_0^t\int_{\mathbb{R}^d}H^2(s, u)dsdu.
\]
According to the proof of Lemma \ref{lemma 4.4}, under $\nu_p^{d, \alpha}$, the random variable $\mathcal{E}(\cdot, t, H, M, N)$ converges to $0$ in probability as $N\rightarrow+\infty$. Hence, there exists a subsequence $\{N_k\}_{k\geq 1}$ of $\{N\}_{N\geq 1}$ such that, under $\nu_p^{d, \alpha}$, $\mathcal{E}(\cdot, t, H, M, N_k)$ converges almost surely, as $k\rightarrow+\infty$, to $0$ for all $M\geq 1$, $t\in \mathbb{Q}$ and $H$ in a given countable dense set of $C_c^\infty([0, T]\times \mathbb{R}^d)$. Hence, there exists $A\subseteq \{0, 1\}^{\mathbb{Z}^d}$ such that $\nu_p^{d, \alpha}(A)=1$ and for each $\eta\in A$, under $\hat{\mathbb{P}}_\eta$,
\[
\mathcal{U}^{N_k}(H)
\]
converges weakly, as $k\rightarrow+\infty$, to $\mathcal{U}(H)$ for all $H\in C_c^\infty([0, T]\times \mathbb{R}^d)$. Consequently, for any $\eta\in A$, when $\eta_0$ is distributed with $\hat{P}_\eta$,
\[
\frac{1}{N_k^{\frac{3}{2}-\frac{d}{2\alpha}}}\left(M_{t_1N_k}^{t_1N_k}, \ldots, M_{t_mN_k}^{t_mN_k}\right)
\]
converges weakly, as $k\rightarrow+\infty$, to $(\Xi_{t_1},\ldots, \Xi_{t_m})$. Then, by the iterated expectation law, when $\eta_0$ is distributed with $\nu_p^{d, \alpha}$, the joint distribution of
\[
\frac{1}{N_k^{\frac{3}{2}-\frac{d}{2\alpha}}}\left(V^{t_1N_k}(0, \eta_0), V^{t_2N_k}(0, \eta_0), \ldots, V^{t_mN_k}(0, \eta_0)\right)
\]
and
$\frac{1}{N_k^{\frac{3}{2}-\frac{d}{2\alpha}}}\left(M_{t_1N_k}^{t_1N_k}, \ldots, M_{t_mN_k}^{t_mN_k}\right)$ converges weakly to the independent coupling of
\[
\left(\Gamma_{t_1}, \ldots, \Gamma_{t_m}\right)
\text{\quad and\quad}
(\Xi_{t_1},\ldots, \Xi_{t_m}).
\]
According to similar arguments, when $\eta_0$ is distributed with $\nu_p^{d, \alpha}$, any subsequence of
\[
\left\{\left(\frac{\left(V^{t_1N}(0, \eta_0), V^{t_2N}(0, \eta_0), \ldots, V^{t_mN}(0, \eta_0)\right)}{N^{\frac{3}{2}-\frac{d}{2\alpha}}},
\frac{\left(M_{t_1N}^{t_1N}, \ldots, M_{t_mN}^{t_mN}\right)}{N^{\frac{3}{2}-\frac{d}{2\alpha}}}\right)\right\}_{N\geq 1}
\]
have a subsequence that converges weakly to the independent coupling of
\[
\left(\Gamma_{t_1}, \ldots, \Gamma_{t_m}\right)
\text{\quad and\quad}
(\Xi_{t_1},\ldots, \Xi_{t_m})
\]
and the proof is complete.
\qed

The following lemma shows that any finite dimensional distribution of
\[
\left\{\frac{1}{N^{\frac{3}{2}-\frac{d}{2\alpha}}}\int_0^{tN}\left(\eta_s(\mathbf{0})-p\right)ds:~0\leq t\leq T\right\}
\]
converges weakly to the corresponding finite dimensional distribution of
\[
\left\{\sqrt{2f_1^{d, \alpha}(\mathbf{0})\frac{\alpha^3p(1-p)\mathcal{C}_{d, \alpha}}{(d-\alpha)(2\alpha-d)(3\alpha-d)}}B_t^{\frac{3}{2}-\frac{d}{2\alpha}}:~0\leq t\leq T\right\}.
\]
\begin{lemma}\label{lemma 4.6}
Let $(d, \alpha)$ satisfy \eqref{equ 4.4} and $\eta_0$ be distributed with $\nu_p^{d, \alpha}$. For any $0\leq t_1<t_2<\ldots<t_m$,
\[
\frac{1}{N^{\frac{3}{2}-\frac{d}{2\alpha}}}\left(\int_0^{t_1N}\left(\eta_s(\mathbf{0})-p\right)ds, \ldots, \int_0^{t_mN}\left(\eta_s(\mathbf{0})-p\right)ds\right)
\]
converges weakly to
\[
\sqrt{2f_1^{d, \alpha}(\mathbf{0})\frac{\alpha^3p(1-p)\mathcal{C}_{d, \alpha}}{(d-\alpha)(2\alpha-d)(3\alpha-d)}}
\left(B_{t_1}^{\frac{3}{2}-\frac{d}{2\alpha}}, \ldots, B_{t_m}^{\frac{3}{2}-\frac{d}{2\alpha}}\right)
\]
as $N\rightarrow+\infty$.
\end{lemma}

\proof[Proof of Lemma \ref{lemma 4.6}]
By Lemmas \ref{lemma 4.1}, \ref{lemma 4.3}, \ref{lemma 4.5} and Equation \eqref{equ 4.3 martingale decomposition},
\[
\frac{1}{N^{\frac{3}{2}-\frac{d}{2\alpha}}}\left(\int_0^{t_1N}\left(\eta_s(\mathbf{0})-p\right)ds, \ldots, \int_0^{t_mN}\left(\eta_s(\mathbf{0})-p\right)ds\right)
\]
converges weakly to the $\mathbb{R}^m$-valued Gaussian random variable
\[
\left(\Gamma_{t_1}+\Xi_{t_1}, \Gamma_{t_2}+\Xi_{t_2}, \ldots, \Gamma_{t_m}+\Xi_{t_m}\right)
\]
as $N\rightarrow+\infty$, where $\{\Gamma_{t_i}\}_{i=1}^m, \{\Xi_{t_i}\}_{i=1}^m$ are defined as in Lemmas \ref{lemma 4.1}, \ref{lemma 4.3} respectively and are independently coupled in the same probability space. Then, to complete the proof, we only need to show that
\begin{align}\label{equ 4.27}
&{\rm Cov}\left(\Gamma_{t_i}+\Xi_{t_i}, \Gamma_{t_j}+\Xi_{t_j}\right)\\
&=\frac{f_1^{d, \alpha}(\mathbf{0})\alpha^3p(1-p)\mathcal{C}_{d, \alpha}}{(d-\alpha)(2\alpha-d)(3\alpha-d)}
\left(t_j^{3-\frac{d}{\alpha}}+t_i^{3-\frac{d}{\alpha}}-(t_j-t_i)^{3-\frac{d}{\alpha}}\right)\notag
\end{align}
for any $t_i\leq t_j$. Since $\{\Gamma_{t_i}\}_{i=1}^m$ and $\{\Xi_{t_i}\}_{i=1}^m$ are independent,
\begin{equation}\label{equ 4.28}
{\rm Cov}\left(\Gamma_{t_i}+\Xi_{t_i}, \Gamma_{t_j}+\Xi_{t_j}\right)
={\rm Cov}(\Gamma_{t_i}, \Gamma_{t_j})+{\rm Cov}(\Xi_{t_i}, \Xi_{t_j}).
\end{equation}
As we have shown in the proof of Lemma \ref{lemma 4.1}, $\Gamma_{t_i}=\mathcal{Y}_0(J_{t_i})$. Therefore,
\begin{align*}
{\rm Cov}(\Gamma_{t_i}, \Gamma_{t_j})=\left(2p(1-p)\sum_{y\neq \mathbf{0}}\frac{\Phi^{d, \alpha}(y)}{\|y\|^{d+\alpha}}\right)\int_0^{+\infty}\int_{\mathbb{R}^d}\mathcal{P}_sJ_{t_i}(u)\mathcal{P}_sJ_{t_j}(u)dsdu.
\end{align*}
Note that, according to the definition of $J_t$,
\begin{align*}
&\mathcal{P}_sJ_{t_i}(u)\mathcal{P}_sJ_{t_j}(u)\\
&=\int_{\mathbb{R}^d}J_{t_i}(v_1)f_s^{d, \alpha}(v_1-u)dv_1\int_{\mathbb{R}^d}J_{t_j}(v_2)f_s^{d, \alpha}(v_2-u)dv_2\\
&=\int_0^{t_i}dr\int_0^{t_j}d\theta\int_{\mathbb{R}^d}dv_1\int_{\mathbb{R}^d}f_r^{d, \alpha}(v_1)f_\theta^{d, \alpha}(v_2)f_s^{d, \alpha}(v_1-u)f_s^{d, \alpha}(v_2-u)dv_2.
\end{align*}
By Kolmogorov-Chapman equation,
\[
\int_{\mathbb{R}^d}\int_{\mathbb{R}^d}\int_{\mathbb{R}^d}f_r^{d, \alpha}(v_1)f_\theta^{d, \alpha}(v_2)f_s^{d, \alpha}(v_1-u)f_s^{d, \alpha}(v_2-u)dv_1dv_2du
=f_{r+\theta+2s}^{d, \alpha}(\mathbf{0}).
\]
Hence,
\[
\int_0^{+\infty}\int_{\mathbb{R}^d}\mathcal{P}_sJ_{t_i}(u)\mathcal{P}_sJ_{t_j}(u)duds
=\int_0^{t_i}dr\int_0^{t_j}d\theta\int_0^{+\infty}f_{r+\theta+2s}^{d, \alpha}(\mathbf{0})ds.
\]
As we have pointed out in the proof of \eqref{equ 4.8},
\[
f_r^{d, \alpha}(0)=r^{-\frac{d}{\alpha}}f_1^{d, \alpha}(\mathbf{0}).
\]
Therefore,
\begin{align*}
&\int_0^{t_i}dr\int_0^{t_j}d\theta\int_0^{+\infty}f_{r+\theta+2s}^{d, \alpha}(0)ds\\
&=f_1^{d, \alpha}(\mathbf{0})\int_0^{t_i}dr\int_0^{t_j}d\theta\int_0^{+\infty}(r+\theta+2s)^{-\frac{d}{\alpha}}ds\\
&=\frac{f_1^{d, \alpha}(\mathbf{0})\alpha^3}{2(d-\alpha)(2\alpha-d)(3\alpha-d)}\left((t_i+t_j)^{3-\frac{d}{\alpha}}-t_i^{3-\frac{d}{\alpha}}
-t_j^{3-\frac{d}{\alpha}}\right)
\end{align*}
and hence, by \eqref{equ 4.20},
\begin{equation}\label{equ 4.29}
{\rm Cov}(\Gamma_{t_i}, \Gamma_{t_j})
=\frac{f_1^{d, \alpha}(\mathbf{0})\alpha^3p(1-p)\mathcal{C}_{d, \alpha}}{(d-\alpha)(2\alpha-d)(3\alpha-d)}\left((t_i+t_j)^{3-\frac{d}{\alpha}}-t_i^{3-\frac{d}{\alpha}}
-t_j^{3-\frac{d}{\alpha}}\right).
\end{equation}
As we have shown in the proof of Lemma \ref{lemma 4.3}, $\Xi_{t_i}=\mathcal{U}(b_{t_i})$ and hence
\[
{\rm Cov}(\Xi_{t_i}, \Xi_{t_j})=\left(2p(1-p)\sum_{y\neq \mathbf{0}}\frac{\Phi^{d, \alpha}(y)}{\|y\|^{d+\alpha}}\right)
\int_0^{t_i}\int_{\mathbb{R}^d}b_{t_j}(s, u)b_{t_i}(s, u)dsdu.
\]
According to the definition of $b_t$ and Kolmogorov-Chapman equation,
\begin{align*}
&\int_0^{t_i}\int_{\mathbb{R}^d}b_{t_j}(s, u)b_{t_i}(s, u)dsdu\\
&=\int_0^{t_i}\int_{\mathbb{R}^d}\left(\int_0^{t_j-s}f_r^{d, \alpha}(u)dr\int_0^{t_i-s}f_v^{d, \alpha}(u)dv\right)dsdu\\
&=\int_0^{t_i}\left(\int_0^{t_j-s}\int_0^{t_i-s}f_{r+v}^{d, \alpha}(\mathbf{0})drdv\right)ds\\
&=f_1^{d,\alpha}(\mathbf{0})\int_0^{t_i}\left(\int_0^{t_i-s}dv\int_0^{t_j-s}(r+v)^{-\frac{d}{\alpha}}dr\right)ds\\
&=\frac{f_1^{d, \alpha}(\mathbf{0})\alpha^3}{(d-\alpha)(2\alpha-d)(3\alpha-d)}
\left(t_j^{3-\frac{d}{\alpha}}+t_i^{3-\frac{d}{\alpha}}-\frac{(t_i+t_j)^{3-\frac{d}{\alpha}}}{2}-\frac{(t_j-t_i)^{3-\frac{d}{\alpha}}}{2}\right).
\end{align*}
Therefore, by \eqref{equ 4.20},
\begin{align}\label{equ 4.30}
&{\rm Cov}(\Xi_{t_i}, \Xi_{t_j})\\
&=\frac{f_1^{d, \alpha}(\mathbf{0})\alpha^3p(1-p)\mathcal{C}_{d, \alpha}}{(d-\alpha)(2\alpha-d)(3\alpha-d)}
\left(2t_j^{3-\frac{d}{\alpha}}+2t_i^{3-\frac{d}{\alpha}}-(t_i+t_j)^{3-\frac{d}{\alpha}}-(t_j-t_i)^{3-\frac{d}{\alpha}}\right).\notag
\end{align}
Equation \eqref{equ 4.27} follows from \eqref{equ 4.28}, \eqref{equ 4.29} and \eqref{equ 4.30}. Since \eqref{equ 4.27} holds, the proof is complete.
\qed

At last we prove Theorem \ref{theorem 2.1 main result d=1, 2} in the case $\alpha\in \left(d/2, \min\{d, 2\}\right)$.

\proof[Proof of Theorem \ref{theorem 2.1 main result d=1, 2} in the case $\alpha\in \left(d/2, \min\{d, 2\}\right)$]
Let $(d, \alpha)$ satisfy \eqref{equ 4.4} and $\eta_0$ be distributed with $\nu_p^{d, \alpha}$. By Lemma \ref{lemma 4.6}, we only need to show that
\[
\left\{\frac{1}{N^{\frac{3}{2}-\frac{d}{2\alpha}}}\int_0^{tN}\left(\eta_s(\mathbf{0})-p\right)ds:~0\leq t\leq T\right\}_{N\geq 1}
\]
are tight under the uniform topology. To check this tightness, we only need to show that there exist $a, b, c>0$ independent of $s, t, N$ such that
\begin{equation}\label{equ 4.31}
\hat{\mathbb{E}}_{\nu_p^{d, \alpha}}\left(\left|\frac{1}{N^{\frac{3}{2}-\frac{d}{2\alpha}}}\int_{sN}^{tN}\left(\eta_r(\mathbf{0})-p\right)dr\right|^{a}\right)
\leq c|t-s|^{b+1}
\end{equation}
for any $0<s<t$ and $N\geq 1$. According to the invariance of $\nu_p^{d, \alpha}$ and the bilinear property of the covariance, for any $t\geq 0$,
\begin{align*}
\hat{\mathbb{E}}_{\nu_p^{d, \alpha}}\left(\left(\int_{0}^{t}\left(\eta_r(\mathbf{0})-p\right)dr\right)^2\right)&=
{\rm Var_{\nu_p^{d, \alpha}}}\left(\int_{0}^{t}\left(\eta_r(\mathbf{0})-p\right)dr\right)\\
&=\int_0^t\int_0^t{\rm Cov}_{\nu_p^{d, \alpha}}\left(\eta_r(\mathbf{0}), \eta_\theta(\mathbf{0})\right)drd\theta\\
&=2\int_0^tdr\int_0^r{\rm Cov}_{\nu_p^{d, \alpha}}\left(\eta_{r-\theta}(\mathbf{0}), \eta_0(\mathbf{0})\right)d\theta\\
&=2\int_0^tdr\int_0^r{\rm Cov}_{\nu_p^{d, \alpha}}\left(\eta_{\theta}(\mathbf{0}), \eta_0(\mathbf{0})\right)d\theta.
\end{align*}
We claim that
\begin{equation}\label{equ 4.32}
{\rm Cov}_{\nu_p^{d, \alpha}}\left(\eta_{\theta}(\mathbf{0}), \eta_0(\mathbf{0})\right)
=p(1-p)\mathcal{C}_{d, \alpha}\int_0^{+\infty}p_{\theta+v}^{d, \alpha}(\mathbf{0}, \mathbf{0})dv.
\end{equation}
We check \eqref{equ 4.32} in Appendix \ref{Appendix A.5}. Then,
\[
\hat{\mathbb{E}}_{\nu_p^{d, \alpha}}\left(\left(\int_{0}^{t}\left(\eta_r(\mathbf{0})-p\right)dr\right)^2\right)
=2\mathcal{C}_{d, \alpha}p(1-p)\int_0^tdr\int_0^r d\theta\int_0^{+\infty}p_{\theta+v}^{d, \alpha}(\mathbf{0}, \mathbf{0})dv.
\]
By Lemma \ref{lemma LCLT}, there exists $C_4<+\infty$ independent of $t$ such that
\[
p_t^{d, \alpha}(\mathbf{0}, \mathbf{0})\leq C_4 t^{-\frac{d}{\alpha}}
\]
for all $t\geq 0$. Therefore,
\begin{align*}
\hat{\mathbb{E}}_{\nu_p^{d, \alpha}}\left(\left(\int_{0}^{t}\left(\eta_r(\mathbf{0})-p\right)dr\right)^2\right)
&\leq 2\mathcal{C}_{d, \alpha}p(1-p)C_4\int_0^tdr\int_0^r d\theta\int_0^{+\infty}(\theta+v)^{-\frac{d}{\alpha}}dv\\
&=\frac{2\mathcal{C}_{d, \alpha}p(1-p)C_4\alpha^3}{(d-\alpha)(2\alpha-d)(3\alpha-d)}t^{3-\frac{d}{\alpha}}.
\end{align*}
Then, according to the invariance of $\nu_p^{d, \alpha}$, \eqref{equ 4.31} holds with $a=2$, $b=2-\frac{d}{\alpha}$ and
\[
c=\frac{2\mathcal{C}_{d, \alpha}p(1-p)C_4\alpha^3}{(d-\alpha)(2\alpha-d)(3\alpha-d)}.
\]
Since \eqref{equ 4.31} holds, the proof is complete.
\qed

\subsection{Proof of Theorem \ref{theorem 2.1 main result d=1, 2}: the case of $\alpha\leq \frac{d}{2}$}\label{subsection 4.2}
In this subsection, we prove Theorem \ref{theorem 2.1 main result d=1, 2} in the case $\alpha\leq \frac{d}{2}$. Hence, throughout this subsection we assume that $d\in \{1, 2, 3\}$ and $\alpha\leq \frac{d}{2}$. We first introduce some notations for later use. For any $x\in \mathbb{Z}^d$, we define
\[
\Upsilon_N(x)=\int_0^{+\infty}e^{-\frac{s}{N}}p_s^{d, \alpha}(\mathbf{0}, x)dx.
\]
Then, for any $\eta\in \{0, 1\}^{\mathbb{Z}^d}$, we define
\[
\mathcal{I}_N(\eta)=\sum_{x\in \mathbb{Z}^d}\left(\eta(x)-p\right)\Upsilon_N(x).
\]
Let $\mathcal{G}$ be defined as in \eqref{equ 1.1 generator of long-range voter model}, then by Kolmogorov-Chapman equation, it is not difficult to show that
\begin{equation}\label{equ 4.33}
\mathcal{G}\mathcal{I}_N(\eta)=\frac{1}{N}\mathcal{I}_N(\eta)-\left(\eta(\mathbf{0})-p\right).
\end{equation}
Reference \cite{Xue2024} gives a detailed proof of an analogue of \eqref{equ 4.33} for the nearest-neighbor case. So in this paper we omit the check of \eqref{equ 4.33}. For any $t\geq 0$, we define
\[
\widetilde{M}_t^N=\mathcal{I}_N(\eta_t)-\mathcal{I}_N(\eta_0)-\int_0^t\mathcal{G}\mathcal{I}_N(\eta_s)ds,
\]
then $\{\widetilde{M}_t^N\}_{t\geq 0}$ is a martingale according to the Dynkin's martingale formula. By \eqref{equ 4.33},
\begin{equation}\label{equ 4.34 martingale decomposition}
\frac{\int_0^{tN}\left(\eta_s(\mathbf{0})-p\right)ds}{\Lambda_{d, \alpha}(N)}=
\frac{\widetilde{M}_{tN}^N}{\Lambda_{d, \alpha}(N)}+\frac{\frac{1}{N}\int_0^{tN}\mathcal{I}_N(\eta_s)ds+\mathcal{I}_N(\eta_0)-\mathcal{I}_N(\eta_{tN})}{\Lambda_{d, \alpha}(N)}.
\end{equation}
The following lemma shows that the term $\frac{\frac{1}{N}\int_0^{tN}\mathcal{I}_N(\eta_s)ds+\mathcal{I}_N(\eta_0)-\mathcal{I}_N(\eta_{tN})}{\Lambda_{d, \alpha}(N)}$ in decomposition \eqref{equ 4.34 martingale decomposition} converges weakly to $0$ as $N\rightarrow+\infty$.

\begin{lemma}\label{lemma 4.7}
Let $d\in \{1, 2, 3\}$, $\alpha\leq \frac{d}{2}$ and $\eta_0$ be distributed with $\nu_p^{d, \alpha}$. For any $t\geq 0$,
\[
\lim_{N\rightarrow+\infty}\frac{\frac{1}{N}\int_0^{tN}\mathcal{I}_N(\eta_s)ds+\mathcal{I}_N(\eta_0)-\mathcal{I}_N(\eta_{tN})}{\Lambda_{d, \alpha}(N)}=0
\]
in probability.
\end{lemma}

\proof[Proof of Lemma \ref{lemma 4.7}]
According to the invariance of $\nu_p^{d, \alpha}$ and Cauchy-Schwarz inequality, we only need to show that
\begin{equation}\label{equ 4.35}
\lim_{N\rightarrow+\infty}\frac{1}{\Lambda^2_{d, \alpha}(N)}\hat{\mathbb{E}}_{\nu_p^{d, \alpha}}\left(\mathcal{I}_N^2(\eta_0)\right)=0.
\end{equation}
By \eqref{equ 1.14}, \eqref{equ 1.15}, \eqref{equ 4.9 four} and Kolmogorov-Chapman equation,
\begin{align*}
\hat{\mathbb{E}}_{\nu_p^{d, \alpha}}\left(\mathcal{I}_N^2(\eta_0)\right)&=\sum_{x\in \mathbb{Z}^d}\sum_{y\in \mathbb{Z}^d}
{\rm Cov}\left(\eta(x), \eta(y)\right)\Upsilon_N(x)\Upsilon_N(y)\\
&=\sum_{x\in \mathbb{Z}^d}\sum_{y\in \mathbb{Z}^d}p(1-p)\mathcal{C}_{d, \alpha}\int_0^{+\infty}p_s^{d, \alpha}(\mathbf{0}, x-y)ds\Upsilon_N(x)\Upsilon_N(y)\\
&=p(1-p)\mathcal{C}_{d, \alpha}\int_0^{+\infty}\int_0^{+\infty}\int_0^{+\infty}e^{-\frac{s_1+s_2}{N}}p_{s+s_1+s_2}^{d, \alpha}(\mathbf{0}, \mathbf{0})ds_1ds_2ds\\
&=p(1-p)\mathcal{C}_{d, \alpha}\int_0^{+\infty}\int_0^{+\infty}re^{-\frac{r}{N}}p_{s+r}^{d, \alpha}(\mathbf{0}, \mathbf{0})drds.
\end{align*}
When $\alpha=\frac{d}{2}$, by Lemma \ref{lemma LCLT}, there exists $C_{15}<+\infty$ independent of $s$ such that
\[
p_s^{d, \frac{d}{2}}(\mathbf{0}, \mathbf{0})\leq C_{15}s^{-2}
\]
for any $s>0$. Hence, when $\alpha=\frac{d}{2}$,
\begin{align*}
\int_0^{+\infty}\int_0^{+\infty}re^{-\frac{r}{N}}p_{s+r}^{d, \frac{d}{2}}(\mathbf{0}, \mathbf{0})drds
&\leq C_{15}\int_0^{+\infty}re^{-\frac{r}{N}}\frac{1}{r}dr=NC_{15}
\end{align*}
and hence
\[
\frac{1}{\Lambda^2_{d, \frac{d}{2}}(N)}\hat{\mathbb{E}}_{\nu_p^{d, \frac{d}{2}}}\left(\mathcal{I}_N^2(\eta_0)\right)=O(1)(\log N)^{-1}.
\]
Consequently, Equation \eqref{equ 4.35} holds when $\alpha=\frac{d}{2}$.

When $\alpha<\frac{d}{2}$, by Lemma \ref{lemma LCLT}, there exists $C_{16}<+\infty$ independent of $s$ such that
\[
p_s^{d, \alpha}(\mathbf{0}, \mathbf{0})\leq C_{16}s^{-\frac{d}{\alpha}}
\]
for any $s>0$. Let $C_{17}=\int_0^{+\infty}p_s^{d, \alpha}(\mathbf{0}, \mathbf{0})ds$, which is finite according to the transience, then
\begin{align*}
&\int_0^{+\infty}\int_0^{+\infty}re^{-\frac{r}{N}}p_{s+r}^{d, \alpha}(\mathbf{0}, \mathbf{0})drds\\
&\leq C_{17}\int_0^1re^{-\frac{r}{N}}dr+C_{16}\int_1^{+\infty}re^{-\frac{r}{N}}\left(\int_r^{+\infty}s^{-\frac{d}{\alpha}}ds\right)dr\\
&\leq C_{17}+O(1)\int_1^{+\infty}r^{2-\frac{d}{\alpha}}e^{-\frac{r}{N}}dr=C_{17}+O(1)N^{3-\frac{d}{\alpha}}\int_{\frac{1}{N}}^{+\infty}v^{2-\frac{d}{\alpha}}e^{-v}dv
\end{align*}
when $\alpha<\frac{d}{2}$. Hence, when $\alpha<\frac{d}{2}$,
\[
\frac{1}{\Lambda^2_{d, \alpha}(N)}\hat{\mathbb{E}}_{\nu_p^{d, \alpha}}\left(\mathcal{I}_N^2(\eta_0)\right)
\leq \frac{O(1)}{N}+O(1)N^{2-\frac{d}{\alpha}}\int_{\frac{1}{N}}^{+\infty}v^{2-\frac{d}{\alpha}}e^{-v}dv.
\]
By L'H\^{o}pital's rule, when $\alpha<\frac{d}{2}$,
\begin{align*}
\lim_{a\rightarrow 0^+}\frac{\int_a^{+\infty}v^{2-\frac{d}{\alpha}}e^{-v}dv}{a^{2-\frac{d}{\alpha}}}
&=\lim_{a\rightarrow 0^+}\frac{a^{2-\frac{d}{\alpha}}e^{-a}}{a^{1-\frac{d}{\alpha}}(\frac{d}{\alpha}-2)}\\
&=\lim_{a\rightarrow 0^+}\frac{ae^{-a}}{\frac{d}{\alpha}-2}=0.
\end{align*}
Hence, Equation \eqref{equ 4.35} also holds when $\alpha<\frac{d}{2}$ and the proof is complete.
\qed

The following lemma is an analogue of Lemma \ref{lemma 4.6}.
\begin{lemma}\label{lemma 4.8}
Let $d\in \{1, 2, 3\}$, $\alpha\leq \frac{d}{2}$ and $\eta_0$ be distributed with $\nu_p^{d, \alpha}$. For any $0\leq t_1<t_2<\ldots<t_m$,
\[
\frac{1}{\Lambda_{d, \alpha}(N)}\left(\int_0^{t_1N}\left(\eta_s(\mathbf{0})-p\right)ds, \ldots, \int_0^{t_mN}\left(\eta_s(\mathbf{0})-p\right)ds\right)
\]
converges weakly to
\[
C_{18}(d, \alpha)\left(\mathcal{W}_{t_1}, \mathcal{W}_{t_2}, \ldots, \mathcal{W}_{t_m}\right)
\]
as $N\rightarrow+\infty$, where
\[
C_{18}(d, \alpha)=
\begin{cases}
\sqrt{2f_1^{d, \frac{d}{2}}(\mathbf{0})\mathcal{C}_{d, \frac{d}{2}}p(1-p)} & \text{~if~}\alpha=\frac{d}{2},\\
\sqrt{2\int_0^{+\infty}\int_0^{+\infty}p_{u+r}^{d, \alpha}(\mathbf{0}, \mathbf{0})dudr\mathcal{C}_{d, \alpha}p(1-p)}
& \text{~if~}\alpha<\frac{d}{2}.
\end{cases}
\]
\end{lemma}

\proof[Proof of Lemma \ref{lemma 4.8}]
By \eqref{equ 4.34 martingale decomposition} and Lemma \ref{lemma 4.7}, we only need to show that
\[
\left\{\frac{1}{\Lambda_{d, \alpha}(N)}\widetilde{M}_{tN}^N:~0\leq t\leq T\right\}
\]
converges weakly, with respect to the Skorohod topology, to $\{C_{18}(d, \alpha)\mathcal{W}_t:~0\leq t\leq T\}$ as $N\rightarrow+\infty$. Let
$\{\langle \widetilde{M}^N\rangle_t\}_{t\geq 0}$ be the quadratic variation process of $\{\widetilde{M}^N_t\}_{t\geq 0}$, then by Theorem 1.4 in Chapter 7 of \cite{Ethier1986}, to complete the proof we only need to show that
\begin{equation}\label{equ 4.37}
\lim_{N\rightarrow+\infty}\frac{1}{\Lambda_{d, \alpha}^2(N)}\langle\widetilde{M}^N\rangle_{tN}=C_{18}^2(d, \alpha)t
\end{equation}
in probability for any $t\geq 0$.

By Dynkin's martingale formula,
\begin{align*}
\langle \widetilde{M}^N\rangle_t&=\int_0^t\sum_{x\in \mathbb{Z}^d}\sum_{y\neq \mathbf{0}}\frac{1}{\|y\|_2^{d+\alpha}}\left(\mathcal{I}_N(\eta_s^{x, y})-\mathcal{I}_N(\eta_s)\right)^2ds\\
&=\int_0^t\sum_{x\in \mathbb{Z}^d}\sum_{y\neq \mathbf{0}}\frac{1}{\|y\|_2^{d+\alpha}}\Upsilon_N^2(x)\left(\eta_s(y)-\eta_s(x)\right)^2ds.
\end{align*}
Then, according to \eqref{equ 1.13}, \eqref{equ 4.20} and the invariance of $\nu_p^{d, \alpha}$,
\[
\hat{\mathbb{E}}_{\nu_p^{d, \alpha}}\langle \widetilde{M}^N\rangle_t=2p(1-p)t\mathcal{C}_{d, \alpha}\sum_{x\in \mathbb{Z}^d}\Upsilon^2_N(x).
\]
When $\alpha<\frac{d}{2}$, by Kolmogorov-Chapman equation,
\begin{align*}
\sum_{x\in \mathbb{Z}^d}\Upsilon^2_N(x)=\int_0^{+\infty}\int_0^{+\infty}e^{-\frac{s_1+s_2}{N}}p_{s_1+s_2}^{d, \alpha}(\mathbf{0}, \mathbf{0})ds_1ds_2
\end{align*}
and hence
\[
\lim_{N\rightarrow+\infty}\sum_{x\in \mathbb{Z}^d}\Upsilon^2_N(x)=\int_0^{+\infty}\int_0^{+\infty}p_{s_1+s_2}^{d, \alpha}(\mathbf{0}, \mathbf{0})ds_1ds_2.
\]
Note that $\int_0^{+\infty}\int_0^{+\infty}p_{s_1+s_2}^{d, \alpha}(\mathbf{0}, \mathbf{0})ds_1ds_2<+\infty$ when $\alpha<\frac{d}{2}$ by Lemma \ref{lemma LCLT}. Consequently, when $\alpha<\frac{d}{2}$,
\begin{align}\label{equ 4.38}
\lim_{N\rightarrow+\infty}\frac{1}{\Lambda_{d, \alpha}^2(N)}\hat{\mathbb{E}}_{\nu_p^{d, \alpha}}\langle\widetilde{M}^N\rangle_{tN}
&=2p(1-p)t\mathcal{C}_{d, \alpha}\int_0^{+\infty}\int_0^{+\infty}p_{s_1+s_2}^{d, \alpha}(\mathbf{0}, \mathbf{0})ds_1ds_2\notag\\
&=C_{18}^2(d, \alpha)t.
\end{align}
When $\alpha=\frac{d}{2}$, for any given $M>0$,
\begin{align*}
\frac{1}{\log N}\sum_{x\in \mathbb{Z}^d}\Upsilon_N^2(x)&=\frac{1}{\log N}\int_0^{+\infty}\int_0^{+\infty}e^{-\frac{s_1+s_2}{N}}
p_{s_1+s_2}^{d, \frac{d}{2}}(\mathbf{0}, \mathbf{0})ds_1ds_2\\
&=\frac{1}{\log N}\int_0^{+\infty}re^{-\frac{r}{N}}p_r^{d, \frac{d}{2}}(\mathbf{0}, \mathbf{0})dr\\
&=\frac{1}{\log N}\int_M^{+\infty}re^{-\frac{r}{N}}p_r^{d, \frac{d}{2}}(\mathbf{0}, \mathbf{0})dr+o(1).
\end{align*}
By Lemma \ref{lemma LCLT}, for any $\epsilon>0$, there exists $M=M(\epsilon)<+\infty$ independent of $s$ such that
\[
f_1^{d, \frac{d}{2}}(\mathbf{0})(1-\epsilon)s^{-2}\leq p_s^{d, \frac{d}{2}}(\mathbf{0}, \mathbf{0})\leq (1+\epsilon)s^{-2}f_1^{d, \frac{d}{2}}(\mathbf{0})
\]
when $s\geq M(\epsilon)$. Then,
\begin{align*}
\int_{M(\epsilon)}^{+\infty}re^{-\frac{r}{N}}p_r^{d, \frac{d}{2}}(\mathbf{0}, \mathbf{0})dr
&\leq (1+\epsilon)f_1^{d, \frac{d}{2}}(\mathbf{0})\int_{M(\epsilon)}^{+\infty}r^{-1}e^{-\frac{r}{N}}dr\\
&=(1+\epsilon)f_1^{d, \frac{d}{2}}(\mathbf{0})\int_{M(\epsilon)/N}^{+\infty}v^{-1}e^{-v}dv
\end{align*}
and
\[
\int_{M(\epsilon)}^{+\infty}re^{-\frac{r}{N}}p_r^{d, \frac{d}{2}}(\mathbf{0}, \mathbf{0})dr\geq
(1-\epsilon)f_1^{d, \frac{d}{2}}(\mathbf{0})\int_{M(\epsilon)/N}^{+\infty}v^{-1}e^{-v}dv.
\]
By L'H\^{o}pital's rule, for any $M>0$,
\begin{align*}
\lim_{u\rightarrow+\infty}\frac{\int_{M/u}^{+\infty}v^{-1}e^{-v}dv}{\log u}
&=\lim_{u\rightarrow+\infty}\frac{-(M/u)^{-1}e^{-M/u}M(-1)u^{-2}}{u^{-1}}\\
&=1.
\end{align*}
As a result,
\[
\limsup_{N\rightarrow+\infty}\frac{1}{\log N}\sum_{x\in \mathbb{Z}^d}\Upsilon_N^2(x)\leq (1+\epsilon)f_1^{d, \frac{d}{2}}(\mathbf{0})
\]
and
\[
\liminf_{N\rightarrow+\infty}\frac{1}{\log N}\sum_{x\in \mathbb{Z}^d}\Upsilon_N^2(x)\geq (1-\epsilon)f_1^{d, \frac{d}{2}}(\mathbf{0}).
\]
Since $\epsilon$ is arbitrary, let $\epsilon\rightarrow 0$ in above two inequalities and then
\[
\lim_{N\rightarrow+\infty}\frac{1}{\log N}\sum_{x\in \mathbb{Z}^d}\Upsilon_N^2(x)=f_1^{d, \frac{d}{2}}(\mathbf{0}).
\]
Consequently, when $\alpha=\frac{d}{2}$,
\begin{equation}\label{equ 4.39}
\lim_{N\rightarrow+\infty}\frac{1}{\Lambda_{d, \frac{d}{2}}^2(N)}\hat{\mathbb{E}}_{\nu_p^{d, \frac{d}{2}}}\langle\widetilde{M}^N\rangle_{tN}
=2p(1-p)t\mathcal{C}_{d, \frac{d}{2}}f_1^{d, \frac{d}{2}}(\mathbf{0})=C_{18}^2\left(d, \frac{d}{2}\right)t.
\end{equation}
By \eqref{equ 4.38} and \eqref{equ 4.39}, to check \eqref{equ 4.37}, we only need to show that
\begin{equation}\label{equ 4.40}
\lim_{N\rightarrow+\infty}{\rm Var}_{\nu_p^{d, \alpha}}\left(\frac{1}{\Lambda_{d, \alpha}^2(N)}\langle\widetilde{M}^N\rangle_{tN}\right)=0.
\end{equation}
According to the bilinear property of the covariance and the dominated convergence theorem, to prove \eqref{equ 4.40}, we only need to check a version of \eqref{equ 4.23} in the case $\alpha\leq \frac{d}{2}$. The proof of \eqref{equ 4.23} given in Subsection \ref{subsection 4.1} only relies on the transience of the long-range random walk, which also applies in the case $\alpha\leq \frac{d}{2}$. Since \eqref{equ 4.40} holds, the proof is complete.
\qed

At last, we prove Theorem \ref{theorem 2.1 main result d=1, 2} in the case $\alpha\leq \frac{d}{2}$. We only give an outline since the proof is similar to that given in \cite{Xue2024} for the nearest-neighbor model on $\mathbb{Z}^4$ and $\mathbb{Z}^5$.

\proof[Proof of Theorem \ref{theorem 2.1 main result d=1, 2} in the case $\alpha\leq \frac{d}{2}$]
Let $d\in \{1, 2, 3\}$ and $\alpha\leq \frac{d}{2}$. By Lemma \ref{lemma 4.8}, we only need to show that
\[
\left\{\frac{1}{\Lambda_{d, \alpha}(N)}\int_0^{tN}\left(\eta_s(\mathbf{0})-p\right)ds:~0\leq t\leq T\right\}_{N\geq 1}
\]
are tight under the uniform topology. To check this tightness, we only need to show that there exist $a, b, c>0$ independent of $s, t, N$ such that
\begin{equation}\label{equ 4.41}
\hat{\mathbb{E}}_{\nu_p^{d, \alpha}}\left(\left|\frac{1}{\Lambda_{d, \alpha}(N)}\int_{sN}^{tN}\left(\eta_r(\mathbf{0})-p\right)dr\right|^{a}\right)
\leq c|t-s|^{b+1}
\end{equation}
for any $s, t$ such that $|t-s|\leq 1$ and $N\geq 1$.

Now we show that \eqref{equ 4.41} holds with $a=4, b=1$ and some $c>0$. By Proposition \ref{proposition voter dual}, it is not difficult to show that
\begin{align}\label{equ 4.42}
&\hat{\mathbb{E}}_{\nu_p^{d, \alpha}}\left(\left|\int_{0}^{tN}\left(\eta_r(\mathbf{0})-p\right)dr\right|^{4}\right)\notag\\
&\leq 72\mathcal{C}_{d, \alpha}\left(\int_0^{tN}d\theta\int_0^{\theta}dv\int_0^{+\infty}p_{r+v}^{d, \alpha}(\mathbf{0}, \mathbf{0})dr\right)^2.
\end{align}
Reference \cite{Xue2024} gives a detailed proof of an analogue of \eqref{equ 4.42} for the nearest-neighbor model, which also applies in the long-range case. So we omit how to check \eqref{equ 4.42}. When $\alpha<\frac{d}{2}$, by Lemma \ref{lemma LCLT},
\[
\int_0^{\theta}dv\int_0^{+\infty}p_{r+v}^{d, \alpha}(\mathbf{0}, \mathbf{0})dr
\leq \int_0^{+\infty}dv\int_0^{+\infty}p_{r+v}^{d, \alpha}(\mathbf{0}, \mathbf{0})dr<+\infty
\]
and hence
\[
\hat{\mathbb{E}}_{\nu_p^{d, \alpha}}\left(\left|\int_{0}^{tN}\left(\eta_r(\mathbf{0})-p\right)dr\right|^{4}\right)
\leq O(1)t^2N^2.
\]
Consequently, Equation \eqref{equ 4.41} holds with $a=4, b=1$.

When $\alpha=\frac{d}{2}$ and $tN\geq 1$, by Lemma \ref{lemma LCLT} and \eqref{equ 4.42},
\begin{align*}
\hat{\mathbb{E}}_{\nu_p^{d, \frac{d}{2}}}\left(\left|\int_{0}^{tN}\left(\eta_r(\mathbf{0})-p\right)dr\right|^{4}\right)
&\leq O(1)\left(\int_1^{tN}d\theta\int_1^{tN}dv\int_0^{+\infty}(v+r)^{-2}dr\right)^2\\
&=O(1)\left(tN\log N\right)^2
\end{align*}
and hence
\[
\hat{\mathbb{E}}_{\nu_p^{d, \frac{d}{2}}}\left(\left|\frac{1}{\sqrt{N\log N}}\int_{0}^{tN}\left(\eta_r(\mathbf{0})-p\right)dr\right|^{4}\right)
\leq O(1)t^2.
\]

When $tN\leq 1$,
\[
\hat{\mathbb{E}}_{\nu_p^{d, \frac{d}{2}}}\left(\left|\int_{0}^{tN}\left(\eta_r(\mathbf{0})-p\right)dr\right|^{4}\right)\leq t^4N^4\leq t^2N^2
\]
and hence
\[
\hat{\mathbb{E}}_{\nu_p^{d, \frac{d}{2}}}\left(\left|\frac{1}{\sqrt{N\log N}}\int_{0}^{tN}\left(\eta_r(\mathbf{0})-p\right)dr\right|^{4}\right)\leq \frac{t^2}{(\log N)^2}.
\]
In conclusion, when $\alpha=\frac{d}{2}$, Equation \eqref{equ 4.41} also holds with $a=4, b=1$. Since \eqref{equ 4.41} holds, the proof is complete.
\qed

\section{Proof of Theorem \ref{theorem 2.2 main result d=3}}\label{section five}
In this section, we prove Theorem \ref{theorem 2.2 main result d=3}. The proof utilizes a similar argument to that given in Subsection \ref{subsection 4.1} and hence we only give an outline to avoid repeating many similar details. We first introduce some notations and definitions. Defined as in Section \ref{section two},
\[
\Lambda_{3, \alpha}(N)=
\begin{cases}
N^{\frac{3}{4}} & \text{~if~} \alpha>2,\\
N^{\frac{3}{4}}\left(\log N\right)^{-\frac{3}{4}} & \text{~if~} \alpha=2.
\end{cases}
\]
Let $\hat{N}_t^{x, y}$, $v(t, x)$, $V^t(s, \eta)$ and $M_s^t$ be defined as in Subsection \ref{subsection 4.1} except that the range of $\alpha$ is $[2, +\infty)$, then we have a $d=3, \alpha\geq 2$ version of \eqref{equ 4.3 martingale decomposition}. The following lemma is an analogue of Lemma \ref{lemma 4.1}.
\begin{lemma}\label{lemma 5.1}
Let $d=3$ and $\alpha\geq 2$.
For any integer $m\geq 1$ and $0<t_1<t_2<\ldots<t_m$, if $\eta_0$ is distributed with $\nu_p^{3, \alpha}$, then
\[
\frac{1}{\Lambda_{3, \alpha}(N)}\left(V^{t_1N}(0, \eta_0), V^{t_2N}(0, \eta_0), \ldots, V^{t_mN}(0, \eta_0)\right)
\]
converges weakly to a $\mathbb{R}^m$-valued Gaussian random variable  $\left(\Gamma_{t_1}^{3, \alpha}, \ldots, \Gamma_{t_m}^{3, \alpha}\right)$ as $N\rightarrow+\infty$.
\end{lemma}
To prove Lemma \ref{lemma 5.1} and compute the variance of $\Gamma_{t}^{3, \alpha}$, we still need a fluctuation theorem of the empirical density field of the model in the case $d=3, \alpha\geq 2$. In detail, we define
\[
\mathcal{Y}_t^{\epsilon, 3, \alpha}(du)=\tilde{h}(\epsilon, \alpha)\sum_{x\in \mathbb{Z}^3}\left(\eta_{t\epsilon^{-1}}(x)-p\right)\delta_{m(\epsilon, \alpha)x}(du),
\]
where
\[
\tilde{h}(\epsilon, \alpha)=
\begin{cases}
\epsilon^{\frac{5}{4}} & \text{~if~} \alpha>2,\\
\epsilon^{\frac{5}{4}}(\log \frac{1}{\epsilon})^{-\frac{3}{4}} & \text{~if~} \alpha=2.
\end{cases}
\]
and
\[
m(\epsilon, \alpha)=
\begin{cases}
\epsilon^{\frac{1}{2}} & \text{~if~} \alpha>2,\\
\epsilon^{\frac{1}{2}}(\log\frac{1}{\epsilon})^{-\frac{1}{2}} & \text{~if~} \alpha=2.
\end{cases}
\]
Note that, to give the fluctuation density field a simpler expression in the case $\alpha=2$, $\epsilon$ in $\mathcal{Y}_t^{\epsilon, 3, \alpha}$ plays the role as that of $\epsilon^\alpha$ in $\mathcal{Y}_t^\epsilon$ defined in Subsection \ref{subsection 4.1}.

When $d=3$ and $\alpha\geq 2$, we have the following generalized O-U process which is an analogue of \eqref{equ 4.5 generalized O-U process},
\begin{equation}\label{equ 5.1 genralized O-U d=3 alpha geq 2}
d\mathcal{Y}_t^{3, \alpha}=\hat{\mathcal{L}}^{3, \alpha}\mathcal{Y}_t^{3, \alpha}dt+Cd\mathcal{B}_t.
\end{equation}
Similarly, Equation \eqref{equ 5.1 genralized O-U d=3 alpha geq 2} has an unique stationary distribution $\pi^{C, 3, \alpha}$ such that
$\mu(f)$ follows the normal distribution with mean zero and variance
\[
C^2\int_0^{+\infty}\int_{\mathbb{R}^d}\left(\mathcal{P}^{3, \alpha}_sf(u)\right)^2dsdu
\]
for any Schwartz function $f$ when the random tempered distribution $\mu$ is distributed with $\pi^{C, 3, \alpha}$, where $\{\mathcal{P}^{3, \alpha}_t\}_{t\geq 0}$ is the Markov semigroup with respect to $\hat{\mathcal{L}}^{3, \alpha}$. We have the following analogue of Lemma \ref{lemma 4.2 fluctuation of long-range voter model}.
\begin{lemma}\label{lemma 5.2 fluctuation of long-range voter model d=3 alpha geq 2}
For $\alpha\geq 2$ and any $T>0$, as $\epsilon\rightarrow 0$,
\[
\left\{\mathcal{Y}_t^{\epsilon, 3, \alpha}:~0\leq t\leq T\right\}
\]
converges weakly, with respect to the Skorohod topology, to
\[
\left\{\mathcal{Y}_t^{3, \alpha}:~0\leq t\leq T\right\},
\]
where $\left\{\mathcal{Y}_t^{3, \alpha}\right\}_{t\geq 0}$ is the solution to Equation \eqref{equ 5.1 genralized O-U d=3 alpha geq 2} with
\[
C=\sqrt{2p(1-p)\sum_{y\neq \mathbf{0}}\frac{\Phi^{3, \alpha}(y)}{\|y\|^{3+\alpha}}}
\]
and $\mathcal{Y}_0^{3, \alpha}$ being distributed with $\pi^{C, 3, \alpha}$.
\end{lemma}

The proof of Lemma \ref{lemma 5.2 fluctuation of long-range voter model d=3 alpha geq 2} utilizes a similar Dynkin's martingale strategy to that of Lemma \ref{lemma 4.2 fluctuation of long-range voter model}, where a crucial step is to check the following claim, which is an analogue of \eqref{equ 4.10}. For $\alpha\geq 2$, there exists a constant $\hat{C}=\hat{C}(\alpha)\in (0, +\infty)$ independent of $x\in \mathbb{Z}^3$ such that
\begin{equation}\label{equ 5.2}
\lim_{\|x\|_2\rightarrow+\infty}\frac{\int_0^{+\infty}p_s^{3, \alpha}(\mathbf{0}, x)ds}{\tilde{Q}(\|x\|_2)}=\hat{C},
\end{equation}
where
\[
\tilde{Q}(u)=
\begin{cases}
u^{-1} & \text{~if~}\alpha>2,\\
\left(u\log u\right)^{-1} & \text{~if~}\alpha=2.
\end{cases}
\]
Equation \eqref{equ 5.2} follows from a similar argument to that in the proof of \eqref{equ 4.10}, which utilizes the inversion formula. By Lemma \ref{lemma 5.2 fluctuation of long-range voter model d=3 alpha geq 2} and an argument similar to that in the proof of Lemma \ref{lemma 4.1}, we have that Lemma \ref{lemma 5.1} holds with
\[
\Gamma_t^{3, \alpha}=\mathcal{Y}_t^{3, \alpha}(J_t^{3, \alpha}),
\]
where
\[
J_t^{3, \alpha}(u)=\int_0^tf_s^{3, \alpha}(u)ds.
\]
Consequently, according to a computation similar to that in the proof of Lemma \ref{lemma 4.6}, for any $0<t_1<t_2$,
\begin{align}\label{equ 5.3}
{\rm Cov}\left(\Gamma_{t_1}^{3, \alpha}, \Gamma_{t_2}^{3, \alpha}\right)
&=\frac{f_1^{3, \alpha}(\mathbf{0})2^3p(1-p)\mathcal{C}_{3, \alpha}}{(3-2)(2\times 2-3)(3\times 2-3)}\left((t_1+t_2)^{3-\frac{3}{2}}-t_1^{3-\frac{3}{2}}-t_2^{3-\frac{3}{2}}\right) \notag\\
&=\frac{8f_1^{3, \alpha}(\mathbf{0})p(1-p)\mathcal{C}_{3, \alpha}}{3}\left((t_1+t_2)^{\frac{3}{2}}-t_1^{\frac{3}{2}}-t_2^{\frac{3}{2}}\right).
\end{align}
Now we deal with the martingale term $M_{tN}^{tN}$. We have the following analogue of Lemma \ref{lemma 4.3}.
\begin{lemma}\label{lemma 5.3}
Let $d=3$ and $\alpha\geq 2$. For any integer $m\geq 1$ and $0<t_1<t_2<\ldots<t_m$, if $\eta_0$ is distributed with $\nu_p^{3, \alpha}$, then
\[
\frac{1}{\Lambda_{3, \alpha}(N)}\left(M_{t_1N}^{t_1N}, \ldots, M_{t_mN}^{t_mN}\right)
\]
converges weakly, as $N\rightarrow+\infty$, to a $\mathbb{R}^m$-valued Gaussian random variable $(\Xi_{t_1}^{3, \alpha},\ldots, \Xi_{t_m}^{3, \alpha})$.
\end{lemma}
To prove Lemma \ref{lemma 5.3}, we define $\mathcal{U}^{N, 3, \alpha}$ as the random measure on $[0, T]\times \mathbb{R}^3$ such that
\begin{align*}
&\mathcal{U}^{N, 3, \alpha}(H)= \\
&\tilde{k}_\alpha(N)\int_{\mathbb{R}^3}\sum_{y\in \mathbb{Z}^3\setminus\{h_\alpha(N)u_{\alpha, N}\}}\left(\int_0^T H(s, u)\left(\eta_{sN-}(y)-\eta_{sN-}(h_\alpha(N)u_{\alpha, N})\right)d\hat{N}_{Ns}^{(h_\alpha(N)u_{\alpha, N}), y}\right)du\\
&=\tilde{k}_\alpha(N)\sum_{x\in \mathbb{Z}^3}\sum_{y\in \mathbb{Z}^3\setminus\{x\}}\int_0^T\left(\int_{\frac{x}{h_\alpha(N)}+\left(-\frac{1}{2h_\alpha(N)}, \frac{1}{2h_\alpha(N)}\right]^3}H(s, u)du\right)(\eta_{sN-}(y)-\eta_{sN-}(x))d\hat{N}_{Ns}^{x, y}
\end{align*}
for any $H\in C_c([0, T]\times \mathbb{R}^3)$, where
\[
h_\alpha(N)=
\begin{cases}
\sqrt{N} & \text{~if~}\alpha>2,\\
\sqrt{N\log N} & \text{~if~}\alpha=2
\end{cases}
\]
defined as in Subsection \ref{subsection 1.2},
\[
\tilde{k}_\alpha(N)=
\begin{cases}
N^{\frac{1}{4}} & \text{~if~}\alpha>2,\\
N^{\frac{1}{4}}\left(\log N\right)^{\frac{3}{4}} & \text{~if~}\alpha=2
\end{cases}
\]
and $u_{\alpha, N}$ is the element in $\mathbb{Z}^3/(h_\alpha(N))$ such that
\[
u-u_{\alpha, N}\in \left(-\frac{1}{2h_\alpha(N)}, \frac{1}{2h_\alpha(N)}\right]^3.
\]
Then we have the following analogue of \eqref{equ 4.19},
\begin{equation}\label{equ 5.4}
\frac{1}{\Lambda_{3, \alpha}(N)}M_{tN}^{tN}=\mathcal{U}^{N, 3, \alpha}(b_t^{N, 3, \alpha}),
\end{equation}
where
\[
b_t^{N, 3, \alpha}(s, u)=\frac{\left(h_\alpha(N)\right)^3}{N}\sum_{x\in \mathbb{Z}^3}v(N(t-s), x)1_{\left\{u\in \frac{x}{h_\alpha(N)}+\left(-\frac{1}{2h_\alpha(N)}, \frac{1}{2h_\alpha(N)}\right]^3, s\leq t\right\}}.
\]
Let $\mathcal{U}$ be defined as in Subsection \ref{subsection 4.1} except that the constraints of $(d, \alpha)$ are $d=3$ and $\alpha\geq 2$, then we have the following analogue of Lemma \ref{lemma 4.4}.
\begin{lemma}\label{lemma 5.4}
Let $d=3, \alpha\geq 2$ and $\eta_0$ be distributed with $\nu_p^{3, \alpha}$. For any
\[
H\in C_c^\infty([0, T]\times\mathbb{R}^3),
\]
$\mathcal{U}^{N, 3, \alpha}(H)$ converges weakly to $\mathcal{U}(H)$ as $N\rightarrow+\infty$.
\end{lemma}
Lemma \ref{lemma 5.4} follows from a similar argument to that in the proof of Lemma \ref{lemma 4.4}, where Proposition \ref{proposition voter dual} plays the key role. By Lemmas \ref{lemma LCLT}, \ref{lemma 5.4} and an argument similar to that in the proof of Lemma \ref{lemma 4.3}, we have that Lemma \ref{lemma 5.3} holds with
\[
\Xi_t^{3, \alpha}=\mathcal{U}(b_t^{3, \alpha}),
\]
where $b_t^{3, \alpha}=J_{t-s}^{3, \alpha}(u)1_{\{s\leq t\}}$. Furthermore, according to a computation similar to that in the proof of Lemma \ref{lemma 4.6}, we have
\begin{align}\label{equ 5.5}
&{\rm Cov}\left(\Xi_{t_1}^{3, \alpha}, \Xi_{t_2}^{3, \alpha}\right)\\
&=\frac{f_1^{3, \alpha}(\mathbf{0})2^3p(1-p)\mathcal{C}_{3, \alpha}}{(3-2)(2\times 2-3)(3\times 2-3)}
\left(2t_2^{3-\frac{3}{2}}+2t_1^{3-\frac{3}{2}}-(t_1+t_2)^{3-\frac{3}{2}}-(t_2-t_1)^{3-\frac{3}{2}}\right)\notag\\
&=\frac{f_1^{3, \alpha}(\mathbf{0})8p(1-p)\mathcal{C}_{3, \alpha}}{3}
\left(2t_2^{\frac{3}{2}}+2t_1^{\frac{3}{2}}-(t_1+t_2)^{\frac{3}{2}}-(t_2-t_1)^{\frac{3}{2}}\right)\notag
\end{align}
for any $0\leq t_1\leq t_2$.

Now we deal with the joint distribution of $M_{tN}^{tN}$ and $V^{tN}(0, \eta_0)$. We have the following analogue of Lemma \ref{lemma 4.5}.
\begin{lemma}\label{lemma 5.5}
Let $d=3, \alpha\geq 2$ and $\eta_0$ be distributed with $\nu_p^{3, \alpha}$. For any $0\leq t_1<t_2<\ldots<t_m$,
the joint distribution of
\[
\frac{1}{\Lambda_{3, \alpha}(N)}\left(V^{t_1N}(0, \eta_0), V^{t_2N}(0, \eta_0), \ldots, V^{t_mN}(0, \eta_0)\right)
\]
and
\[
\frac{1}{\Lambda_{2, \alpha}(N)}\left(M_{t_1N}^{t_1N}, \ldots, M_{t_mN}^{t_mN}\right)
\]
converges weakly, as $N\rightarrow+\infty$, to the independent coupling of
\[
\left(\Gamma_{t_1}^{3, \alpha}, \ldots, \Gamma_{t_m}^{3, \alpha}\right)
\text{\quad and\quad}
(\Xi_{t_1}^{3, \alpha},\ldots, \Xi_{t_m}^{3, \alpha}),
\]
which are given in Lemmas \ref{lemma 5.1} and \ref{lemma 5.3} respectively.
\end{lemma}
The proof of Lemma \ref{lemma 5.5} is similar to that of Lemma \ref{lemma 4.5} and hence we omit the details. By \eqref{equ 5.3}, \eqref{equ 5.5} and Lemma \ref{lemma 5.5}, we have the following analogue of Lemma \ref{lemma 4.6}.
\begin{lemma}\label{lemma 5.6}
Let $d=3, \alpha\geq 2$ and $\eta_0$ be distributed with $\nu_p^{3, \alpha}$. For any $0\leq t_1<t_2<\ldots<t_m$,
\[
\frac{1}{\Lambda_{3, \alpha}(N)}\left(\int_0^{t_1N}\left(\eta_s(\mathbf{0})-p\right)ds, \ldots, \int_0^{t_mN}\left(\eta_s(\mathbf{0})-p\right)ds\right)
\]
converges weakly to
\[
\sqrt{\frac{16}{3}f_1^{3, \alpha}(\mathbf{0})\mathcal{C}_{3, \alpha}p(1-p)}
\left(B_{t_1}^{\frac{3}{4}}, \ldots, B_{t_m}^{\frac{3}{4}}\right)
\]
as $N\rightarrow+\infty$.
\end{lemma}

At last, we prove Theorem \ref{theorem 2.2 main result d=3}.

\proof[Proof of Theorem \ref{theorem 2.2 main result d=3}]
Let $d=3, \alpha\geq 2$ and $\eta_0$ be distributed with $\nu_p^{3, p}$.
By Lemma \ref{lemma 5.6}, we only need to show that
\[
\left\{\frac{1}{\Lambda_{3, \alpha}(N)}\int_0^{tN}\left(\eta_s(\mathbf{0})-p\right)ds:~0\leq t\leq T\right\}_{N\geq 1}
\]
are tight under the uniform topology. To check this tightness, we only need to show that there exist $a, b, c>0$ independent of $s, t, N$ such that
\begin{equation}\label{equ 5.6}
\hat{\mathbb{E}}_{\nu_p^{3, \alpha}}\left(\left|\frac{1}{\Lambda_{3, \alpha}(N)}\int_{sN}^{tN}\left(\eta_r(\mathbf{0})-p\right)dr\right|^{a}\right)
\leq c|t-s|^{b+1}
\end{equation}
for any $0<s<t$ such that $|t-s|<1$ and $N\geq 2$.
According to a computation similar to that in the proof of Theorem \ref{theorem 2.1 main result d=1, 2},
\[
\hat{\mathbb{E}}_{\nu_p^{3, \alpha}}\left(\left(\int_{0}^{t}\left(\eta_r(\mathbf{0})-p\right)dr\right)^2\right)
=2\mathcal{C}_{3, \alpha}p(1-p)\int_0^tdr\int_0^r d\theta\int_0^{+\infty}p_{\theta+v}^{3, \alpha}(\mathbf{0}, \mathbf{0})dv.
\]
By Lemma \ref{lemma LCLT}, there exists $C_5=C_5(\alpha)<+\infty$ independent of $t$ such that
\begin{equation}\label{equ 5.7}
p_t^{3, \alpha}(\mathbf{0}, \mathbf{0})\leq C_5\left(h_\alpha(t)\right)^{-3}
\end{equation}
for all $t\geq 0$. Therefore, when $\alpha>2$,
\begin{align*}
\hat{\mathbb{E}}_{\nu_p^{3, 2}}\left(\left(\int_{0}^{t}\left(\eta_r(\mathbf{0})-p\right)dr\right)^2\right)
&\leq 2\mathcal{C}_{3, \alpha}p(1-p)C_5\int_0^tdr\int_0^r d\theta\int_0^{+\infty}(\theta+v)^{-\frac{3}{2}}dv\\
&=\frac{16}{3}\mathcal{C}_{3, \alpha}p(1-p)C_5t^{\frac{3}{2}}.
\end{align*}
Then, when $\alpha>2$, Equation \eqref{equ 5.6} holds with $a=2$, $b=\frac{1}{2}$ and
\[
c=\frac{16}{3}\mathcal{C}_{3, \alpha}p(1-p)C_5.
\]
When $\alpha=2$, for $t<1$ and $N\leq \frac{1000}{t}$,
\begin{align*}
\int_0^{tN}dr\int_0^r d\theta\int_0^{+\infty}p_{\theta+v}^{3, 2}(\mathbf{0}, \mathbf{0})dv
\leq C_6\int_0^{tN}dr\int_0^r1d\theta\leq C_6t^2N^2,
\end{align*}
where
\[
C_6=\int_0^{+\infty}p_{v}^{3, 2}(\mathbf{0}, \mathbf{0})dv<+\infty
\]
according to the transience. Hence, when $tN\leq 1000$,
\[
\hat{\mathbb{E}}_{\nu_p^{3, 2}}\left(\left(\frac{1}{\Lambda_{3, 2}(N)}\int_{0}^{tN}\left(\eta_r(\mathbf{0})-p\right)dr\right)^2\right)
\leq C_7 t^2N^2\frac{\left(\log N\right)^{\frac{3}{2}}}{N^{\frac{3}{2}}}\leq C_8 t^{\frac{3}{2}}\left(\log \frac{1000}{t}\right)^{\frac{3}{2}},
\]
where $C_7=2\mathcal{C}_{3, 2}p(1-p)C_6$ and $C_8=\sqrt{1000}C_7$. Since $(\log\frac{1000}{t})^{\frac{3}{2}}/(t^{-0.1})\rightarrow 0$ as $t\rightarrow 0$, there exists $C_9<+\infty$ such that
\[
\left(\log\frac{1000}{t}\right)^{\frac{3}{2}}\leq C_9 t^{-0.1}
\]
for all $t\in (0, 1)$. Hence, when $t\leq 1$ and $tN\leq 1000$,
\[
\hat{\mathbb{E}}_{\nu_p^{3, 2}}\left(\left(\frac{1}{\Lambda_{3, 2}(N)}\int_{0}^{tN}\left(\eta_r(\mathbf{0})-p\right)dr\right)^2\right)
\leq  C_{10}t^{1.4},
\]
where $C_{10}=C_{8}C_9$. By \eqref{equ 5.7} and integration-by-parts formula, it is not difficult to show that there exists $C_{11}<+\infty$ independent of $s$ such that
\begin{equation}\label{equ 5.8}
\int_0^sdr\int_0^r d\theta\int_0^{+\infty}p_{\theta+v}^{3, 2}(\mathbf{0}, \mathbf{0})dv\leq C_{11}\frac{s^{\frac{3}{2}}}{\left(\log s\right)^{\frac{3}{2}}}
\end{equation}
for any $s\geq 1000$. We check \eqref{equ 5.8} in Appendix \ref{Appendix A.6}. By \eqref{equ 5.8}, when $t<1$ and $tN\geq 1000$,
\begin{align*}
\hat{\mathbb{E}}_{\nu_p^{3, 2}}\left(\left(\frac{1}{\Lambda_{3, 2}(N)}\int_{0}^{tN}\left(\eta_r(\mathbf{0})-p\right)dr\right)^2\right)
&\leq C_{12}\frac{\left(\log N\right)^{\frac{3}{2}}}{N^{\frac{3}{2}}}\frac{(tN)^{\frac{3}{2}}}{\left(\log(tN)\right)^{\frac{3}{2}}}\\
&=C_{12}t^{1.5}\left(\frac{\log N}{\log N-\log\frac{1}{t}}\right)^{\frac{3}{2}},
\end{align*}
where $C_{12}=2C_{11}\mathcal{C}_{3, \alpha}p(1-p)$. Since $N\geq \frac{1000}{t}$,
\[
\left(\frac{\log N}{\log N-\log\frac{1}{t}}\right)^{\frac{3}{2}}\leq \left(\frac{1}{\log 1000}\right)^{\frac{3}{2}}\left(\log \frac{1000}{t}\right)^{\frac{3}{2}}\leq C_{13}t^{-0.1},
\]
where $C_{13}=\left(\frac{1}{\log 1000}\right)^{\frac{3}{2}}C_9$. Hence, when $t<1$ and $tN\geq 1000$,
\[
\hat{\mathbb{E}}_{\nu_p^{3, 2}}\left(\left(\frac{1}{\Lambda_{3, 2}(N)}\int_{0}^{tN}\left(\eta_r(\mathbf{0})-p\right)dr\right)^2\right)
\leq  C_{14}t^{1.4},
\]
where $C_{14}=C_{12}C_{13}$. In conclusion, when $\alpha=2$, Equation \eqref{equ 5.6} holds with $a=2$, $b=0.4$ and
\[
c=\max\{C_{10}, C_{14}\}.
\]
Since \eqref{equ 5.6} holds for both cases, the proof is complete.
\qed

\section{Proofs of Theorems \ref{theorem 2.3 main result d=4} and \ref{theorem 2.4 main result d geq 5}}\label{section six}
In this section, we prove Theorems \ref{theorem 2.3 main result d=4} and \ref{theorem 2.4 main result d geq 5}. Our proofs utilize an argument similar to that given in Subsection \ref{subsection 4.2}, so we only give an outline to avoid repeating similar details. Let $\Upsilon_N, \mathcal{I}_N, \widetilde{M}_t^N$ be defined as in Subsection \ref{subsection 4.2} except that the range of $d$ is $\{4, 5, \ldots\}$, then we have a version of \eqref{equ 4.34 martingale decomposition} in the case $d\geq 4$. The following lemma is an analogue Lemma \ref{lemma 4.7}.

\begin{lemma}\label{lemma 6.1}
Let $d\geq 4$ and $\eta_0$ be distributed with $\nu_p^{d, \alpha}$. For any $t\geq 0$,
\[
\lim_{N\rightarrow+\infty}\frac{\frac{1}{N}\int_0^{tN}\mathcal{I}_N(\eta_s)ds+\mathcal{I}_N(\eta_0)-\mathcal{I}_N(\eta_{tN})}{\Lambda_{d, \alpha}(N)}=0
\]
in probability.
\end{lemma}

\proof[Proof of Lemma \ref{lemma 6.1}]
By Lemma \ref{lemma LCLT}, when $d=4$ and $\alpha>2$,
\[
p_s^{4, \alpha}(\mathbf{0}, \mathbf{0})=O(1)s^{-2}
\]
for sufficiently large $s$. Hence, the proof of Lemma \ref{lemma 6.1} in the case $d=4, \alpha>2$ is the same as that of Lemma \ref{lemma 4.7} in the case $\alpha=\frac{d}{2}$.

By Lemma \ref{lemma LCLT}, when
\[
(d, \alpha)\in \{4\}\times (0, 2)\bigcup \{5, 6, \ldots\}\times (0, +\infty),
\]
there exists $\varepsilon>0$ such that
\[
p_s^{d, \alpha}(\mathbf{0}, \mathbf{0})\leq O(1)s^{-(2+\varepsilon)}
\]
for sufficiently large $s$. Hence, the proof of Lemma \ref{lemma 6.1} in this case is the same as that of Lemma \ref{lemma 4.7} in the case $\alpha<\frac{d}{2}$.

We only need to deal with the case $d=4, \alpha=2$. According to Lemma \ref{lemma LCLT} and an argument similar to that in the proof of Lemma \ref{lemma 4.7},
\begin{align*}
\frac{1}{N}\hat{\mathbb{E}}_{\nu_p^{4, 2}}\left(\mathcal{I}_N^2(\eta_0)\right)
&\leq \frac{O(1)}{N}+\frac{O(1)}{N}\int_{10}^{+\infty}re^{-\frac{r}{N}}\left(\int_r^{+\infty}\frac{1}{s^2(\log s)^2}ds\right)dr\\
&\leq \frac{O(1)}{N}+\frac{O(1)}{N}\int_{10}^{+\infty}re^{-\frac{r}{N}}\frac{1}{r(\log r)^2}dr\\
&=\frac{O(1)}{N}+\frac{O(1)}{N}\int_{10}^{+\infty}e^{-\frac{r}{N}}\frac{1}{(\log r)^2}dr\\
&=\frac{O(1)}{N}+O(1)\int^{+\infty}_{\frac{10}{N}}e^{-u}\frac{1}{(\log (Nu))^2}du.
\end{align*}
For any $\epsilon>0$,
\begin{align*}
\int^{+\infty}_{\frac{10}{N}}e^{-u}\frac{1}{(\log (Nu))^2}du
\leq \frac{1}{(\log (N\epsilon))^2}\int_\epsilon^{+\infty}e^{-u}du+\frac{1}{\left(\log 10\right)^2}\int_0^\epsilon e^{-u}du
\end{align*}
for sufficiently large $N$ and hence
\[
\limsup_{N\rightarrow+\infty}\int^{+\infty}_{\frac{10}{N}}e^{-u}\frac{1}{(\log (Nu))^2}du\leq \frac{1}{\left(\log 10\right)^2}\int_0^\epsilon e^{-u}du.
\]
Since $\epsilon$ is arbitrary, let $\epsilon\rightarrow 0$ in the above inequality and then
\[
\lim_{N\rightarrow+\infty}\frac{1}{N}\hat{\mathbb{E}}_{\nu_p^{4, 2}}\left(\mathcal{I}_N^2(\eta_0)\right)=0,
\]
which completes the proof.
\qed

The following lemma is an analogue of Lemma \ref{lemma 4.8}.
\begin{lemma}\label{lemma 6.2}
Let $d\geq 4$ and $\eta_0$ be distributed with $\nu_p^{d, \alpha}$. For any $0\leq t_1<t_2<\ldots<t_m$,
\[
\frac{1}{\Lambda_{d, \alpha}(N)}\left(\int_0^{t_1N}\left(\eta_s(\mathbf{0})-p\right)ds, \ldots, \int_0^{t_mN}\left(\eta_s(\mathbf{0})-p\right)ds\right)
\]
converges weakly to
\[
C_{19}(d, \alpha)\left(\mathcal{W}_{t_1}, \mathcal{W}_{t_2}, \ldots, \mathcal{W}_{t_m}\right)
\]
as $N\rightarrow+\infty$, where
\[
C_{19}(d, \alpha)=
\begin{cases}
\sqrt{2\int_0^{+\infty}\int_0^{+\infty}p_{u+r}^{4, \alpha}(\mathbf{0}, \mathbf{0})dudr\mathcal{C}_{4, \alpha}p(1-p)} & \text{~if~}d=4, \alpha\leq 2,\\
\sqrt{2f_1^{4, \alpha}(\mathbf{0})\mathcal{C}_{4, \alpha}p(1-p)}
& \text{~if~}d=4, \alpha>2,\\
\sqrt{2\int_0^{+\infty}\int_0^{+\infty}p_{u+r}^{d, \alpha}(\mathbf{0}, \mathbf{0})dudr\mathcal{C}_{d, \alpha}p(1-p)}& \text{~if~}d\geq 5.
\end{cases}
\]
\end{lemma}

\proof[Proof of Lemma \ref{lemma 6.2}]
By Lemma \ref{lemma LCLT}, when
\[
(d, \alpha)\in \{4\}\times (0, 2]\bigcup \{5, 6, \ldots\}\times (0, +\infty),
\]
we have
\[
\int_0^{+\infty}\int_0^{+\infty}p_{s+v}^{d, \alpha}(\mathbf{0}, \mathbf{0})dsdv<+\infty.
\]
Hence, the proof of Lemma \ref{lemma 6.2} in this case is the same as that of Lemma \ref{lemma 4.8} in the case $\alpha<\frac{d}{2}$.

By Lemma \ref{lemma LCLT}, when $d=4, \alpha>2$,
we have
\[
p_s^{4, \alpha}(\mathbf{0}, \mathbf{0})=(1+o(1))s^{-2}f_1^{4, 2}(\mathbf{0})
\]
for sufficiently large $s$. Hence, the proof of Lemma \ref{lemma 6.2} in this case is the same as that of Lemma \ref{lemma 4.8} in the case $\alpha=\frac{d}{2}$.
\qed

At last, we prove Theorems \ref{theorem 2.3 main result d=4} and \ref{theorem 2.4 main result d geq 5}.

\proof[Proofs of Theorems \ref{theorem 2.3 main result d=4} and \ref{theorem 2.4 main result d geq 5}]
By Lemma \ref{lemma 6.2}, we only need to show that
\[
\left\{\frac{1}{\Lambda_{d, \alpha}(N)}\int_0^{tN}\left(\eta_s(\mathbf{0})-p\right)ds:~0\leq t\leq T\right\}_{N\geq 1}
\]
are tight under the uniform topology.

By Lemma \ref{lemma LCLT}, when
\[
(d, \alpha)\in \{4\}\times (0, 2]\bigcup \{5, 6, \ldots\}\times (0, +\infty),
\]
we have
\[
\int_0^{+\infty}\int_0^{+\infty}p_{s+v}^{d, \alpha}(\mathbf{0}, \mathbf{0})dsdv<+\infty.
\]
Hence, the check of the tightness in this case is the same as that in the case $d\in \{1, 2, 3\}$ and $\alpha<\frac{d}{2}$, which is given in Subsection \ref{subsection 4.2}.

By Lemma \ref{lemma LCLT}, when $d=4, \alpha>2$,
we have
\[
p_s^{4, \alpha}(\mathbf{0}, \mathbf{0})=O(1)s^{-2}
\]
for sufficiently large $s$. Hence, the check of the tightness in this case is the same as that in the case $d\in \{1, 2, 3\}$ and $\alpha=\frac{d}{2}$, which is given in Subsection \ref{subsection 4.2}.
\qed

\appendix{}
\section{Appendix}
\subsection{Proof of \eqref{equ 4.8}}\label{Appendix A.1}
Now we check \eqref{equ 4.8}.

\proof[Proof of \eqref{equ 4.8}] Assuming that $(d, \alpha)$ satisfies \eqref{equ 4.4}. By Kolmogorov-Chapman equation,
\[
\int_{\mathbb{R}^d}J_t^2(u)du=\int_0^t\int_0^tf_{s+v}^{d, \alpha}(\mathbf{0})dsdv
\]
and
\begin{align}\label{equ A.1}
\int_{\mathbb{R}^d}\left(\tilde{J}^N_t(u)\right)^2du&=\sum_{x\in \mathbb{Z}^d}N^{-\frac{d}{\alpha}}\left(\tilde{J}^N_t\left(\frac{x}{N^{\frac{1}{\alpha}}}\right)\right)^2\notag\\
&=\sum_{x\in \mathbb{Z}^d}N^{-\frac{d}{\alpha}}N^{\frac{2d}{\alpha}}\int_0^t\int_0^tp_{sN}^{d, \alpha}(\mathbf{0}, x)p_{vN}^{d, \alpha}(\mathbf{0}, x)dsdv\notag\\
&=N^{\frac{d}{\alpha}}\int_0^t\int_0^tp_{(s+v)N}^{d, \alpha}(\mathbf{0}, \mathbf{0})dsdv.
\end{align}
Since $\frac{1}{N^{\frac{1}{\alpha}}}X_{tN}^{d, \alpha}$ converges weakly to $\hat{X}_t^{d, \alpha}$, we have that $\frac{1}{t^{\frac{1}{\alpha}}}\hat{X}_t^{d, \alpha}$ and $\hat{X}_1^{d, \alpha}$ have the same distribution. Hence,
\[
f_t^{d, \alpha}(\mathbf{0})=t^{-\frac{d}{\alpha}}f_1^{d, \alpha}(\mathbf{0})
\]
and
\begin{align*}
\int_{\mathbb{R}^d}J_t^2(u)du&=\int_0^t\int_0^tf_{s+v}^{d, \alpha}(\mathbf{0})dsdv\\
&\leq f_1^{d, \alpha}(\mathbf{0})\int_0^{t}\left(\int_0^{+\infty}(s+v)^{-\frac{d}{\alpha}}dv\right)ds=O(1)\int_0^{t}s^{1-\frac{d}{\alpha}}ds.
\end{align*}
Since $\alpha>\frac{d}{2}$, $\int_{\mathbb{R}^d}J_t^2(u)du<+\infty$ and then \eqref{equ 4.8} follows from Lemma \ref{lemma LCLT} and \eqref{equ A.1}. \qed

\quad

\subsection{Proof of \eqref{equ 4.10}}\label{Appendix A.2}

Now we check \eqref{equ 4.10}.

\proof[Proof of \eqref{equ 4.10}]

We only deal with the case of $d=1$ and $\frac{1}{2}<\alpha<1$, since other cases follow from similar arguments.

By \eqref{equ 1.3 characteristic of X_t} and \eqref{equ 1.4 inverse formula random walk},
\begin{align}\label{equ A.2}
\int_0^{+\infty}p_s^{1, \alpha}(0, x)ds&=\frac{1}{2\pi}\int_{-\pi}^{+\pi}e^{-\sqrt{-1}x\theta}\left(\int_0^{+\infty}\psi_s^{d, \alpha}(\theta)ds\right)d\theta \notag\\
&=\frac{1}{2\pi}\int_{-\pi}^{+\pi}e^{-\sqrt{-1}x\theta}\frac{1}{\sum_{y\in \mathbb{Z}\setminus\{0\}}\frac{1-\cos(\theta y)}{|y|^{1+\alpha}}}d\theta\notag\\
&=\frac{1}{2\pi}\int_{-\pi}^{+\pi}\frac{\cos(x\theta)}{\sum_{y\in \mathbb{Z}\setminus\{0\}}\frac{1-\cos(\theta y)}{|y|^{1+\alpha}}}d\theta \notag\\
&=\frac{1}{2\pi|x|}\int_{-\pi|x|}^{\pi|x|}\frac{\cos\theta}{\sum_{y\in \mathbb{Z}\setminus\{0\}}\frac{1-\cos\left(\frac{\theta}{|x|} y\right)}{|y|^{1+\alpha}}}d\theta \notag\\
&=\frac{1}{2\pi|x|^{1-\alpha}}\int_{-\pi|x|}^{\pi|x|}\frac{\cos\theta}{G^{1, \alpha}\left(\frac{\theta}{|x|}\right)|\theta|^\alpha}d\theta,
\end{align}
where $G^{1, \alpha}$ is defined as in the proof of Lemma \ref{lemma LCLT}. As we have shown in the proof of Lemma \ref{lemma LCLT},
\[
\lim_{u\rightarrow 0}G^{1, \alpha}(u)=\tilde{G}^{1, \alpha}>0
\]
and hence
\begin{equation}\label{equ A.3}
\lim_{N\rightarrow+\infty}\int_{-M}^{M}\frac{\cos\theta}{G^{1, \alpha}\left(\frac{\theta}{|x|}\right)|\theta|^\alpha}d\theta=
\int_{-M}^{M}\frac{\cos\theta}{\tilde{G}^{1, \alpha}|\theta|^\alpha}d\theta
\end{equation}
for any $M>0$. We claim that
\begin{equation}\label{equ A.4}
\lim_{M\rightarrow+\infty}\limsup_{|x|\rightarrow+\infty}\left|\int_M^{\pi|x|}\frac{\cos\theta}{G^{1, \alpha}\left(\frac{\theta}{|x|}\right)|\theta|^\alpha}d\theta\right|=0.
\end{equation}
By \eqref{equ A.3} and \eqref{equ A.4}, \eqref{equ 4.10} holds with
\[
C_2=\frac{1}{2\pi}\int_{-\infty}^{\infty}\frac{\cos\theta}{\tilde{G}^{1, \alpha}|\theta|^\alpha}d\theta.
\]
At last, we check \eqref{equ A.4}. According to the formula of integration by parts,
\[
\int_M^{\pi|x|}\frac{\cos\theta}{G^{1, \alpha}\left(\frac{\theta}{|x|}\right)|\theta|^\alpha}d\theta
=\frac{\sin\theta }{G^{1, \alpha}\left(\frac{\theta}{|x|}\right)|\theta|^\alpha}\Bigg|_M^{\pi|x|}-\int_M^{\pi|x|}\sin\theta\frac{d}{d\theta}\left(\frac{\theta^{-\alpha}}{G^{1, \alpha}\left(\frac{\theta}{|x|}\right)}\right)d\theta.
\]
As we have shown in the proof of Lemma \ref{lemma LCLT}, $G^{1, \alpha}$ is continuous and strictly positive on $[0, \pi]$. Hence,
\begin{equation}\label{equ A.4 two}
0<\min_{0\leq u\leq \pi}G^{1, \alpha}(u)\leq \max_{0\leq u\leq \pi}G^{1, \alpha}(u)<+\infty.
\end{equation}
Hence, to check \eqref{equ A.4}, we only need to show that
\begin{equation}\label{equ A.5}
\lim_{M\rightarrow+\infty}\limsup_{|x|\rightarrow+\infty}\int_M^{\pi|x|}\left|\frac{d}{d\theta}\left(\frac{\theta^{-\alpha}}{G^{1, \alpha}\left(\frac{\theta}{|x|}\right)}\right)\right|d\theta=0.
\end{equation}
Note that
\[
\frac{d}{d\theta}\left(\frac{\theta^{-\alpha}}{G^{1, \alpha}\left(\frac{\theta}{|x|}\right)}\right)=\frac{-\alpha\theta^{-\alpha-1}}{G^{1, \alpha}\left(\frac{\theta}{|x|}\right)}+\frac{\theta^{-\alpha}\frac{d}{d\theta}G^{1, \alpha}\left(\frac{\theta}{|x|}\right)}{\left(G^{1, \alpha}\left(\frac{\theta}{|x|}\right)\right)^2}.
\]
According to \eqref{equ A.4 two} and the fact that $\int_1^{+\infty}\theta^{-\alpha-1}d\theta<+\infty$, to check \eqref{equ A.5}, we only need to show that
\begin{equation}\label{equ A.7}
\lim_{M\rightarrow+\infty}\limsup_{|x|\rightarrow+\infty}\int_M^{\pi|x|}\theta^{-\alpha}\left|\frac{d}{d\theta}G^{1, \alpha}\left(\frac{\theta}{|x|}\right)\right|d\theta=0.
\end{equation}
According to the definition of $G^{1, \alpha}$ given in the proof of Lemma \ref{lemma LCLT},
\[
\theta^{-\alpha}\frac{d}{d\theta}G^{1, \alpha}\left(\frac{\theta}{|x|}\right)=-\alpha\theta^{-\alpha-1}G^{1, \alpha}\left(\frac{\theta}{|x|}\right)
+\frac{2\theta^{-2\alpha}}{|x|^{1-\alpha}}\sum_{y=1}^{+\infty}\frac{\sin\left(\frac{\theta}{|x|}y\right)}{y^\alpha}.
\]
Hence, to check \eqref{equ A.7}, we only need to show that
\begin{equation}\label{equ A.8}
\lim_{M\rightarrow+\infty}\limsup_{|x|\rightarrow+\infty}\frac{1}{|x|^{1-\alpha}}\int_M^{\pi|x|}\theta^{-2\alpha}
\left|\sum_{y=1}^{+\infty}\frac{\sin\left(\frac{\theta}{|x|}y\right)}{y^\alpha}\right|d\theta=0.
\end{equation}
According to the fact that
\[
\lim_{u\rightarrow 0}\left|u\sum_{y=1}^{+\infty}\frac{\sin\left(uy\right)}{(uy)^\alpha}\right|=\left|\int_0^{+\infty}\frac{\sin x}{x^\alpha}dx\right|<+\infty,
\]
there exists $\epsilon>0$ such that
\[
\left|\sum_{y=1}^{+\infty}\frac{\sin\left(\frac{\theta}{|x|}y\right)}{y^\alpha}\right|\leq 2\left(\frac{\theta}{|x|}\right)^{\alpha-1}\left|\int_0^{+\infty}\frac{\sin x}{x^\alpha}dx\right|
\]
when $0<\theta<\epsilon |x|$. Hence,
\begin{align*}
\frac{1}{|x|^{1-\alpha}}\int_M^{\epsilon|x|}\theta^{-2\alpha}
\left|\sum_{y=1}^{+\infty}\frac{\sin\left(\frac{\theta}{|x|}y\right)}{y^\alpha}\right|d\theta
&\leq \frac{O(1)}{|x|^{1-\alpha}}\int_M^{+\infty}\theta^{-\alpha-1}d\theta|x|^{1-\alpha}\\
&=O(1)\int_M^{+\infty}\theta^{-\alpha-1}d\theta.
\end{align*}
Consequently, to check \eqref{equ A.8}, we only need to show that
\begin{equation}\label{equ A.9}
\lim_{|x|\rightarrow+\infty}\frac{1}{|x|^{1-\alpha}}\int_{\epsilon|x|}^{\pi|x|}\theta^{-2\alpha}
\left|\sum_{y=1}^{+\infty}\frac{\sin\left(\frac{\theta}{|x|}y\right)}{y^\alpha}\right|d\theta=0.
\end{equation}
For $\epsilon\leq u\leq \pi$ and integer $m\geq 1$, let $S_m(u)=\sum_{k=1}^m\sin(ku)$, then
\[
\sum_{k=1}^m\frac{\sin(ku)}{k^\alpha}=\sum_{k=1}^m\frac{S_k(u)-S_{k-1}(u)}{k^\alpha}=\sum_{k=1}^{m-1}S_k(u)\left(\frac{1}{k^\alpha}-\frac{1}{(k+1)^\alpha}\right)
+\frac{S_m(u)}{m^\alpha}.
\]
Note that
\begin{align*}
\sin\left(\frac{u}{2}\right)S_m(u)&=\sum_{k=1}^m\left(\cos\left(\left(k-\frac{1}{2}\right)u\right)-\cos\left(\left(k+\frac{1}{2}\right)u\right)\right)\\
&=\cos\left(\frac{u}{2}\right)-\cos\left(\left(m+\frac{1}{2}\right)u\right)
\end{align*}
and hence
\[
|S_m(u)|\leq \frac{2}{|\sin\left(\frac{u}{2}\right)|}.
\]
Consequently,
\[
\left|\sum_{k=1}^m\frac{\sin(ku)}{k^\alpha}\right|\leq \frac{2}{|\sin\left(\frac{u}{2}\right)|}\left(1-\frac{1}{m^\alpha}\right)+\frac{2}{|\sin\left(\frac{u}{2}\right)|m^\alpha}
\]
and hence
\[
\left|\sum_{k=1}^{+\infty}\frac{\sin(ku)}{k^\alpha}\right|\leq \frac{2}{C_3}
\]
for $\epsilon\leq u\leq \pi$, where $C_3=\min_{\epsilon\leq u\leq \pi}\sin\left(\frac{u}{2}\right)>0$. Therefore,
\begin{equation}\label{equ A.10}
\int_{\epsilon|x|}^{\pi|x|}\theta^{-2\alpha}
\left|\sum_{y=1}^{+\infty}\frac{\sin\left(\frac{\theta}{|x|}y\right)}{y^\alpha}\right|d\theta
\leq \frac{2}{C_3}\int_{\epsilon|x|}^{\pi|x|}\theta^{-2\alpha}d\theta.
\end{equation}
Since $2\alpha>1$ and $1-\alpha>0$, \eqref{equ A.9} follows from \eqref{equ A.10} and the proof is complete.  \qed

\quad

\subsection{Proof of \eqref{equ 4.16}}\label{Appendix A.3}

Now we check \eqref{equ 4.16}

\proof[Proof of \eqref{equ 4.16}]

According to the bilinear property of the covariance and the dominated convergence theorem, to check \eqref{equ 4.16}, we only need to show that
\begin{equation}\label{equ A.11}
\lim_{\|x\|_2\rightarrow +\infty}\left|{\rm Cov}_{\nu_p^{d, \alpha}}\left(\left(\eta_0(x+y_1)-\eta_0(x)\right)^2, \left(\eta_0(y_2)-\eta_0(\mathbf{0})\right)^2\right)\right|=0
\end{equation}
for any given $y_1, y_2\in \mathbb{Z}^d$.

According to the definition of $\nu_p^{d, \alpha}$,
\begin{align*}
&{\rm Cov}_{\nu_p^{d, \alpha}}\left(\left(\eta_0(x+y_1)-\eta_0(x)\right)^2, \left(\eta_0(y_2)-\eta_0(\mathbf{0})\right)^2\right)\\
&=\lim_{t\rightarrow+\infty}{\rm Cov}_{\mu_p}\left(\left(\eta_t(x+y_1)-\eta_t(x)\right)^2, \left(\eta_t(y_2)-\eta_t(\mathbf{0})\right)^2\right).
\end{align*}
By Proposition \ref{proposition voter dual},
\[
{\rm Cov}_{\mu_p}\left(\left(\eta_t(x+y_1)-\eta_t(x)\right)^2, \left(\eta_t(y_2)-\eta_t(\mathbf{0})\right)^2\right)=
{\rm \uppercase\expandafter{\romannumeral1}}-{\rm \uppercase\expandafter{\romannumeral2}},
\]
where
\[
{\rm \uppercase\expandafter{\romannumeral1}}
=\mathbb{E}\left(\hat{\mathbb{E}}_{\mu_p}
\left(\left(\eta_0\left(\mathcal{X}_{t, x+y_1}^{d, \alpha, 0}\right)-\eta_0\left(\mathcal{X}_{t, x}^{d, \alpha, 0}\right)\right)^2\left(\eta_0\left(\mathcal{X}_{t, y_2}^{d, \alpha, 0}\right)-\eta_0\left(\mathcal{X}_{t, \mathbf{0}}^{d, \alpha, 0}\right)\right)^2\right)\right)
\]
and
\begin{align*}
&{\rm \uppercase\expandafter{\romannumeral2}}\\
&=\mathbb{E}\left(\hat{\mathbb{E}}_{\mu_p}\left(\left(\eta_0\left(\mathcal{X}_{t, x+y_1}^{d, \alpha, 0}\right)-\eta_0\left(\mathcal{X}_{t, x}^{d, \alpha, 0}\right)\right)^2\right)\right)\mathbb{E}\left(\hat{\mathbb{E}}_{\mu_p}\left(\left(\eta_0\left(\mathcal{X}_{t, y_2}^{d, \alpha, 0}\right)-\eta_0\left(\mathcal{X}_{t, \mathbf{0}}^{d, \alpha, 0}\right)\right)^2\right)\right).
\end{align*}
We define $\left\{\left(\tilde{\mathcal{X}}_{t, y_2}^{d, \alpha, 0}, \tilde{\mathcal{X}}_{t, \mathbf{0}}^{d, \alpha, 0}\right)\right\}_{t\geq 0}$ as a copy of $\left\{\left(\mathcal{X}_{t, y_2}^{d, \alpha, 0}, \mathcal{X}_{t, \mathbf{0}}^{d, \alpha, 0}\right)\right\}_{t\geq 0}$ such that
\[
\left\{\left(\tilde{\mathcal{X}}_{t, y_2}^{d, \alpha, 0}, \tilde{\mathcal{X}}_{t, \mathbf{0}}^{d, \alpha, 0}\right)\right\}_{t\geq 0}
\]
is independent of $\left\{\left(\mathcal{X}_{t, x+y_1}^{d, \alpha, 0}, \mathcal{X}_{t, x}^{d, \alpha, 0}\right)\right\}_{t\geq 0}$ and
\[
\left(\tilde{\mathcal{X}}_{t, y_2}^{d, \alpha, 0}, \tilde{\mathcal{X}}_{t, \mathbf{0}}^{d, \alpha, 0}\right)=\left(\mathcal{X}_{t, y_2}^{d, \alpha, 0}, \mathcal{X}_{t, \mathbf{0}}^{d, \alpha, 0}\right)
\]
for $0\leq t\leq \tilde{\tau}_{x, y_1, y_2}$, where
\[
\tilde{\tau}_{x, y_1, y_2}=\inf\left\{t:~t\geq 0, \left\{\mathcal{X}_{t, x+y_1}^{d, \alpha, 0}, \mathcal{X}_{t, x}^{d, \alpha, 0}\right\}\bigcap\left\{\mathcal{X}_{t, y_2}^{d, \alpha, 0}, \mathcal{X}_{t, \mathbf{0}}^{d, \alpha, 0}\right\}\neq \emptyset\right\}.
\]
Then,
\begin{align*}
{\rm \uppercase\expandafter{\romannumeral2}}
=\mathbb{E}\left(\hat{\mathbb{E}}_{\mu_p}\left(\left(\eta_0\left(\mathcal{X}_{t, x+y_1}^{d, \alpha, 0}\right)-\eta_0\left(\mathcal{X}_{t, x}^{d, \alpha, 0}\right)\right)^2\right)\hat{\mathbb{E}}_{\mu_p}\left(\left(\eta_0\left(\tilde{\mathcal{X}}_{t, y_2}^{d, \alpha, 0}\right)-\eta_0\left(\tilde{\mathcal{X}}_{t, \mathbf{0}}^{d, \alpha, 0}\right)\right)^2\right)\right).
\end{align*}
On the event $\{\tilde{\tau}_{x, y_1, y_2}>t\}$, $\left\{\mathcal{X}_{t, x+y_1}^{d, \alpha, 0}, \mathcal{X}_{t, x}^{d, \alpha, 0}\right\}\bigcap\left\{\mathcal{X}_{t, y_2}^{d, \alpha, 0}, \mathcal{X}_{t, \mathbf{0}}^{d, \alpha, 0}\right\}=\emptyset$ and hence
\begin{align*}
&\hat{\mathbb{E}}_{\mu_p}
\left(\left(\eta_0\left(\mathcal{X}_{t, x+y_1}^{d, \alpha, 0}\right)-\eta_0\left(\mathcal{X}_{t, x}^{d, \alpha, 0}\right)\right)^2\left(\eta_0\left(\mathcal{X}_{t, y_2}^{d, \alpha, 0}\right)-\eta_0\left(\mathcal{X}_{t, \mathbf{0}}^{d, \alpha, 0}\right)\right)^2\right)\\
&=\hat{\mathbb{E}}_{\mu_p}
\left(\left(\eta_0\left(\mathcal{X}_{t, x+y_1}^{d, \alpha, 0}\right)-\eta_0\left(\mathcal{X}_{t, x}^{d, \alpha, 0}\right)\right)^2\right)\hat{\mathbb{E}}_{\mu_p}\left(\left(\eta_0\left(\mathcal{X}_{t, y_2}^{d, \alpha, 0}\right)-\eta_0\left(\mathcal{X}_{t, \mathbf{0}}^{d, \alpha, 0}\right)\right)^2\right)\\
&=\hat{\mathbb{E}}_{\mu_p}
\left(\left(\eta_0\left(\mathcal{X}_{t, x+y_1}^{d, \alpha, 0}\right)-\eta_0\left(\mathcal{X}_{t, x}^{d, \alpha, 0}\right)\right)^2\right)\hat{\mathbb{E}}_{\mu_p}\left(\left(\eta_0\left(\tilde{\mathcal{X}}_{t, y_2}^{d, \alpha, 0}\right)-\eta_0\left(\tilde{\mathcal{X}}_{t, \mathbf{0}}^{d, \alpha, 0}\right)\right)^2\right).
\end{align*}
Consequently,
\begin{equation*}
\left|{\rm \uppercase\expandafter{\romannumeral1}}-{\rm \uppercase\expandafter{\romannumeral2}}\right|
\leq 2\mathbb{P}\left(\tilde{\tau}_{x, y_1, y_2}\leq t\right)
\end{equation*}
and hence
\begin{equation}\label{equ A.12}
\left|{\rm Cov}_{\nu_p^{d, \alpha}}\left(\left(\eta_0(x+y_1)-\eta_0(x)\right)^2, \left(\eta_0(y_2)-\eta_0(\mathbf{0})\right)^2\right)\right|
\leq 2\mathbb{P}\left(\tilde{\tau}_{x, y_1, y_2}<\infty\right).
\end{equation}
According to the definition of $\tilde{\tau}_{x, y_1, y_2}$,
\begin{align*}
&\mathbb{P}\left(\tilde{\tau}_{x, y_1, y_2}<\infty\right)\\
&\leq \left(1-\Phi^{d, \alpha}(x+y_1-y_2)\right)+\left(1-\Phi^{d, \alpha}(x-y_2)\right)
+\left(1-\Phi^{d, \alpha}(x+y_1)\right)+\left(1-\Phi^{d, \alpha}(x)\right).
\end{align*}
For $(d, \alpha)$ satisfying \eqref{equ 4.4}, $\{X_t^{d, \alpha}\}_{t\geq 0}$ is transient and hence
\[
\lim_{\|x\|_2\rightarrow+\infty}\left(1-\Phi^{d, \alpha}(x)\right)=0.
\]
Consequently, \eqref{equ A.11} follows from \eqref{equ A.12} and the proof is complete.
\qed

\quad

\subsection{Proof of \eqref{equ 4.20}}\label{Appendix A.4}

Now we check \eqref{equ 4.20}.

\proof[Proof of \eqref{equ 4.20}] According to the Strong Markov property of the random walk,
\[
1-\Phi(y)=\mathcal{C}_{d, \alpha}\int_0^{+\infty}p_s^{d, \alpha}(\mathbf{0}, y)ds.
\]
Hence,
\[
\sum_{y\in \mathbb{Z}^d\setminus\{\mathbf{0}\}}\frac{1-\Phi^{d, \alpha}(y)}{\|y\|_2^{d+\alpha}}=\mathcal{C}_{d, \alpha}\int_0^{+\infty}\left(\sum_{y\in \mathbb{Z}^d\setminus\{\mathbf{0}\}}\frac{p_s^{d, \alpha}(\mathbf{0}, y)}{\|y\|_2^{d+\alpha}}\right)ds.
\]
By Kolmogorov-Chapman equation,
\[
\frac{d}{ds}p_s^{d, \alpha}(\mathbf{0}, \mathbf{0})=-p_s^{d, \alpha}(\mathbf{0}, \mathbf{0})\left(\sum_{y\in \mathbb{Z}^d\setminus\{\mathbf{0}\}}\frac{1}{\|y\|_2^{d+\alpha}}\right)+\sum_{y\in \mathbb{Z}^d\setminus\{\mathbf{0}\}}\frac{p_s^{d, \alpha}(\mathbf{0}, y)}{\|y\|_2^{d+\alpha}}
\]
and hence
\begin{align*}
&\int_0^{+\infty}\left(\sum_{y\in \mathbb{Z}^d\setminus\{\mathbf{0}\}}\frac{p_s^{d, \alpha}(\mathbf{0}, y)}{\|y\|_2^{d+\alpha}}\right)ds\\
&=\int_0^{+\infty}\frac{d}{ds}p_s^{d, \alpha}(\mathbf{0}, \mathbf{0})ds+\left(\sum_{y\in \mathbb{Z}^d\setminus\{\mathbf{0}\}}\frac{1}{\|y\|_2^{d+\alpha}}\right)
\int_0^{+\infty}p_s^{d, \alpha}(\mathbf{0}, \mathbf{0})ds\\
&=(0-1)+(\mathcal{C}_{d, \alpha})^{-1}\left(\sum_{y\in \mathbb{Z}^d\setminus\{\mathbf{0}\}}\frac{1}{\|y\|_2^{d+\alpha}}\right).
\end{align*}
Therefore,
\[
\sum_{y\in \mathbb{Z}^d\setminus\{\mathbf{0}\}}\frac{1-\Phi^{d, \alpha}(y)}{\|y\|_2^{d+\alpha}}=-\mathcal{C}_{d, \alpha}+\sum_{y\in \mathbb{Z}^d\setminus\{\mathbf{0}\}}\frac{1}{\|y\|_2^{d+\alpha}}
\]
and hence \eqref{equ 4.20} holds.
\qed

\quad

\subsection{Proof of \eqref{equ 4.32}}\label{Appendix A.5}

Now we check \eqref{equ 4.32}

\proof[Proof of \eqref{equ 4.32}]
By \eqref{equ 1.14},
\[
{\rm Cov}_{\nu_p^{d, \alpha}}\left(\eta_{\theta}(\mathbf{0}), \eta_0(\mathbf{0})\right)
=\hat{\mathbb{E}}_{\nu_p^{d, \alpha}}\left(\eta_{\theta}(\mathbf{0})\eta_0(\mathbf{0})\right)-p^2.
\]
By \eqref{equ voter dual}, for any $\eta\in \{0, 1\}^{\mathbb{Z}^d}$,
\[
\hat{\mathbb{E}}_\eta\eta_{\theta}(\mathbf{0})=\sum_{x\in \mathbb{Z}^d}p^{d, \alpha}_\theta(\mathbf{0}, x)\eta(x).
\]
Then, by the iterated expectation law,
\[
\hat{\mathbb{E}}_{\nu_p^{d, \alpha}}\left(\eta_{\theta}(\mathbf{0})\eta_0(\mathbf{0})\right)=\sum_{x\in \mathbb{Z}^d}p_\theta^{d, \alpha}(\mathbf{0}, x)
\hat{\mathbb{E}}_{\nu_p^{d, \alpha}}\left(\eta_0(x)\eta_0(\mathbf{0})\right).
\]
Hence, by \eqref{equ 1.14} and \eqref{equ 1.15},
\begin{align*}
{\rm Cov}_{\nu_p^{d, \alpha}}\left(\eta_{\theta}(\mathbf{0}), \eta_0(\mathbf{0})\right)&=
\sum_{x\in \mathbb{Z}^d}p_\theta^{d, \alpha}(\mathbf{0}, x)\left(\hat{\mathbb{E}}_{\nu_p^{d, \alpha}}\left(\eta_0(x)\eta_0(\mathbf{0})\right)-p^2\right)\\
&=\sum_{x\in \mathbb{Z}^d}p_\theta^{d, \alpha}(\mathbf{0}, x){\rm Cov}_{\nu_p^{d, \alpha}}\left(\eta_0(x), \eta_0(\mathbf{0})\right)\\
&=p(1-p)\sum_{x\in \mathbb{Z}^d}p_\theta^{d, \alpha}(\mathbf{0}, x)(1-\Phi^{d, \alpha}(x)).
\end{align*}
According to the strong Markov property of the random walk and Kolmogorov-Chapman equation,
\begin{align*}
\sum_{x\in \mathbb{Z}^d}p_\theta^{d, \alpha}(\mathbf{0}, x)(1-\Phi^{d, \alpha}(x))
&=\mathcal{C}_{d, \alpha}\sum_{x\in \mathbb{Z}^d}p_{\theta}^{d, \alpha}(\mathbf{0}, x)\int_0^{+\infty}p^{d, \alpha}_r\left(\mathbf{0}, x\right)dr\\
&=\mathcal{C}_{d, \alpha}\int_0^{+\infty}p^{d, \alpha}_{\theta+r}\left(\mathbf{0}, \mathbf{0}\right)dr.
\end{align*}
Hence, \eqref{equ 4.32} holds.
\qed

\subsection{Proof of \eqref{equ 5.8}}\label{Appendix A.6}
Now we check \eqref{equ 5.8}
\proof[Proof of \eqref{equ 5.8}]
We only need to show that
\begin{equation}\label{equ A.13}
\limsup_{s\rightarrow+\infty}\frac{\int_0^sdr\int_0^r d\theta\int_0^{+\infty}p_{\theta+v}^{3, 2}(\mathbf{0}, \mathbf{0})dv}{\frac{s^{\frac{3}{2}}}{\left(\log s\right)^{\frac{3}{2}}}}<+\infty.
\end{equation}
Let $C_6=\int_0^{+\infty}p_{v}^{3, 2}(\mathbf{0}, \mathbf{0})dv$ be defined as in Section \ref{section five}. Since
\[
\int_0^{1000}dr\int_0^r d\theta\int_0^{+\infty}p_{\theta+v}^{3, 2}(\mathbf{0}, \mathbf{0})dv\leq C_61000^2
\]
and
\[
\int_{1000}^sdr\int_0^{1000} d\theta\int_0^{+\infty}p_{\theta+v}^{3, 2}(\mathbf{0}, \mathbf{0})dv\leq C_61000(s-1000)
\]
for $s\geq 1000$, to check \eqref{equ A.13}, we only need to show that
\begin{equation}\label{equ A.14}
\limsup_{s\rightarrow+\infty}\frac{\int_{1000}^sdr\int_{1000}^r d\theta\int_0^{+\infty}p_{\theta+v}^{3, 2}(\mathbf{0}, \mathbf{0})dv}{\frac{s^{\frac{3}{2}}}{\left(\log s\right)^{\frac{3}{2}}}}<+\infty.
\end{equation}
By \eqref{equ 5.7} and the integration-by-parts formula, for $\theta\geq 1000$,
\begin{align*}
\int_0^{+\infty}p_{\theta+v}^{3, 2}(\mathbf{0}, \mathbf{0})dv
&\leq C_5\int_\theta^{+\infty}v^{-\frac{3}{2}}\left(\log v\right)^{-\frac{3}{2}}dv\\
&=C_5\left(2v^{-\frac{1}{2}}\left(\log v\right)^{-\frac{3}{2}}\Big|_{+\infty}^\theta
-3\int_\theta^{+\infty}v^{-\frac{3}{2}}(\log v)^{-\frac{5}{2}}dv\right)\\
&\leq 2C_5\theta^{-\frac{1}{2}}(\log \theta)^{-\frac{3}{2}}.
\end{align*}
For $r>1000$, by integration-by-parts formula,
\begin{align*}
\int_{1000}^r\theta^{-\frac{1}{2}}(\log \theta)^{-\frac{3}{2}} d\theta
&=2\theta^{\frac{1}{2}}(\log \theta)^{-\frac{3}{2}}\Big|_{1000}^r
-2\int_{1000}^r\theta^{\frac{1}{2}}\left(-\frac{3}{2}\right)(\log \theta)^{-\frac{5}{2}}\frac{1}{\theta}d\theta\\
&\leq 2r^{\frac{1}{2}}\left(\log r\right)^{-\frac{3}{2}}+\frac{3}{\log 1000}\int_{1000}^r\theta^{-\frac{1}{2}}(\log \theta)^{-\frac{3}{2}}d\theta
\end{align*}
and hence
\[
\int_{1000}^r\theta^{-\frac{1}{2}}(\log \theta)^{-\frac{3}{2}} d\theta\leq C_{15}r^{\frac{1}{2}}\left(\log r\right)^{-\frac{3}{2}},
\]
where $C_{15}=\frac{2}{1-\frac{3}{\log 1000}}$. For $s>1000$, according to a similar argument,
\[
\int_{1000}^sr^{\frac{1}{2}}\left(\log r\right)^{-\frac{3}{2}}dr\leq C_{16} s^{\frac{3}{2}}\left(\log s\right)^{-\frac{3}{2}},
\]
where $C_{16}=\frac{2}{3(1-\frac{1}{\log 1000})}$. In conclusion,
\[
\limsup_{s\rightarrow+\infty}\frac{\int_{1000}^sdr\int_{1000}^r d\theta\int_0^{+\infty}p_{\theta+v}^{3, 2}(\mathbf{0}, \mathbf{0})dv}{\frac{s^{\frac{3}{2}}}{\left(\log s\right)^{\frac{3}{2}}}}\leq 2C_5C_{15}C_{16}
\]
and hence \eqref{equ A.14} holds. Since \eqref{equ A.14} holds, the proof is complete.
\qed

\quad

\textbf{Acknowledgments.}
The author is grateful to Professor Zhao Linjie for useful comments about hydrodynamic limits of long-range interacting particle systems. The author is grateful to the financial
support from the National Natural Science Foundation of China with grant number 12371142.

{}

\begin{thebibliography}{}
\bibitem{Birkner2007}Birkner, M. and Z\"{a}hle, I. (2007). A functional CLT for the occupation time
of a state-dependent branching random walk. \emph{The Annals of Probability} \textbf{35}, 2063-2090.
\bibitem{Bramson1988}Bramson, M., Cox, J. T. and Griffeath, D. (1988). Occupation time large deviations of the voter
model. \emph{Probability Theory and Related Fields} \textbf{77}, 401-413.
\bibitem{Cox1983}Cox, J. T. and Griffeath, D. (1983). Occupation time limit theorems for the voter model. \emph{The Annals of Probability} \textbf{11}, 876-893.
\bibitem{Ethier1986}Ethier, N. and Kurtz, T. (1986). \emph{Markov Processes: Characterization and Convergence.} John Wiley and Sons, Hoboken, NJ, USA.
\bibitem{Kipnis1987}Kipnis, C. (1987). Fluctuations des temps d'occupation d'un site dans l'exclusion simple sym\'{e}trique. \emph{Annales de l'IHP Probabilit\'{e}s et statistiques} \textbf{23}(1), 21-35.
\bibitem{Lawler2010} Lawler, G., and Limic, V. (2010). \emph{Random Walk: a modern introduction.} Cambridge.
\bibitem{Lig1985}Liggett, T. M. (1985). \emph{Interacting Particle Systems.} Springer, New York.
\bibitem{Lig1999}Liggett, T. M. (1999). \emph{Stochastic interacting systems: contact, voter and exclusion processes.}
Springer, New York.
\bibitem{Maillard2009}Maillard, G, and Mountford, T. (2009). Large deviations for voter model occupation times in two dimensions. \emph{Annales de l'Institut Henri Poincar\'{e}-Probabilit\'{e}s et Statistiques} \textbf{45}, 577-588.
\bibitem{Presutti1983}Presutti, E. and Spohn, H. (1983). Hydrodynamics of the voter model. \emph{The Annals of Probability} \textbf{11}, 867-875.
\bibitem{Xue2024}Xue, X. (2024). Sample path central limit theorem for the occupation time of the voter model on a lattice. Arxiv: 2412.05064.
\end{thebibliography}
\end{document}